\documentclass{amsart}
\usepackage{bbm,amsmath,amsfonts,amssymb}

\newtheorem{theorem}{Theorem}[section]
\newtheorem{lemma}[theorem]{Lemma}
\newtheorem{corollary}[theorem]{Corollary}
\newtheorem{proposition}[theorem]{Proposition}
\newtheorem{definition}[theorem]{Definition}

\newtheorem{remark}[theorem]{Remark}

\newcounter{figures}[section]

\def\bZ{{\mathbb Z}}

\def\bN{{\mathbb N}}

\def\bC{{\mathbb C}}
\def\bR{{\mathbb R}}

\def\bT{{\mathbb T}}
\def\bS{{\mathbb S}}

\def\cD{\mathcal{D}}
\def\cE{\mathcal{E}}

\def\cH{\mathcal{H}}

\def\cL{\mathcal{L}}
\def\cM{\mathcal{M}}
\def\cN{\mathcal{N}}

\def\cP{\mathcal{P}}

\def\cS{\mathcal{S}}

\def\cX{\mathcal{X}}
\def\cY{\mathcal{Y}}

\def\LL{L}
\def\supp{{\rm{supp}\, }}
\def\eps{{\varepsilon}}
\def\Id{{\rm Id}}
\def\ONE{{\mathbbm 1}}
\def\tONE{{\tilde{\mathbbm 1}}}
\def\diam{{\rm diam \,}}
\def\ph{\varphi}

\def\Css{{K(\sigma)}}
\def\CPhi{{K(\sigma_*)}}
\def\Cdiam{{c^\diamond}}

\def\ww{w}
\def\spec{{Spec \,}}
\def\R{\mathbb{R}}
\def\UCB{{\rm UCB}}

\def\tB{\tilde{B}}
\def\tb{\tilde{b}}
\def\tF{\tilde{F}}
\def\tf{\tilde{f}}
\def\W{w_{\alpha, \beta}}
\def\gg{\gamma}

\def\XX{{\cX}}
\def\sig{\sigma}
\def\dd{{\delta}}
\def\ddj{{\dd_j}}

\def\LLam{\Gamma}
\def\TTheta{\Theta}
\def\TT{T}
\def\SSS{S}

\def\lam{\lambda}
\def\cd{{c_\diamond}}
\def\cf{{c_\flat}}
\def\bb{\beta}
\def\bbb{\beta}
\def\cn{{c_\natural}}
\def\cs{{c_\sharp}}
\def\ct{\tilde{c}}
\def\tc{{c_*}}
\def\dst{{\delta_\star}}
\def\kap{\kappa}
\def\hh{F}
\def\kaph{{\hat\kappa}}
\def\ch{{\hat c}}
\def\bet{\zeta}
\def\kk{a}
\def\FF{F}
\def\Cstar{{C^\star}}
\def\cstar{{c^\star}}
\def\cast{{c_*}}
\def\cdg{{\hat{c}}}

\begin{document}

\title[Heat kernel based decomposition of spaces]
{Heat kernel based decomposition of spaces \\of distributions in the framework \\of Dirichlet spaces}

\author{Gerard Kerkyacharian}
\address{
Laboratoire de Probabilit\'{e}s et Mod\`{e}les Al\'{e}atoires, CNRS-UMR 7599,
Universit\'{e}e Paris Diderot, Batiment Sophie Germain, 
Avenue de France, Paris 75013, France}
\email{kerk@math.univ-paris-diderot.fr}

\author{Pencho Petrushev}
\address{Department of Mathematics\\University of South Carolina\\
Columbia, SC 29208}
\email{pencho@math.sc.edu}




\thanks{P. Petrushev has been supported by NSF Grant DMS-1211528.}

\begin{abstract}
Classical and nonclassical Besov and Triebel-Lizorkin spaces with complete range of indices
are developed in the general setting of Dirichlet space with a doubling measure and
local scale-invariant Poincar\'e inequality.
This leads to Heat kernel with small time Gaussian bounds and H\"{o}lder continuity,
which play a central role in this article.
Frames with band limited elements of sub-exponential space localization are developed,
and frame and heat kernel characterizations of Besov and Triebel-Lizorkin spaces are established.
This theory, in particular, allows to develop Besov and Triebel-Lizorkin spaces
and their frame and heat kernel characterization
in the context of Lie groups,
Riemannian manifolds,
and other settings.
\end{abstract}

\maketitle

\tableofcontents

\section{Introduction}\label{introduction}
\setcounter{equation}{0}

Spaces of functions or distributions play a prominent role in various areas of mathematics
such as harmonic analysis, PDEs, approximation theory, probability theory and statistics
and their applications.
The main purpose of this article is to develop the theory of Besov and Triebel-Lizorkin spaces
with full set of indices in the general setting
of strictly local regular Dirichlet spaces with doubling measure and local scale-invariant
Poincar\'e inequality, leading to a markovian heat kernel with small time
Gaussian bounds and H\"{o}lder continuity.
The gist of our method is to have the freedom of dealing with different geometries,
on compact and noncompact sets,
and with nontrivial weights,
and at the same time to allow for the development and frame decomposition of
Besov and Triebel-Lizorkin spaces with complete range of indices,
and therefore to cover a great deal of classical and nonclassical settings.
As an application, our theory allows to develop in full Besov and Triebel-Lizorkin spaces
and their frame decomposition in the setup of Lie groups or homogeneous spaces with polynomial volume growth,
complete Riemannian manifolds with Ricci curvature bounded from below and satisfying
the volume doubling condition, and various other nonclassical setups.

There are many forerunners of the ideas in this article which we even do not try to list here.
Our development can be viewed as a generalization of the Littlewood-Paley theory developed by
Frazier and Jawerth in the classical setting on $\R^n$ in \cite{F-J1, F-J2}, see also \cite{F-J-W}.
More recently, Besov and Triebel-Lizorkin spaces and their frame characterization
were developed in nonclassical settings such as on the sphere \cite{NPW2}
and more general homogeneous spaces \cite{PesG},
on the interval with Jacobi weights \cite{KPX1},
on the ball with weights \cite{KPX2}, and in the context of Hermite \cite{PX3}
and Laguerre expansions \cite{KPPX}.

This is a follow-up paper to \cite{CKP}, where we laid down some of the ground work needed for
the developments in this paper.
We adhere to the framework and notation established in \cite{CKP}, which we recall
in the following, beginning with {\bf the setting}:

I. We assume that $(M,\rho,\mu)$ is a metric measure space satisfying the conditions:
$(M, \rho)$ is a locally compact  metric space with distance $\rho(\cdot, \cdot)$
and $\mu$ is a~positive Radon measure
such that the following {\em volume doubling condition} is valid
\begin{equation}\label{doubling-0}
0 < \mu(B(x,2r) ) \le c_0\mu(B(x,r))<\infty
\quad\hbox{for all $x \in M$ and $r>0$,}
\end{equation}
where $B(x,r)$ is the open ball centered at $x$ of radius $r$ and $c_0>1$ is a constant.
Note that (\ref{doubling-0}) readily implies
\begin{equation}\label{doubling}
\mu(B(x,\lambda r) ) \le c_0\lambda^d \mu(B(x,r))
\quad\hbox{for $x \in M$, $r>0$, and $\lambda >1$.}
\end{equation}
Here $d=\log_2 c_0 >0$ is a constant playing the role of a dimension, but one should not confuse it with dimension.


II. The main assumption is that the local geometry of the space $(M,\rho,\mu)$ is related to
a self-adjoint positive operator $L$ on $L^2(M, d\mu)$,
mapping real-valued to real-valued functions,
such that the associated semigroup $P_t=e^{-tL}$ consists of integral operators with
(heat) kernel $p_t(x,y)$ obeying the conditions:

$\bullet$ {\it Small time  Gaussian upper bound}:
\begin{equation}\label{Gauss-local}
|p_t(x,y)| \le \frac{ \Cstar\exp \{-\frac{\cstar\rho^2(x,y)}t\}}{\sqrt{\mu(B(x,\sqrt t))\mu(B(y, \sqrt t))}}
\quad\hbox{for} \;\;x,y\in M,\,0<t\le 1.
\end{equation}

$\bullet$ {\em H\"{o}lder continuity}: There exists a constant $\alpha>0$ such that
\begin{equation}\label{lip}
\big|  p_t(x,y) - p_t(x,y')  \big|
\le \Cstar\Big(\frac{\rho(y,y')}{\sqrt t}\Big)^\alpha
\frac{\exp\{-\frac{\cstar\rho^2(x,y)}t \}}{\sqrt{\mu(B(x,\sqrt t))\mu(B(y, \sqrt t))}}
\end{equation}
for $x, y, y'\in M$ and $0<t\le 1$, whenever $\rho(y,y')\le \sqrt{t}$.

\smallskip

$\bullet$ {\it Markov property}:
\begin{equation}\label{hol3}
\int_M p_t(x,y) d\mu(y) \equiv 1
\quad\hbox{for $t >0$.}
\end{equation}
Above $\Cstar, \cstar>0$ are structural constants which along with $c_0$ will affect most of the
constants in the sequel.

In certain situations, we shall assume one or both of the following
{\em additional conditions}:

$\bullet$ {\em Reverse doubling condition}: There exists a constant
$c>1$ such that
\begin{equation}\label{reverse-doubling}
\mu(B(x,2r) ) \ge c \mu(B(x,r))
\quad\hbox{for $x \in M$ and $0< r \le \frac{\diam M}{3}$.}
\end{equation}

$\bullet$ {\em Non-collapsing condition}:
There exists a~constant $c>0$ such that
\begin{equation}\label{non-collapsing}
\inf_{x\in M}\mu(B(x,1) )\ge c.
\end{equation}
It will be explicitly indicated where each of these two conditions is required.

As is shown in \cite{CKP} a natural {\bf realization of the above setting} appears
in the general framework of Dirichlet spaces.
It turns out that in the setting of strictly local regular Dirichlet spaces with a complete
intrinsic metric (see \cite{BH,FUKU,Ouhabaz,ALB,Sturm0, Sturm1,Sturm, BM, BM1, Davies})
it suffices to only verify the local Poincar\'e inequality and the global doubling condition
on the measure and then our general theory applies in full.
We refer the reader to \S1.2 in \cite{CKP} for the details.

The point is that situations where our theory applies are quite common, which becomes evident from
{\bf the examples} given in \cite{CKP}.
We next describe them briefly.

$\bullet$ {\em Uniformly elliptic divergence form operators on $\R^d$.}
Given a uniformly elliptic symmetric matrix-valued function $\{a_{i,j}(x)\}$ depending on $x\in\R^d$,
one can define an operator
$
L=-\sum_{i,j=1}^d\frac{\partial}{\partial x_i}\left(a_{i,j}\frac{\partial}{\partial x_j}\right)
$
on $L^2(\R^d,dx)$  via the associated quadratic form.
The uniform ellipticity condition yields that the intrinsic metric associated with this operator
is equivalent to the Euclidean distance.
The Gaussian upper and lower bounds on the heat kernel in this setting
are due to Aronson and the H\"older regularity of the solutions is due to Nash \cite{N}.

$\bullet$ {\em Domains in $\R^d$.}
Uniformly elliptic divergence form operators on domains in $\R^d$ can be developed by choosing boundary conditions.
In this case the upper bounds of the heat kernels are well understood (see e.g. \cite{Ouhabaz}).
The Gaussian lower bounds is much more complicated to establish and one has to choose Neumann conditions
and impose regularity assumptions on the domain. We refer the reader to \cite{GS} for more details.

$\bullet$ {\em Riemannian manifolds and Lie groups.}
The local Poincar\'e inequality and doubling condition are verified for
the Laplace-Beltrami operator of a Riemannian manifold with non-negative Ricci curvature \cite{LY},
also for manifolds with Ricci curvature bounded from below if one assumes in addition that they satisfy
the volume doubling property,
also for manifolds that are quasi-isometric to such a~manifold \cite{G1, SD, Parma},
also for co-compact covering manifolds whose deck transformation group has polynomial growth \cite{SD, Parma},
for sublaplacians on polynomial growth Lie groups \cite{VSC, Robinson} and their homogeneous spaces \cite{M}.
Observe that the case of the sphere endowed with the natural Laplace-Beltrami operator
treated in \cite{NPW1,NPW2} and the case of more general compact homogeneous spaces endowed with
the Casimir operator considered in \cite{PesG} fall into the above category.
One can also consider variable coefficients operators on Lie groups, see \cite{SStr}.

We refer the reader to \cite[Section 2.1]{GS} for more details on the above examples
and to \cite{Davies, Grig, Ouhabaz, S, VSC} as general references for the heat kernel.


$\bullet$ {\em Heat kernel on $[-1, 1]$ generated by the Jacobi operator.}
In this case
$M =[-1, 1]$ with $d\mu(x) = \W(x)dx$, where
$\W(x)=(1-x)^\alpha (1+x)^\beta$, $\alpha , \beta >-1$,
is the classical Jacobi weight, and $L$ is the Jacobi operator.
As is well-known, e.g. \cite{Sz},
$LP_k =\lambda_k P_k$, where $P_k$ ($k\ge 0$) is the $k$th degree (normalized) Jacobi polynomial and
$\lambda_k = k(k+\alpha + \beta +1)$.
As is shown in \cite{CKP} in this case the general theory applies,
resulting in a complete strictly local Dirichlet space with an intrinsic metric
$
\rho(x, y)= |\arccos x - \arccos y|.
$
It is also shown that the respective scale-invariant Poincar\'e inequality is valid
and the measure $\mu$ obeys the doubling condition.
Therefore, this example fits in the general setting described above and our theory applies
and covers completely the results
in \cite{KPX1, PX1}.


The development of weighted 
spaces on the unit ball in $\R^d$
in \cite{KPX2, PX2} also fits in our general setting.
The treatment of this and other examples will be the theme of a future work.


In this article we advance on several fronts.
We refine considerably one of the main results in \cite{CKP}
which asserts that in the general setting described above for any compactly supported function $f\in C^\infty(\bR)$
obeying $f^{(2\nu+1)}(0)=0$, $\nu\ge 0$, the operator $f(\sqrt{L})$ has a kernel $f(\sqrt{L})(x, y)$ of
nearly exponential space localization (see Theorem~\ref{thm:main-local-kernels} below).
Furthermore, we show that for appropriately selected functions $f$ of this sort with ``small" derivatives
$f(\sqrt{L})(x, y)$ has sub-exponential space localization:
$$
|f(\delta \sqrt{L})(x, y)|
\le \frac{c\exp\big\{-\kap\big(\frac{\rho(x, y)}{\delta}\big)^{1-\eps}\big\}}
{\sqrt{\mu(B(x, \delta))\mu(B(y, \delta))}}.
$$
We also show that the class of integral operators with sub-exponentially localized kernels is an algebra,
which plays a crucial role in the development of frames.
We~make a substantial improvement in the scheme for constriction of frames from \cite{CKP}
which enable us to construct duals of sub-exponential space localization.
These advances allow us to generalize in full the theory of Frazier-Jawerth \cite{F-J1, F-J2, F-J-W}.


To introduce Besov spaces in the general setting of this article we follow the well known idea
\cite{Peetre, Triebel1, Triebel2}
of using spectral decompositions induced by a self-adjoint positive operator.
Consider real-valued functions $\varphi_0, \varphi\in C^\infty(\bR_+)$ such that
$\supp \varphi_0 \subset   [0, 2]$, $\varphi_0^{(2\nu+1)}(0) = 0$ for $\nu\ge 0$,
$\supp \varphi \subset   [1/2, 2]$, and
$|\varphi_0(\lambda)| +\sum_{j\ge 1} |\varphi(2^{-j}\lambda)| \ge c >0$ on~$\bR_+$.
Set $\varphi_j(\lambda):= \varphi(2^{-j}\lambda)$ for $j\ge 1$.
The possibly anisotropic geometry of $M$ is the reason for introducing two types of Besov spaces (\S\ref{besov-spaces}):

(i) The ``classical" Besov space  $B_{pq}^{s}=B_{pq}^{s}(L)$ is defined as the set of all
distributions $f$ such that
$$
\|f\|_{B_{pq}^{s}} :=
\Big(\sum_{j\ge 0} \Big(2^{s j}
\|\varphi_j(\sqrt{L}) f(\cdot)\|_{\LL^p}
\Big)^q\Big)^{1/q} <\infty, \quad\hbox{and}
$$

(ii) The ``nonclassical" Besov space $\tB_{pq}^{s}=\tB_{pq}^{s}(L)$ is defined by the norm
$$
\|f\|_{\tB_{pq}^{s}} :=
\Big(\sum_{j\ge 0} \Big(
\| |B(\cdot, 2^{-j})|^{-s/d}
\varphi_j(\sqrt{L}) f(\cdot)\|_{\LL^p}
\Big)^q\Big)^{1/q}.
$$
Our main motivation for introducing the spaces $\tB_{pq}^{s}$ lies in nonlinear approximation (\S\ref{Nonlin-app}).
However, we believe that these spaces capture well the geometry of the underlying space $M$
and will play an important role in other situations.

``Classical" Triebel-Lizorkin spaces $F_{pq}^{s}= F_{pq}^{s}(L)$
are defined by means of the norms
$$
\|f\|_{F_{pq}^s} :=
\Big\|\Big(\sum_{j\ge 0} \Big(
2^{js}
|\varphi_j(\sqrt L) f(\cdot)|
\Big)^q\Big)^{1/q}\Big\|_{\LL^p},
$$
while their ``nonclassical" version  $\tF_{pq}^{s}= \tF_{pq}^{s}(L)$
is introduced through the norms
$$
\|f\|_{\tF_{pq}^s} :=
\Big\|\Big(\sum_{j\ge 0} \Big(
|B(\cdot, 2^{-j})|^{-s/d}
|\varphi_j(\sqrt L) f(\cdot)|
\Big)^q\Big)^{1/q}\Big\|_{\LL^p}.
$$
It is important that our setting, though general, permits to develop Besov and Triebel-Lizorkin spaces
with complete range of indices, e.g. $s\in\bR$, $0<p, q \le \infty$, in the case of Besov spaces.
We only consider inhomogeneous Besov and Triebel-Lizorkin spaces here for this enables us to
treat simultaneously the compact and noncompact cases.
Their homogeneous version, however, can be developed in a similar manner.

One of the main results in this article is the frame decomposition of the Besov and Triebel-Lizorkin spaces
in the spirit of the $\varphi$-transform decomposition in the classical case by Frazier and Jawerth~\cite{F-J1, F-J2}.
To~show the flavor of these results, let
$\{\psi_\xi\}_{\xi\in\cX}$ and $\{\tilde\psi_\xi\}_{\xi\in \cX}$
be the pair of dual frames constructed here, indexed by a multilevel set $\cX=\cup_{j\ge 0}\cX_j$.
Then the decomposition of e.g. $\tB_{pq}^{s}$ takes the form~(\S\ref{sec:frame-dec-B-spaces})
$$
\|f\|_{\tB_{pq}^s} \sim \Big(\sum_{j\ge 0} \Bigl[\sum_{\xi\in \cX_j}
\Big(|B(\xi, b^{-j})|^{-s/d}
\|\langle f,\tilde\psi_{\xi}\rangle\psi_{\xi}\|_{p}\Big)^p\Bigr]^{q/p}\Bigr)^{1/q}.
$$
We also establish characterization of Besov and Triebel-Lizorkin spaces in terms of the heat kernel.
For instance, for $\tB_{pq}^{s}$ we have for $m>s$ (\S\ref{sec:heat-dec-B-spaces})
$$
\|f\|_{\tB_{pq}^{s}}
\sim \||B(\cdot, 1)|^{-s/d}e^{-L}f\|_p
+ \Big(\int_0^1\big\||B(\cdot, t^{1/2})|^{-s/d}(tL)^{m/2} e^{-tL}f\big\|_p^q \frac{dt}{t}
\Big)^{1/q}.
$$

As will be shown our theory covers completely the classical case on $\bR^d$ and on the torus $\bT^d$
as well as the above mentioned cases on the sphere \cite{NPW2} and more general homogeneous spaces \cite{PesG},
on the interval \cite{KPX1}, and on the ball \cite{KPX2}.
Our theory also applies in full in the various situations briefly indicated above.
Others are yet to be identified or developed.
Related interesting issues such as atomic decompositions and interpolation will not be treated here.

The metric measure space $(M, \rho, \mu)$ (with the doubling condition) from the setting of this article
is a space of homogeneous type in the sense of Coifman and Weiss \cite{CW}.
The theory of Besov and Triebel-Lizorkin spaces on general homogeneous spaces is well developed by now,
see e.g. \cite{HMY1, HMY2, MY, YZ}. The principle difference between this theory and our theory is
that the smoothness of the spaces in the former theory is limited ($|s|<\eps$).
Yet, it is a reasonable question to explore the relationship between these two theories.
We do not attempt to address this issue here.

For Hardy spaces $H^p$ associated with non-negative self-adjoint operators
under the general assumption of the Davies-Gaffney estimate we refer the reader to \cite{DY, JY, HLMMY}.

\smallskip

\noindent
{\bf The organization} of this paper is as follows:
In \S2 we present some technical results and background.
In \S3 we refine and extend the functional calculus results from \cite{CKP}.
In~\S4 we develop an improved version of the construction of frames from \cite{CKP}
which produces frame elements of sub-exponential space localization.
In \S 5 we introduce distributions in the setting of this paper and establish some of
their main properties and decomposition.
In \S 6 we introduce classical and nonclassical inhomogeneous Besov spaces and give their characterization
in terms of the heat kernel and the frames from \S 4.
We also show the application of Besov spaces to nonlinear approximation from frames.
In \S 7 we develop classical and nonclassical inhomogeneous Triebel-Lizorkin spaces in the underlying setting
and establish their characterization in terms of the heat kernel and the frames from \S 4.
We also present identification of some Triebel-Lizorkin spaces.

\smallskip

\noindent
{\bf Notation.}
Throughout this article we shall use the notation
$|E|:= \mu(E)$
and $\ONE_E$ will denote the characteristic function of  $E\subset M$,
$\|\cdot\|_p=\|\cdot\|_{\LL^p}:=\|\cdot\|_{L^p(M, d\mu)}$.
$\UCB$ will stand for the space of all uniformly continuous and bounded functions on~$M$.
We shall denote by $C^\infty_0(\bR_+)$ the set of all compactly supported $C^\infty$ functions
on $\bR_+:=[0, \infty)$.
In some cases ``$\sup$" will mean ``${\rm ess \,sup}$", which will be clear from the contex.
Positive constants will be denoted by $c$, $C$, $c_1$, $c'$, $\dots$ and they may vary
at every occurrence. Most of them will depend on the basic structural constants $c_0, \Cstar, \cstar$
from (\ref{doubling-0})-(\ref{lip}). This dependence usually will not be indicated explicitly.
Some important constants will be denoted by $\cn, \cs, \cf, \dots$
and they will remain unchanged throughout.
The notation $a\sim b$ will stand for $c_1\le a/b\le c_2$.

\section{Background}\label{sec:background}
\setcounter{equation}{0}

In this section we collect a number of technical results that will be needed in the sequel.
Most of the nontrivial of them are proved in \cite{CKP}.

\subsection{Estimates and facts related to the doubling and other conditions}

Since $B(x,r) \subset B(y, \rho(y,x) +r)$, then (\ref{doubling}) yields
\begin{equation}\label{D2}
|B(x, r)| \le c_0\Big(1+ \frac{\rho(x,y)}{r}\Big)^d  |B(y, r)|,
\quad x, y\in M, \; r>0.
\end{equation}

The reverse doubling condition (\ref{reverse-doubling}) implies
\begin{equation}\label{D3}
|B(x,\lambda r)| \ge c^{-1}\lambda^\bet |B(x, r)|,
\quad \lambda > 1, \; r>0,\; \hbox{$0< \lambda r<\frac{\diam M}{3}$,}
\end{equation}
where $c>1$ is the constant from (\ref{reverse-doubling}) and $\bet=\log_2 c >0$.
In this article, the reverse doubling condition will be used:
for lower bound estimates of the $L^p$ norms of operator kernels (\S\ref{sec:kernel-norms})
and frame elements (\S\ref{sec:frames}),
in nonlinear approximation from frames (\S \ref{Nonlin-app}),
and in the identification of some Triebel-Lizorkin spaces (\S\ref{sec:identify-Fspaces}).

The non-collapsing condition (\ref{non-collapsing}) and (\ref{doubling}) yield
\begin{equation}\label{est-muB-unify1}
\inf_{x\in M}|B(x, r)| \ge \check{c} r^d,
\quad 0< r \le 1,\quad \check{c} = const.
\end{equation}
The non-collapsing condition is needed in establishing the $L^p \to L^q$ boundedness
of integral operators (\S\ref{max-operator})
and for embedding results for Besov spaces (\S\ref{sec:B-embedding}).

As shown in \cite{CKP} the following clarifying statements hold:

\smallskip

\noindent
(a) $\mu(M) <\infty$ if and only if $\diam M <\infty$.
Moreover, if $\diam M =D <\infty$, then
\begin{equation}\label{est-muB}
\inf_{x\in M}|B(x, r)| \ge cr^d|M|D^{-d},
\quad  0< r \leq D.
\end{equation}

\smallskip

\noindent
(b)
If $M$ is connected, then the reverse doubling condition (\ref{reverse-doubling}) is valid.
Therefore, it is not quite restrictive.

\smallskip

\noindent
(c) In general, $|B(x, r)|$ can be much larger than $O(r^d)$ for certain points $x\in M$
as is evident from the example on $[-1, 1]$ with the heat kernel induced by the Jacobi operator.

\smallskip

The following symmetric functions will govern the localization of most operator kernels
in the sequel:
\begin{equation}\label{def-Dxy}
D_{\delta, \sigma}(x,y)
:= \big(|B(x, \delta)||B(y, \delta)|\big)^{-1/2}
\Big(1+ \frac{\rho(x,y)}{\delta}\Big)^{-\sigma},
\quad x,y \in M.
\end{equation}
Observe that (\ref{doubling}) and (\ref{D2}) readily imply
\begin{equation}\label{E1}
 D_{\delta, \sigma}(x,y)
 \le c|B(x, \delta)|^{-1}\Big(1+ \frac{\rho(x,y)}{\delta}\Big)^{-\sigma+d/2}.
\end{equation}
Furthermore, for $0<p<\infty$ and $\sigma>d(1/2+1/p)$
\begin{equation}\label{INT}
\|D_{\delta,\sigma}(x,\cdot)\|_p
=\Big(\int_M \big[D_{\delta,\sigma}(x,y)\big]^p d\mu(y)\Big)^{1/p}
\le c_p |B(x,\delta)|^{1/p-1}.
\end{equation}
and
\begin{equation}\label{Comp}
\int_M D_{\delta,\sigma}(x, u) D_{\delta,\sigma}(u,y) d\mu(u)
\le cD_{\delta,\sigma}( x, y)\quad \hbox{if } \; \sigma >  2d.
\end{equation}

The above two estimates follow readily by the following lemma
which will also be useful.


\begin{lemma}\label{lem:tech-est}
$(a)$ For $\sigma > d$ and $\delta >0$
\begin{equation}\label{tech1}
\int_M(1+\delta^{-1}\rho(x, y))^{-\sigma}d\mu(y) \le
c_1|B(x, \delta)|,
\quad x\in M.
\end{equation}
%

$(b)$ If $\sigma>d$, then for $x, y\in M$ and $\delta>0$
\begin{align}\label{tech2}
\int_M\frac{1}{(1+\delta^{-1}\rho(x, u))^{\sigma}(1+\delta^{-1}\rho(y, u))^{\sigma}} d\mu(u)
&\le
 2^\sigma c_1 \frac{|B(x, \delta)|+|B(y, \delta)|}{(1+\delta^{-1}\rho(x, y))^{\sigma}}\notag\\
&\le
c_2\frac{|B(x, \delta)|}{(1+\delta^{-1}\rho(x, y))^{\sigma-d}}.
\end{align}

$(c)$ If $\sigma  > 2d $, then for $x, y\in M$ and $\delta>0$
\begin{align}\label{tech3}
\int_M\frac{1}{|B(u, \delta)|(1+\delta^{-1}\rho(x, u))^{\sigma}(1+\delta^{-1}\rho(y, u))^{\sigma}} d\mu(y)
\le
\frac{ c_3}{(1+\delta^{-1}\rho(x, y))^{\sigma}},
\end{align}
and if in addition $0<\delta \le 1$, then
\begin{align}\label{tech5}
\int_M\frac{1}{|B(u, \delta)|(1+\delta^{-1}\rho(x, u))^{\sigma}(1+\rho(y, u))^{\sigma}} d\mu(y)
\le \frac{c_4}{(1+\rho(x, y))^{\sigma}}.
\end{align}

\end{lemma}

\noindent
{\bf Proof.}
Estimates (\ref{tech1})-(\ref{tech3}) are proved in \cite{CKP} (see Lemma~2.3).
Estimate (\ref{tech5}) follows easily from (\ref{tech3}).
Indeed, denote by $J$ the integral in (\ref{tech5}).
If $\rho(x, y) \le 1$, then using (\ref{D2}), (\ref{tech1}), and that $\sigma>2d$ we get
\begin{align*}
J \le\int_M\frac{d\mu(y)}{|B(u, \delta)|(1+\delta^{-1}\rho(x, u))^{\sigma}}
\le c_0\int_M\frac{d\mu(y)}{|B(x, \delta)|(1+\delta^{-1}\rho(x, u))^{\sigma-d}}
\le c
\end{align*}
which implies (\ref{tech5}).
If $\rho(x, y) > 1$, then using (\ref{tech3})
\begin{align*}
J &\le\int_M\frac{\delta^{-\sigma}d\mu(y)}
{|B(u, \delta)|(1+\delta^{-1}\rho(x, u))^{\sigma}(1+\delta^{-1}\rho(y, u))^{\sigma}}
\le \frac{c\delta^{-\sigma}}{(1+\delta^{-1}\rho(x, y))^{\sigma}},
\end{align*}
which yields (\ref{tech5}).
$\qed$

\subsection{Maximal \boldmath $\delta$-nets}\label{sec:max-d-nets}

In the construction of frames in the general setting of this article there is an underlying
sequence of maximal $\delta$-nets $\{\cX_j\}_{j\ge 0}$ on $M$:
{\em We say that $\cX\subset M$ is a $\delta$-net on $M$ $(\delta>0)$ if $\rho(\xi, \eta) \ge \delta$
for all $\xi, \eta\in\cX$,
and $\cX\subset M$ is a maximal $\delta$-net on $M$
if $\cX$ is a $\delta$-net on $M$ that cannot be enlarged.
}

We next summarize the basic properties of maximal $\delta$-nets \cite[Proposition 2.5]{CKP}:
{\em A maximal $\delta$-net on $M$ always exists and
if $\cX$ is a maximal $\delta$-net on $M$, then
\begin{equation}\label{net-prop}
M=\cup_{\xi\in\cX} B(\xi, \delta)
\quad\hbox{and}\quad
B(\xi, \delta/2)\cap B(\eta, \delta/2)=\emptyset
\quad\hbox{if}\;\; \xi\ne\eta, \; \xi, \eta\in\cX.
\end{equation}
Furthermore, $\cX$ is countable or finite
and there exists a~disjoint partition $\{A_\xi\}_{\xi\in\cX}$ of $M$ consisting
of measurable sets such that
\begin{equation}\label{B-Axi-B}
B(\xi, \delta/2) \subset A_\xi \subset B(\xi, \delta), \quad \xi\in\cX.
\end{equation}
}

Discrete versions of estimates (\ref{Comp})-(\ref{tech5}) are valid \cite{CKP}.
In particular, assuming that $\cX$ is a maximal $\delta$-net on $M$ and
$\{A_\xi\}_{\xi\in\cX}$ is a~companion disjoint partition of $M$ as above,
then
\begin{equation}\label{discr-tech11}
\sum_{\xi\in\cX} \big(1+\delta^{-1}\rho(x, \xi)\big)^{-2d-1} \le c,
\end{equation}
and if $\sigma \ge 2d+1$
\begin{equation}\label{discr-comp}
\sum_{\xi\in\cX} |A_\xi| D_{\delta,\sigma}(x, \xi) D_{\delta,\sigma}(y,\xi)
\le cD_{\delta,\sigma}( x, y).
\end{equation}
Furthermore, if $\dst \ge \delta$, then
\begin{equation}\label{basic-est}
\sum_{\xi\in \cX} \frac{|A_\xi|}{|B(\xi, \dst)|}\big(1+\dst^{-1}\rho(x, \xi)\big)^{-2d-1} \le c.
\end{equation}

\subsection{Maximal and integral operators}\label{max-operator}

The maximal operator will be an important tool for proving various estimates.
We shall use its version $\cM_t$ ($t>0$) defined by
\begin{equation}\label{def:max-op}
\cM_tf(x):=\sup_{B\ni x}\left(\frac1{|B|}\int_B |f|^t\, d\mu \right)^{1/t},
\quad x\in M,
\end{equation}
where the sup is over all balls  $B\subset M$ such that $x\in B$.

Since $\mu$ is a Radon measure on $M$ which satisfies the doubling condition (\ref{doubling})
the general theory of maximal operators applies
and the Fefferman-Stein vector-valued maximal inequality holds
(\cite{Stein}, see also \cite{GLY}):
If $0<p<\infty, 0<q\le\infty$, and
$0<t<\min\{p,q\}$ then for any sequence of functions
$\{f_\nu\}$ on $M$
\begin{equation}\label{max-ineq}
\Big\|\Bigl(\sum_{\nu}|\cM_tf_\nu(\cdot)|^q\Bigr)^{1/q} \Big\|_{\LL^p}
\le c\Big\|\Bigl(\sum_{\nu}| f_\nu(\cdot)|^q\Bigr)^{1/q}\Big\|_{\LL^p}.
\end{equation}
We shall also need the following ``integral version" of this inequality:
Let $p, q, t$ be as above.   
Then for any measurable function
$F: M\times [0,1] \to \bC$ with respect to the product measure $d\mu \times du$
one has
\begin{equation}\label{max-ineq-integral}
\Big\|\Big(
\int_{0}^{1}\Big[\cM_t(F(\cdot, u))(\cdot)\Big]^q\frac{du}{u}
\Big)^{1/q}\Big\|_{\LL^p}
\le c \Big\|\Big(
\int_{0}^{1}|F(\cdot, u)|^q \frac{du}{u}
\Big)^{1/q}\Big\|_{\LL^p}.
\end{equation}
An elaborate proof of estimate (\ref{max-ineq}) in the general setting of homogeneous type spaces
is given in \cite{GLY}.
The same proof can be easily adapted for the proof of estimate (\ref{max-ineq-integral}).
We omit the details.


\begin{remark}\label{rem:max-ineq}
{\rm
The vector-valued maximal inequality $(\ref{max-ineq})$ is usually stated and used
with $t=1$ and $p, q>1$. We find the maximal inequality in the form given in $(\ref{max-ineq})$
with $0<t<\min\{p,q\}$ more convenient. It follows immediately from the case
$t=1$ and $p, q>1$. The same observation is valid for inequality $(\ref{max-ineq-integral})$.
}
\end{remark}

A lower bound estimate on the maximal operator of the characteristic function $\ONE_{B(y, r)}$
of the ball $B(y, r)$ will be needed: 

\begin{equation}\label{max-ineq-2}
(\cM_t\ONE_{B(y, r)})(x)
\ge c\Big(1+\frac{\rho(x, y)}{r}\Big)^{-d/t},
\quad x \in M.
\end{equation}
This estimate follows easily from the doubling condition (\ref{doubling}).

\smallskip


The localization of the kernels of most integral operators that will appear in the sequel
will be controlled by the quantities $D_{\delta,\sigma}(x,y)$, defined in (\ref{def-Dxy}).
We next give estimates on the norms of such operators.


\begin{proposition}\label{prop:young}
Let $H$ be an integral operator with kernel $H(x,y)$, i.e.
$$
Hf(x)=\int_M H(x,y) f(y) d\mu(y),
\;\;\hbox{and let}\;\;\;
|H(x,y)|\le c'D_{\delta,\sigma}(x,y)
$$
for some $0<\delta\le 1$ and $\sigma \ge 2d+1$.
Then we have:

$(i)$ For  $1\le p \le \infty$
\begin{equation}\label{Lp-Lp-bound}
\|Hf\|_p \le  cc'\|f\|_p,
\quad f\in \LL^p.
\end{equation}

$(ii)$ Assuming the non-collapsing condition $(\ref{non-collapsing})$
and $1\leq p \leq q \leq \infty$
\begin{equation}\label{young-1}
\|Hf\|_q \le  cc'\delta^{d(\frac 1q - \frac 1p)}\|f\|_p,
\quad f\in \LL^p.
\end{equation}
\end{proposition}

\noindent
{\bf Proof.}
By (\ref{INT}) there exists a constant $c>0$ such that
$$
\sup_{x\in M} \|H(x, \cdot)\|_{L^1} \le c
\quad\hbox{and}\quad
\sup_{y\in M} \|H(\cdot, y)\|_{L^1} \le c
$$
Then (\ref{Lp-Lp-bound}) follows by the Schur lemma.
The proof of (\ref{young-1}) is given in \cite[Proposition 2.6]{CKP}.
$\qed$

The following useful result for products of integral
and non-integral operators is shown in \cite{CKP}.


\begin{proposition}\label{prop:prod-oper}
In the general setting of a doubling metric measure space $(M, \rho, \mu)$,
let $U, V: \LL^2 \to \LL^2$ be integral operators
and suppose that for some $ 0<  \delta  \leq 1$ and $\sigma \ge d+1$ we have
\begin{equation}\label{local-UV}
|U(x,y)| \leq c_1D_{\delta, \sigma}(x,y)
\quad\hbox{and}\quad  |V(x,y)|
\le c_2D_{\delta,\sigma}(x,y).
\end{equation}
Let $R: \LL^2\to \LL^2$ be a bounded operator, not necessarily an integral operator.
Then $U R  V $ is an integral operator
with the following bound on its kernel
\begin{equation}\label{local-URV}
|U  R  V (x,y)|
\le \|U(x,\cdot)\|_2  \| R \|_{2 \to 2}\|V(\cdot,y)\|_2
\le  \frac{cc_1c_2\|R \|_{2 \to 2}}{\big(|B(x, \delta)||B(y, \delta)|\big)^{1/2}}.
\end{equation}
\end{proposition}

\subsection{Compactly supported cut-off functions with small derivatives}\label{sec:cut-off-functions}

In the construction of frames
we shall need compactly supported $C^\infty$ functions with smallest possible derivatives.
Such functions are developed in \cite{IPX1, IPX2}.


\begin{definition}\label{cutoff-d1}
A real-valued function $\varphi\in C^\infty (\R_+)$ is said to be an admissible cut-off function if $\varphi \ne 0$,
 $\supp \varphi \subset [0, 2]$ and $\varphi^{(m)}(0)=0$ for $m\ge 1$.
Furthermore, $\varphi$ is said to be admissible of
type $(a)$, $(b)$ or $(c)$ if $\varphi$ is admissible and in addition obeys the respective condition:

$(a)$ $\varphi(t) = 1$, $t\in [0, 1]$,

$(b)$ $\supp \varphi \subset [1/2, 2]$ or

$(c)$ $\supp \varphi \subset [1/2, 2]$ and $\sum_{j=0}^\infty |\varphi(2^{-j}t)|^2=1$ for $t\in [1,\infty)$.
\end{definition}

The following proposition will be instrumental in the construction of frames.


\begin{proposition}\cite{IPX2}\label{prop:cutoff-1}
For any $0<\eps \le 1$ there exists a cut-off function $\varphi$ of type $(a)$, $(b)$ or $(c)$
such that $\|\varphi\|_\infty\le 1$ and
\begin{equation}\label{cut-off1}
\|\varphi^{(k)}\|_\infty \le 8\big(16\eps^{-1}k^{1+\eps}\big)^k
\quad\forall k\in \bN.
\end{equation}

\end{proposition}

Observe that, as shown in \cite{IPX2},
Proposition~\ref{prop:cutoff-1} is sharp in the sense that there is no cut-off function $\varphi$ such that
$
\|\varphi^{(k)}\|_\infty \le \gamma (\tilde \gamma k)^k
$
for all $k\in \bN$ no matter how large $\gamma, \tilde\gamma >0$ might be.
For more information about cut-off functions with ``small" derivatives we refer the reader to \cite{IPX2}.

\subsection{Key implications of the heat kernel properties}

We first observe that as $L$ is a nonnegative self-adjoint operator and maps real-valued to real-valued functions,
then for any real-valued, measurable and bounded function $f$ on $\R_+$ the operator $f(L):=\int_0^\infty f(\lambda)dE_\lambda$
is bounded on $L^2$, self-adjoint, and maps real-valued functions to real-valued functions.
Moreover, if $f(L)$ is an integral operator, then its kernel $f(L)(x, y)$ is real-valued and
$f(L)(y, x)=f(L)(x, y)$.

The main results in this paper will rely on the functional calculus induced by the heat kernel.
We~shall further refine the functional calculus developed in \cite{CKP} by
improving the assumptions and constant in the main space localization estimate (see Theorem 3.4 in \cite{CKP}).
Our new proof will utilize two basic ingredients:
(i) The finite speed propagation property for the solution of the associated to $L$ wave equation, and
(ii) A non-smooth functional calculus estimate. 


In this theory, the following {\bf Davies-Gaffney estimate} for the heat kernel plays a significant role:
\begin{equation}\label{davies-gaffney}
|\langle P_t f_1, f_2\rangle| \le \exp\Big\{-\frac{\cdg r^2}{t}\Big\} \|f_1\|_2\|f_2\|_2,
\quad t>0,
\end{equation}
for all open sets $U_j \subset M$ and $f_j\in \LL^2(M)$ with $\supp f_j\subset U_j$, $j=1, 2$,
where $r:=\rho(U_1, U_2)$ and $\cdg>0$ is a constant.
We next show that the Davies-Gaffney estimate is a consequence of the conditions on
the heat kernel stipulated in \S \ref{introduction}.


\begin{proposition}\label{prop:davies-gaffney}
In the setting described in the introduction
the Davies-Gaffney estimate $(\ref{davies-gaffney})$ holds with $\cdg =\cstar$,
where $\cstar$ is the constant from $(\ref{Gauss-local})$.
\end{proposition}

\noindent
{\bf Proof.}
We shall proceed in the spirit of [7].
The proof relies on the following version of the Phragm\'{e}n-Lindel\"{o}f theorem (see [7], Proposition 2.2):
{\em
Suppose the function $F$ is holomorphic in $\bC_+:=\{z: {\rm Re}\, z >0\}$ and verifies the following conditions:
$$
|F(z)|\le B \quad\hbox{for} \quad  z\in \bC_+ \quad\hbox{and}\quad
|F(t)|\le Ae^{-\frac{\gamma}{t}}\quad\hbox{for}\quad t>0,
$$
where $A, B, \gamma >0$ are constants.
Then
$|F(z)|\le Be^{-{\rm Re}\,\frac{\gamma}{z}}$ for $z\in \bC_+$.
}

This theorem can be slightly improved as follows:

\noindent
{\bf Claim:}
{\em
Suppose $F$ is holomorphic in $\bC_+:=\{z: {\rm Re}\, z >0\}$
and there exist constants $B, \gamma >0$ such that
$|F(z)|\le B$ for $z\in \bC_+$,
and for any $\varepsilon >0$ there exists a constant $A_\varepsilon>0$ such that
$|F(t)|\le A_\varepsilon e^{-\frac{\gamma-\varepsilon}{t}}$ for $0<t\le 1$.
Then
$|F(z)|\le Be^{-{\rm Re}\,\frac{\gamma}{z}}$ for $z\in \bC_+$.
}

\noindent
{\bf Proof.}
Clearly,
$|F(t)|\le B \le B e^{\gamma-\varepsilon} e^{-\frac{\gamma-\varepsilon}{t}}$
for $t>1$ and $0<\varepsilon <\gamma$.
Therefore,
$$
|F(t)|\le B \le \max\{A_\varepsilon, B e^{\gamma-\varepsilon}\} e^{-\frac{\gamma-\varepsilon}{t}}
\quad \hbox{for all $t>0$ and $0<\varepsilon <\gamma$.}
$$

By applying the Phragm\'{e}n-Lindel\"{o}f theorem from above we conclude that for any $\varepsilon >0$
we have
$|F(z)|\le Be^{-{\rm Re}\,\frac{\gamma-\varepsilon}{z}}$ for $z\in \bC_+$.
This readily implies the claim.
$\qed$

\smallskip

As shown in [7, Lemma 3.1] to prove Proposition~\ref{prop:davies-gaffney}
it suffices to prove (\ref{davies-gaffney}) in the case when $U_1$, $U_2$
are arbitrary balls $B_1$, $B_2$.
Let $B_j=B(a_j, r_j)$ and $f_j\in L^2(B_j)$, $j=1, 2$. Write $r:=\rho(B_1, B_2)$.

Since $L$ is a nonnegative self-adjoint operator,
$P_z$ is holomorphic in $\bC_+$ and $\|P_z\|_{2\to 2} \le 1$ $\forall z\in \bC_+$.
Therefore, the function $F(z):= \langle P_zf_1, f_2\rangle$ is holomorphic in $\bC_+$ and
\begin{equation}\label{est-Fz}
|F(z)| \le \|f_1\|_2\|f_2\|_2 \quad \hbox{for $z\in \bC_+$.}
\end{equation}

We next estimate $|F(t)|$ for $0<t\le 1$.
Using (1.3) we infer
\begin{align*}
&|\langle P_tf_1, f_2\rangle|
\le \int_{B_1}\int_{B_2} |P_t(x, y)||f_1(x)||f_2(y)| d\mu(x) d\mu(y)\\
& \quad
\le C^\star e^{-\frac{c^\star r^2}{t}}
\int_{B_1}\frac{|f_1(x)|}{|B(x, \sqrt{t})|^{1/2}} d\mu(x)
\int_{B_2}\frac{|f_2(x)|}{|B(y, \sqrt{t})|^{1/2}} d\mu(y)\\
& \quad
\le C^\star e^{-\frac{c^\star r^2}{t}} \|f_1\|_2 \|f_2\|_2
\Big(\int_{B_1}|B(x, \sqrt{t})|^{-1} d\mu(x)\Big)^{1/2}
\Big(\int_{B_2}|B(x, \sqrt{t})|^{-1} d\mu(x)\Big)^{1/2},
\end{align*}
where in the last inequality we used the Cauchy-Schwartz inequality.

If $\sqrt{t} \le r_1$ and $x\in B_1$, then $B(a, r_1) \subset B(x, 2r_1)$ and hence
using (\ref{doubling}) we obtain
$|B(a, r_1)| \le |B(x, \frac{2r_1\sqrt{t}}{\sqrt{t}})|
\le c_0\big(\frac{2r_1}{\sqrt{t}}\big)^d |B(x, \sqrt{t})|$,
which leads to 
\begin{equation}\label{est-int-B1}
\Big(\int_{B_1}|B(x, \sqrt{t})|^{-1} d\mu(x)\Big)^{1/2}
\le c \big(\frac{r_1}{\sqrt{t}}\big)^{d/2}
\le c_\varepsilon e^{\frac{\varepsilon}{t}}
\quad \hbox{$\forall \varepsilon >0$ and $\sqrt{t} < r_1$,}
\end{equation}
where $c_\varepsilon >0$ is a sufficiently large constant.

If $\sqrt{t} > r_1$ and $x\in B_1=B(a, r_1)$, then $B(a, r_1) \subset B(x, 2\sqrt{t})$
and hence by using (\ref{doubling-0}) we get
$|B_1| \le c_0|B(x, \sqrt{t})|$, implying
$\int_{B_1}|B(x, \sqrt{t})|^{-1} d\mu(x) \le c_0$.

Therefore, estimate (\ref{est-int-B1}) 
holds for all $0<t\le 1$, if $c_\varepsilon$ is sufficiently large.

Exactly in the same manner one shows that for every $\varepsilon >0$
there exists a constant $\tilde c_\varepsilon>0$ such that
$\big(\int_{B_2}|B(y, \sqrt{t})|^{-1} d\mu(y)\big)^{1/2} \le \tilde c_\varepsilon e^{\frac{\varepsilon}{t}}$
for $0<t\le 1$.


Putting the above estimates together we conclude that for every $\varepsilon >0$
there exists a constant $A_\varepsilon>0$ such that
$|F(t)| \le A_\varepsilon e^{-\frac{c^\star r^2-\varepsilon}{t}}$ for $0<t\le 1$.

Finally, from this and (\ref{est-Fz}), employing the claim from above, we conclude that
$|F(z)| \le e^{-{\rm Re} \{\frac{c^\star r^2}{z}\}}$ for $z\in \bC_+$,
which completes the proof.
$\qed$

\smallskip



\smallskip

In going further, observe that as proved in \cite{CS} (Theorem 3.4),
the Davies-Gaffney estimate (\ref{davies-gaffney}) implies
(in fact, it is equivalent to)
the {\bf finite speed propagation property}:
\begin{equation}\label{finite-speed}
\big\langle \cos(t\sqrt{L})f_1, f_2 \big\rangle=0,
\quad 0< \ct t<r,
\quad \ct:=\frac{1}{2\sqrt{\cstar}},
\end{equation}
for all open sets $U_j \subset M$, $f_j\in \LL^2(M)$, $\supp f_j\subset U_j$,
$j=1, 2$, where $r:=\rho(U_1, U_2)$.

\smallskip

We next use this to derive important information about
the kernels of operators of the form $f(\delta\sqrt{L})$ whenever $\hat f$ is band limited.
Here $\hat f(\xi):=\int_\R f(t)e^{-it\xi}dt$.


\begin{proposition}\label{prop:finite-sp}
Let $f$ be even, $\supp \hat f \subset [-A, A]$ for some $A>0$,
and $\hat f\in W^m_1$ for some $m>d$, i.e. $\|\hat f^{(m)}\|_1 <\infty$.
Then for $\delta>0$ and $x, y\in M$
\begin{equation}\label{finite-speed-2}
f(\delta\sqrt{L})(x, y) = 0
\quad\hbox{if}\quad
\ct \delta A < \rho(x, y).
\end{equation}
\end{proposition}

\noindent
{\bf Proof.}
From functional calculus and the Fourier inversion formula
$$
f(\delta\sqrt{L}) =\frac{1}{\pi}\int_0^A\hat f(\xi)\cos(\xi\delta\sqrt{L}) d\xi.
$$
Fix $x, y\in M$, $x\ne y$, and let $\ct\delta A<\rho(x, y)$.
Choose $\epsilon>0$ so that  $\ct\delta A < \rho (x, y)-2\epsilon$
and let
$g_1:= |B(x, \epsilon)|^{-1}\ONE_{B(x, \epsilon)}$ and
$g_2:= |B(y, \epsilon)|^{-1}\ONE_{B(y, \epsilon)}$.
Then from above and (\ref{finite-speed}) we derive
\begin{equation}\label{finite-speed-3}
\big\langle f(\delta\sqrt{L}) g_1, g_2\big\rangle
= \frac{1}{\pi}\int_0^A \hat f(\xi)\big\langle \cos(\xi\delta\sqrt{L})g_1, g_2\big\rangle d\xi = 0,
\end{equation}
using that $\ct\delta A < \rho (x, y)-2\epsilon \le \rho\big(B(x, \epsilon), B(y, \epsilon)\big)$.
On the other hand, it easily follow from Theorem 3.7 in \cite{CKP} that the kernel of $f(\delta\sqrt{L})$
is continuous and, therefore,
$$
\big\langle f(\delta\sqrt{L}) g_1, g_2 \big\rangle
=\int_M\int_M f(\delta\sqrt{L})(u, v)g_1(u)g_2(v) d\mu(u)d\mu(v)
\to f(\delta\sqrt{L})(x, y)
$$
as $\epsilon \to 0$.
This and (\ref{finite-speed-3}) imply (\ref{finite-speed-2}).
$\qed$

\smallskip

Another important ingredient for our further development will be the following (Theorem~3.7 in \cite{CKP})


\begin{proposition}\label{prop:rough-kernels}
Let $f$ be a bounded measurable function on $\bR_+$ with
$\supp f \subset [0,\tau]$ for some $\tau \ge 1$.
Then $f(\sqrt L)$ is an integral operator with kernel $f(\sqrt L)(x, y)$ satisfying
\begin{equation}\label{rough1}
| f(\sqrt L) (x,y) |
\le  \frac{\cf\| f \|_\infty}
{\big(|B(x, \tau^{-1})|| B(y, \tau^{-1})|\big)^{1/2}},\quad x,y\in M,
\end{equation}
where $\cf>0$ depends only on the constants $c_0, \Cstar, \cstar$ from $(\ref{doubling-0})-(\ref{Gauss-local})$.
\end{proposition}

This proposition also follows by the properties of the heat kernel $p_t(x, y)$ from~\S\ref{introduction}.


\begin{remark}\label{rem:davies-gaffney}
{\rm
As is well known the Davies-Gaffney estimate (\ref{davies-gaffney}) is weaker than assuming
the Gaussian bound (\ref{Gauss-local}) on the heat kernel
and also estimate (\ref{rough1}) is weaker than (\ref{Gauss-local}).
However, it can be shown by combining results from \cite{CS} and \cite{Ouhabaz}
that the Davies-Gaffney estimate (\ref{davies-gaffney}), estimate (\ref{rough1}),
and the doubling condition (\ref{doubling-0}) imply (\ref{Gauss-local}).
Therefore, deriving in the next section the main localization estimate (\ref{main-loc-ker1}) of the functional calculus
by using the finite speed propagation property (\ref{finite-speed}) and (\ref{rough1}) instead of (\ref{Gauss-local})
we essentially do not weaken our assumptions.
}
\end{remark}

\section{Smooth functional calculus induced by the heat kernel}\label{sec:smooth-func-calculus}
\setcounter{equation}{0}

We shall make heavy use in this paper of the functional calculus developed in \cite{CKP}
in the setting described in the introduction.
We next improve and extend some basic results from \S 3 in \cite{CKP}.

\subsection{Kernel localization and H\"{o}lder continuity}\label{sec:func-calculus}

We first establish an improved version of Theorem~3.4 in \cite{CKP}.
The main new feature is the improved control on the constants, which will be important
for our subsequent developments.


\begin{theorem}\label{thm:main-local-kernels}
Let $f\in C^k(\bR_+)$, $k \ge d+1$, $\supp f\subset [0, R]$ for some $R\ge 1$, and
$f^{(2\nu+1)}(0)=0$ for $\nu \ge 0$ such that $2\nu+1 \le k$.
Then $f(\delta \sqrt L)$, $0<\delta\le 1$, is an integral operator with kernel
$f(\delta \sqrt L)(x, y)$ satisfying
\begin{equation}\label{main-loc-ker1}
\big|f(\delta \sqrt L)(x, y)\big| \le c_k D_{\delta, k}(x,y)
\quad \hbox{and}
\end{equation}
\begin{equation}\label{main-loc-ker2}
\big|f(\delta \sqrt L)(x, y)  -  f(\delta \sqrt L)(x,y')\big|
\le c_k'\Big(\frac{\rho(y,y')}{\delta}\Big)^\alpha D_{\delta, k}(x,y)
\;\; \hbox{if}\;\; \rho (y, y') \le \delta.
\end{equation}
Here
$D_{\delta, k}(x,y)$ is from $(\ref{def-Dxy})$,
\begin{equation}\label{const-ck}
c_k = c_k(f)= R^d\big[(c_1k)^k\|f\|_{L^\infty} + (c_2R)^k\|f^{(k)}\|_{L^\infty}\big],
\end{equation}
where
$c_1, c_2>0$ depend only on the constants $c_0, \Cstar, \cstar$ from $(\ref{doubling-0})-(\ref{lip})$ and
$c_k'=c_3c_kR^\alpha$
with $c_3>0$ depending only on $c_0, \Cstar, \cstar$ and $k$;
as before $\alpha>0$ is the constant from $(\ref{lip})$.
Furthermore,
\begin{equation}\label{main-loc-ker3}
\int_M  f(\delta \sqrt L)(x, y) d\mu(y) = f(0).
\end{equation}
\end{theorem}


\begin{remark}\label{rem:even}
{\rm
The condition $f^{(2\nu+1)}(0)=0$ for $\nu \ge 0$ such that $2\nu+1 \le k$
simply says that if $f$ is extended as an even function to $\R$ $(f(-\lambda)= f(\lambda))$,
then $f\in C^k(\bR)$.
}
\end{remark}


\noindent
{\bf Proof.}
It suffices to prove the theorem in the case $R=1$. Then in general it follows by rescaling.

Assume that $f$ satisfies the hypotheses of the theorem with $R=1$ and denote again by $f$
its even extension to $\R$.
As already observed in Remark~\ref{rem:even}, $f\in C^k(\R)$.
The idea of the proof is to approximate $f$ by a band limited function $f_A$
and then utilize Propositions~\ref{prop:finite-sp}-\ref{prop:rough-kernels}.

Set
$$
\hat\phi:= \ONE_{[-\frac{1}{2}-\delta, \frac{1}{2}+\delta]}*\underbrace{H_\delta*\cdots *H_\delta}_{k+1},
\quad\hbox{where}\;\; H_\delta:=(2\delta)^{-1}\ONE_{[-\delta, \delta]},
\quad \delta:= \frac{1}{2(k+2)}.
$$
Clearly,
$\hat\phi$ is even, $\supp \hat\phi\subset [-1, 1]$, $0\le \hat\phi \le 1$,
$\hat\phi(\xi)=1$ for $\xi\in [-1/2, 1/2]$,
and
\begin{equation}\label{der-phi}
\|\hat\phi^{(\nu)}\|_\infty\le \delta^{-\nu} \le (2(k+2))^\nu \le (4k)^\nu
\quad\hbox{for }\; \nu=0, 1, \dots, k+1.
\end{equation}
The last inequality follows just as in \cite[Theorem 1.3.5]{H}.

Denote
$\phi(t):=(2\pi)^{-1}\int_\R\hat\phi(\xi) e^{i\xi t}d\xi$
and set $\phi_A(t):=A\phi(At)$, $A>0$.
Then $\widehat{\phi_A}(\xi)= \hat\phi(\xi/A)$
and hence $\supp \widehat{\phi_A}\subset [-A, A]$.

Now, consider the function
$f_A:= f*\phi_A$.
Clearly, $\widehat{f_A}=\hat f \widehat{\phi_A}$, which implies
$\supp \widehat{f_A}\subset [-A, A]$.
Since $f$ and $\phi$ are even, then $f_A$ is even.
Furthermore,
\begin{align*}
f(t)-f_A(t)
&= (2\pi)^{-1}\int_\R\hat f(\xi)(1-\hat\phi(\xi/A))e^{i\xi t} d\xi\\
&= (2\pi)^{-1}A^{-k}\int_\R \xi^k \hat f(\xi)\hat \hh(\xi/A)e^{i\xi t} d\xi,
\end{align*}
where
$\hat \hh(\xi) = (1-\hat \phi(\xi))\xi^{-k}$.
Set
$\hh_A(t):=A\hh(At)$ and note that $\widehat{\hh_A}(\xi)= \hat \hh(\xi/A)$.
Also, observe that
$\widehat{f^{(k)}}(\xi) = (i\xi)^k\hat f(\xi)$.
From all of the above we derive
\begin{equation}\label{est-f-fA}
\|f-f_A\|_\infty \le A^{-k}\|f^{(k)}*\hh_A\|_\infty
\le A^{-k}\|f^{(k)}\|_\infty\|\hh_A\|_{L^1}.
\end{equation}
Clearly,
$$
t^2\hh(t)= \frac{i^2}{2\pi} \int_\R \Big(\frac{d}{d\xi}\Big)^2\hat \hh(\xi) e^{i\xi t} d\xi
\quad\hbox{and}\quad
\Big|\Big(\frac{d}{d\xi}\Big)^2\hat \hh(\xi)\Big| \le c^k (1+|\xi|)^{-k-2},
$$
and hence
$
|\hh(t)| \le c_1^k(1+|t|)^{-2},
$
which leads to
$\|\hh_A\|_{L^1} = \|\hh\|_{L^1} \le c^k$,
where $c>1$ is an absolute constant.
From this and (\ref{est-f-fA}) we get
\begin{equation}\label{est-f-fA-2}
\|f-f_A\|_\infty \le c^k A^{-k}\|f^{(k)}\|_\infty.
\end{equation}

We next estimate $|f(t)-f_A(t)|$ for $t>1$.
For this we need an estimate on the localization of $|\phi_A(t)|$.
Since $\supp \hat\phi \subset [-1, 1]$
we have
$
\phi(t) = \frac{1}{2\pi}\int_{-1}^1 \hat\phi(\xi)e^{i\xi t} d\xi
$
and integrating by parts $k+1$ times we obtain
$$
\phi(t) = \frac{(-1)^{k+1}}{2\pi(it)^{k+1}}\int_{-1}^1 \hat\phi^{(k+1)}(\xi)e^{i\xi t} d\xi.
$$
Therefore, using (\ref{der-phi})
$$
|t^{k+1}\phi(t)| \le \|\hat \phi^{(k+1)}\|_\infty
\le (4k)^{k+1}.
$$
In turn, this and the obvious estimate $\|\phi\|_\infty \le 2$
imply
$
|\phi(t)| \le (c'k)^{k}(1+|t|)^{-k-1},
$
where $c'>4$ is an absolute constant.
Hence,
\begin{equation}\label{est-phiA}
|\phi_A(t)| \le c(k)A(1+A|t|)^{-k-1},
\quad c(k)=(c'k)^{k}.
\end{equation}
Using this and $\supp f\subset [-1, 1]$ we obtain for $t>1$
\begin{align*}
|f(t)-f_A(t)|
&= |f_A(t)| = |f*\phi_A(t)| \le \int_{-1}^1|f(y)||\phi_A(y-t)|dy\\
& \le \|f\|_\infty \int_{t-1}^{t+1}|\phi_A(u)|du
\le c(k)\|f\|_\infty \int_{t-1}^{t+1}A(1+Au)^{-k-1}du\\
&\le c(k)\|f\|_\infty \int_{A(t-1)}^\infty(1+v)^{-k-1}du
\le c(k)A^{-k}\|f\|_\infty(t-1)^{-k}.
\end{align*}
This yields
\begin{equation}\label{est-f-fA-3}
|f(t)-f_A(t)| \le (3c'k)^{k}A^{-k}\|f\|_\infty(t+1)^{-k}
\quad\hbox{for $t\ge 2$.}
\end{equation}

In our next step we utilize Proposition~\ref{prop:rough-kernels}.
For this we need to apply a decomposition of unity argument.
Choose $\ph_0\in C^{\infty}(\bR_+)$ so that $\supp \ph_0\subset [0, 2]$, $0\le \ph_0 \le 1$, and
$\ph_0(\lambda)=1$ for $\lambda \in [0, 1]$.
Let $\ph(\lambda):=\ph_0(\lambda)-\ph_0(2\lambda)$.
Note that $\ph\in C^{\infty}(\bR)$ and $\supp \ph \subset [1/2, 2]$.
Set
$\ph_j(\lambda):=\ph(2^{-j}\lambda)$, $j\ge 1$.
Then
$\sum_{j \ge 0}\ph_j(\lambda)=1$ for $\lambda \in \bR_+$ and hence
$$
f(\lambda)- f_A(\lambda) = \sum_{j \ge 0}[f(\lambda)- f_A(\lambda)]\ph_j(\lambda),
$$
which implies
\begin{equation}\label{phi-decomp}
f(\delta \sqrt L) - f_A(\delta \sqrt L)
= \sum_{j \ge 0}\big[f(\delta \sqrt L)- f_A(\delta \sqrt L)\big]\ph_j(\delta \sqrt L), \;\; \delta>0,
\end{equation}
where the convergence is strong (in the $L^2 \to L^2$ operator norm).

Let $x, y\in M$, $x\ne y$, and assume $\rho(x, y) \ge \delta$.
Choose $A>0$ so that
\begin{equation}\label{choose-A}
\frac{\rho(x, y)}{2\delta}\le \ct A <\frac{\rho(x, y)}{\delta}.
\end{equation}
Since $\supp \widehat{f_A} \subset [-A, A]$ and $\widehat{f_A}\in W^{k+1}_1$ with $k\ge d+1$,
by Proposition~\ref{prop:finite-sp}
$f_A(\delta \sqrt L)(x, y)=0$ and hence
$$
f(\delta \sqrt L)(x, y)=f(\delta \sqrt L)(x, y) - f_A(\delta \sqrt L)(x, y).
$$
Denote briefly
$\FF_j(\lambda):= (f(\lambda)- f_A(\lambda))\ph_j(\lambda)$.
Then the above and (\ref{phi-decomp}) lead to
$$
|f(\delta \sqrt L)(x, y)|
\le \sum_{j \ge 0} |\FF_j(\delta \sqrt L)(x, y)|.
$$
Note that, $\supp \FF_0 = \supp \varphi_0 \subset [0, 2]$
and
$\supp \FF_j = \supp \varphi_j \subset [2^{j-1}, 2^{j+1}]$, $j\ge 1$.

For $j=0, 1$ we use (\ref{est-f-fA-2}) to obtain
$\|\FF_j\|_\infty \le c^k A^{-k}\|f^{(k)}\|_\infty$
and applying Proposition~\ref{prop:rough-kernels}
$$
|\FF_j(\delta \sqrt L)(x, y)|
\le \frac{\cf c^k A^{-k}\|f^{(k)}\|_\infty}{\big(|B(x, \delta/2)|| B(y, \delta/2)|\big)^{1/2}}
\le \frac{c_0\cf (2c\ct)^k \|f^{(k)}\|_\infty}
{\big(|B(x, \delta)|| B(y, \delta)|\big)^{1/2}\big(1+\frac{\rho(x, y)}{\delta}\big)^k},
$$
where we used (\ref{choose-A}) and
$|B(\cdot, \delta)|\le c_0|B(\cdot, \delta/2)|$ by (\ref{doubling-0}).

For $j\ge 2$ we use (\ref{est-f-fA-3}) to obtain
$\|\FF_j\|_\infty \le (3c' k)^k A^{-k}\|f\|_\infty 2^{-k(j-1)}$
and applying again Proposition~\ref{prop:rough-kernels} we get
\begin{align*}
|\FF_j(\delta \sqrt L)(x, y)|
& \le \frac{\cf(6 c'\ct k)^k\|f\|_\infty 2^{-k(j-1)}}
{\big(|B(x, \delta 2^{-j-1})|| B(y, \delta 2^{-j-1})|\big)^{1/2}\big(1+\frac{\rho(x, y)}{\delta}\big)^k}\\
&  \le \frac{c_0\cf(6 c'\ct k)^k\|f\|_\infty 2^{-k(j-1)}2^{d(j+1)}}
{\big(|B(x, \delta)|| B(y, \delta)|\big)^{1/2}\big(1+\frac{\rho(x, y)}{\delta}\big)^k}.
\end{align*}
Here we used again (\ref{choose-A}) and
$|B(\cdot, \delta)|\le c_02^{(j+1)d}|B(\cdot, \delta 2^{-j-1})|$ by (\ref{doubling}).

We sum up the above estimates taking into account that $k\ge d+1$ and obtain
$$
|f(\delta \sqrt L)(x, y)|
\le \frac{(c_1 k)^{k}\|f\|_\infty + c_2^k\|f^{(k)}\|_\infty}
{\big(|B(x, \delta)|| B(y, \delta)|\big)^{1/2}\big(1+\frac{\rho(x, y)}{\delta}\big)^k}
\quad\hbox{if } \rho(x, y) \ge \delta.
$$
Whenever $\rho(x, y) < \delta$, this estimate is immediate from Proposition~\ref{prop:rough-kernels} with $c\|f\|_\infty$
in the numerator.
The proof of estimate (\ref{main-loc-ker1}) is complete.

For the proof of (\ref{main-loc-ker2}) we write
$$
f(\delta \sqrt L)(x, y) = \int_M f(\delta \sqrt L)e^{\delta^2 L}(x, u)e^{-\delta^2 L}(u, y) d\mu(u)
$$
and proceed further exactly as in the proof of (3.3) in \cite{CKP} using (\ref{main-loc-ker1}) and
the H\"{o}lder continuity of the heat kernel, stipulated in (\ref{lip}).
$\qed$


\begin{remark}\label{rem:kd}
It is readily seen that Theorem~\ref{thm:main-local-kernels} holds under the slightly weaker condition
$k>d$ rather than $k\ge d+1$, but then the constants $c_k, c_k'$ will also depend on $k-d$.
\end{remark}

Now, we would like to make a step forward and free the function $f$ in the hypothesis of
Theorem~\ref{thm:main-local-kernels} from the restriction of being compactly supported.


\begin{theorem}\label{thm:S-local-kernels}
Suppose $f\in C^k(\bR_+)$, $k\ge d+1$,
$$ 
\hbox{
$|f^{(\nu)}(\lambda)|\leq C_k(1+\lambda)^{-r}$ for $\lambda>0$ and $0\le \nu\le k$, where $r\ge k+d+1$,
}
$$ 
and $f^{(2\nu+1)}(0)=0$ for $\nu\ge 0$ such that $2\nu+1 \le k$.
Then $f(\delta \sqrt L)$ is an integral operator with kernel $f(\delta \sqrt L)(x, y)$
satisfying $(\ref{main-loc-ker1})$-$(\ref{main-loc-ker2})$,
where the constants $c_k,c_k'$ depend on $k, d, \alpha$,
but also depend linearly on $C_k$.
\end{theorem}

\noindent
{\bf Proof.}
As in the proof of Theorem~\ref{thm:main-local-kernels},
choose $\ph_0\in C^{\infty}(\bR_+)$ so that $0\le \ph_0 \le 1$,
$\ph_0(\lambda)=1$ for $\lambda \in [0, 1]$, and $\supp \ph_0\subset [0, 2]$.
Let $\ph(\lambda):=\ph_0(\lambda)-\ph_0(2\lambda)$
and set $\ph_j(\lambda):=\ph(2^{-j}\lambda)$, $j\ge 1$.
Clearly,
$\sum_{j \ge 0}\ph_j(\lambda)=1$ for $\lambda \in \bR_+$
and hence
\begin{equation}\label{f-decomp}
f(\lambda)= \sum_{j \ge 0}f(\lambda)\ph_j(\lambda)
\Longrightarrow
f(\delta \sqrt L)= \sum_{j \ge 0}f(\delta \sqrt L)\ph_j(\delta \sqrt L), \;\; \delta>0,
\end{equation}
where the convergence is strong.
Set $h_j(\lambda):=f(2^j\lambda)\ph(\lambda)$, $j\ge 0$, and $h_0(\lambda):= f(\lambda)\ph_0(\lambda)$.
Then $h_j(2^{-j}\delta\sqrt L)=f(\delta \sqrt L)\ph_j(\delta \sqrt L)$.

By the hypotheses of the theorem it follows that for $j\ge 1$
\begin{align*}
\|h_j^{(k)}\|_{L^\infty}
\le c 2^{jk}\max_{0\le \ell\le k} \|f^{(\ell)}(2^j\cdot)\|_{L^\infty[1/2, 2]}
\le c2^{jk}2^{-jr}
\le c2^{-j(d+1)}
\end{align*}
and $\|h_j\|_{L^\infty} \le c2^{-jr} \le c2^{-j(d+1)}$.
We use this and Theorem~\ref{thm:main-local-kernels} to conclude that
$f(\delta \sqrt L)\ph_j(\delta \sqrt L)$ is an integral operator with kernel
satisfying
\begin{align*}
\big|f(\delta \sqrt L)\ph_j(\delta \sqrt L)(x, y)\big|
&= |h_j(2^{-j}\delta\sqrt L)(x, y)|
\le c\frac{2^{-j(d+1)}\big(1+\delta^{-1}2^j\rho(x, y)\big)^{-k}}
{\big(|B(x,\delta2^{-j})||B(y,\delta2^{-j})|\big)^{1/2}}\\
&\le c\frac{2^{-j}\big(1+\delta^{-1}2^j\rho(x, y)\big)^{-k}}
{\big(|B(x,\delta)||B(y,\delta)|\big)^{1/2}}.
\end{align*}
Here for the latter estimate we used (\ref{doubling}).
Exactly as above we derive a similar estimate when $j=0$.
Finally, summing up we obtain
$$
|f(\delta \sqrt L)(x, y)|
\le c\big(|B(x,\delta)||B(y,\delta)|\big)^{-1/2}
\sum_{j\ge 0}2^{-j}\big(1+\delta^{-1}2^j\rho(x, y)\big)^{-k}
\le cD_{\delta, k}(x, y),
$$
which proves (\ref{main-loc-ker1}).
The proof of (\ref{main-loc-ker2}) goes along similar lines and will be omitted.
$\qed$


\begin{corollary}\label{thm:LS-local-kernels}
Suppose $f\in C^\infty(\bR_+)$,
$|f^{(\nu)}(\lambda)|\leq C_{\nu, r}(1+\lambda)^{-r}$ for all $\nu, r\ge 0$ and $\lambda>0$,
and $f^{(2\nu+1)}(0)=0$ for $\nu \ge 0$.
Then for any $m\ge 0$ and $\delta>0$ the operator
$L^m f(\delta \sqrt L)$ is an integral operator with kernel $L^m f(\delta \sqrt L)(x, y)$
having the property that for any $\sigma>0$ there exists a constant $c_{\sigma,m}>0$ such that
\begin{equation}\label{LS-loc-ker1}
\big|L^mf(\delta \sqrt L)(x, y)\big| \le c_{\sigma,m} \delta^{-2m}D_{\delta, \sigma}(x,y)
\quad \hbox{and}
\end{equation}
\begin{equation}\label{LS-loc-ker2}
\big|L^mf(\delta \sqrt L)(x, y)  -  L^mf(\delta \sqrt L)(x,y')\big|
\le c_{\sigma,m}\delta^{-2m}\Big(\frac{\rho(y,y')}{\delta}\Big)^\alpha D_{\delta, k}(x,y),
\end{equation}
whenever $\rho (y, y') \le \delta$.
\end{corollary}

\noindent
{\bf Proof.}
Let $h(\lambda):=\lambda^{2m}f(\lam)$. Then $h(\delta\sqrt{L})= \delta^{2m}L^m f(\delta\sqrt{L})$.
It is easy to see that
$h^{(2\nu+1)}(0)=0$ for all $\nu \ge 0$.
Then the corollary follows readily by Theorem~\ref{thm:S-local-kernels} applied to $h$.
$\qed$

\subsection{Band-limited sub-exponentially localized kernels}\label{sec:local-kernels}

The kernels of operators of the form $\varphi(\delta \sqrt{L})$ with
sub-exponential space localization and $\varphi\in C^\infty_0(\R_+)$
will be the main building blocks in constructing our frames.


\begin{theorem}\label{thm:band-sub-exp}
For any $0<\eps<1$ there exists a cut-off function $\varphi$ of any type,
$(a)$ or $(b)$ or $(c)$ $($Definition~$\ref{cutoff-d1}$$)$
such that for any $\delta >0$
\begin{equation}\label{band-sub-exp}
|\varphi(\delta \sqrt{L})(x, y)|
\le \frac{c_1\exp\big\{-\kap\big(\frac{\rho(x, y)}{\delta}\big)^{1-\eps}\big\}}
{\big(|B(x, \delta)||B(y, \delta)|\big)^{1/2}},
\quad x, y\in M,
\end{equation}
and
\begin{equation}\label{holder-sub-exp}
|\varphi(\delta \sqrt{L})(x, y)-\varphi(\delta \sqrt{L})(x, y')|
\le \frac{c_2\big(\frac{\rho(y,y')}{\delta}\big)^\alpha
\exp\big\{-\kap\big(\frac{\rho(x, y)}{\delta}\big)^{1-\eps}\big\}}
{\big(|B(x, \delta)||B(y, \delta)|\big)^{1/2}}
\;\;\hbox{if } \rho(y, y')\le \delta,
\end{equation}
where $c_1, \kap >0$ depend only on $\eps$ and
the constants $c_0, \Cstar, \cstar$ from $(\ref{doubling-0})-(\ref{lip})$;
$c_2>0$ also depends on $\alpha$.
Furthermore, for any $m\in \bN$
\begin{equation}\label{band-sub-exp-m}
|L^m\varphi(\delta \sqrt{L})(x, y)|
\le \frac{c_3\delta^{-2m}\exp\Big\{-\kap\big(\frac{\rho(x, y)}{\delta}\big)^{1-\eps}\Big\}}
{\big(|B(x, \delta)||B(y, \delta)|\big)^{1/2}},
\quad x, y\in M,
\end{equation}
with $c_3>0$ depending on $\eps, c_0, \Cstar, \cstar$, and $m$.
\end{theorem}

\noindent
{\bf Proof.}
Let $0<\eps<1$. Then by Proposition~\ref{prop:cutoff-1} there exists a cut-off function $\varphi$
of any type ($(a)$ or $(b)$ or $(c)$) such that
$\|\varphi^{(k)}\|_\infty \le (ck)^{k(1+\eps)}$ for all $k\in \bN$
and $\|\varphi\|_\infty\le 1$.
Now, using Theorem~\ref{thm:main-local-kernels} we obtain
$$
|\varphi(\delta \sqrt{L})(x, y)|
\le \frac{(Ck)^{k(1+\eps)}}
{\big(|B(x, \delta)||B(y, \delta)|\big)^{1/2}\big(1+\delta^{-1}\rho(x, y)\big)^k}
\quad \forall x, y\in M, \forall k\in \bN.
$$
Here $C>1$ depends only on $\eps, c_0, \Cstar, \cstar$.
From this we infer
$$
|\varphi(\delta \sqrt{L})(x, y)|
\le \frac{e^{-k}}{\big(|B(x, \delta)||B(y, \delta)|\big)^{1/2}}
\quad\hbox{if}\quad \delta^{-1}\rho(x, y) \ge e(Ck)^{1+\eps}=: \cast k^{1+\eps}.
$$
Assume $\delta^{-1}\rho(x, y) \ge 4\cast$ and choose $k\in\bN$ so that
$k\le \big(\frac{\delta^{-1}\rho(x, y)}{\cast}\big)^{1/(1+\eps)} < k+1$.
Then from above
\begin{align*}
|\varphi(\delta \sqrt{L})(x, y)|
\le \frac{\exp\big\{-\frac{1}{2}\big(\frac{\delta^{-1}\rho(x, y)}{\cast}\big)^{1/(1+\eps)}\big\}}
{\big(|B(x, \delta)||B(y, \delta)|\big)^{1/2}}
\le \frac{\exp\big\{-\kap\big(\frac{\rho(x, y)}{\delta}\big)^{1-\eps}\big\}}
{\big(|B(x, \delta)||B(y, \delta)|\big)^{1/2}},
\end{align*}
provided $\delta^{-1}\rho(x, y) \ge 4\cast$.
In the case $\delta^{-1}\rho(x, y) < 4\cast$ we get from Proposition~\ref{prop:rough-kernels}
$$
|\varphi(\delta \sqrt{L})(x, y)|
\le \frac{c'}{\big(|B(x, \delta)||B(y, \delta)|\big)^{1/2}}
\le \frac{c\exp\big\{-\kap\big(\frac{\rho(x, y)}{\delta}\big)^{1-\eps}\big\}}
{\big(|B(x, \delta)||B(y, \delta)|\big)^{1/2}}
$$
with $c=c'\exp\{4\kap \cast\}$.
This completes the proof of (\ref{band-sub-exp}).

For the proof of (\ref{holder-sub-exp}) we shall use the representation
$\varphi(\delta \sqrt{L})=\varphi(\delta \sqrt{L})e^{\delta^2 L}e^{-\delta^2 L}$.
Let $h(\lambda):=\varphi(\lambda)e^{\lambda^2}$.
Rough calculation shows that
$\|(d/d\lambda)^ke^{\lambda^2}\|_{L^\infty[-2,2]} \le (ck)^k$
and applying Leibniz rule
$\|h^{(k)}\|_\infty \le (ck)^{k(1+\eps)}$, $\forall k\in\bN$.
Now, as in the proof of (\ref{band-sub-exp})
$$
|h(\delta \sqrt{L})(x, y)|
\le \frac{c\exp\Big\{-\kap\big(\frac{\rho(x, y)}{\delta}\big)^{1-\eps}\Big\}}
{\big(|B(x, \delta)||B(y, \delta)|\big)^{1/2}}.
$$
Just as in the proof of Theorem 3.1 in \cite{CKP},
using this and the H\"{o}lder continuity of the heat kernel we obtain
whenever $\rho(x, y)\le \delta$
\begin{align*}
&|\varphi(\delta \sqrt{L})(x, y)-\varphi(\delta \sqrt{L})(x, y')|
\le \int_M |h(\delta \sqrt{L})(x, u)||p_{\delta^2}(u, y)-p_{\delta^2}(u, y')| d\mu(u)\\
& \qquad\qquad \le  \frac{c\big(\frac{\rho(y,y')}{\delta}\big)^\alpha}
{\big(|B(x, \delta)||B(y, \delta)|\big)^{1/2}}
\int_M \frac{\exp\big\{-\kap\big(\frac{\rho(x, y)}{\delta}\big)^{1-\eps}
- c\big(\frac{\rho(x, y)}{\delta}\big)^2\big\}}
{|B(u, \delta)|}d\mu(u)\\
& \qquad\qquad \le  \frac{c\big(\frac{\rho(y,y')}{\delta}\big)^\alpha
\exp\big\{-\kap\big(\frac{\rho(x, y)}{\delta}\big)^{1-\eps}\big\}}
{\big(|B(x, \delta)||B(y, \delta)|\big)^{1/2}}.
\end{align*}
Here for the last estimate we used inequality (\ref{prod-est}) below.
This confirms (\ref{holder-sub-exp}).

To show (\ref{band-sub-exp-m}) consider the function $\psi(\lambda)=\lambda^{2m}\varphi(\lambda)$.
Using the fact that
$\|\varphi^{(k)}\|_\infty \le (ck)^{k(1+\eps)}$
it is easy to see that
$\|\psi^{(k)}\|_\infty \le 2^{2m}(2m)!(2ck)^{k(1+\eps)}$, $\forall k\in \bN$ and
$\|\psi\|_\infty \le 2^{2m}$.
Also, it is easy to see that $\psi^{(2\nu+1)}(0)=0$ for all $\nu\ge 0$.
Then just as above it follows that $|\psi(\delta \sqrt{L})(x, y)|$ satisfies (\ref{band-sub-exp})
with a slightly bigger constant on the right multiplied in addition by $2^{2m}(2m)!$.
On~the other hand, $\psi(\delta\sqrt{L})=\delta^{2m}L^m\varphi(\delta\sqrt{L})$
and (\ref{band-sub-exp-m}) follows.
$\qed$


\begin{remark}\label{rem:sharpness}
{\rm
As shown in \cite{IPX1}, in general, estimate (\ref{band-sub-exp}) is no longer
valid with $\eps=0$ for an admissible cut-off function $\varphi$ no matter what the selection of the constants
$c_1, \kap>0$ may be.
}
\end{remark}

\subsection{The algebra of operators with sub-exponentially localized kernels}

\begin{definition}
We denote by $\cL(\bbb, \kap)$ with $0<\bbb<1$ and $\kap>0$
the set of all operators of the form $f(\delta \sqrt{L})$, where
$f:\R \to \bC$ is such that the operator $f(\delta \sqrt{L})$ is an integral operator
with kernel $f(\delta \sqrt{L})(x, y)$ obeying
\begin{equation}\label{sub-exp-op}
|f(\delta \sqrt{L})(x, y)|
\le \frac{C\exp\big\{-\kap\big(\frac{\rho(x, y)}{\delta}\big)^\bbb\big\}}
{\big(|B(x, \delta)||B(y, \delta)|\big)^{1/2}},
\quad x, y\in M, \; \delta>0,
\end{equation}
for some constant $C>0$.
We introduce the norm $\|f(\delta \sqrt{L})\|_\star :=\inf C$
on $\cL(\bbb, \kap)$.
\end{definition}
We shall use the abbreviated notation
\begin{equation}\label{ddef-E}
E_{\delta, \kap}(x, y)
:=\frac{\exp\big\{-\kap\big(\frac{\rho(x, y)}{\delta}\big)^\bbb\big\}}
{\big(|B(x, \delta)||B(y, \delta)|\big)^{1/2}}.
\end{equation}

It will be critical for our development of frames to show that the class
$\cL(\bbb, \kap)$ is an algebra:


\begin{theorem}\label{thm:algebra}
$(a)$
If the operators $f_1(\delta \sqrt{L})$ and $f_2(\delta \sqrt{L})$ belong to $\cL(\bbb, \kap)$, i.e.
\begin{equation}\label{loc-f-j}
|f_j(\delta \sqrt{L})(x, y)| \le c_jE_{\delta, \kap}(x, y), \quad j=1, 2,
\end{equation}
then the operator $f_1(\delta \sqrt{L})f_2(\delta \sqrt{L})$ also belongs to $\cL(\bbb, \kap)$
and
\begin{equation}\label{loc-ff}
|f_1(\delta \sqrt{L})f_2(\delta \sqrt{L})(x, y)| \le \cn c_1c_2E_{\delta, \kap}(x, y),
\end{equation}
for some constant $\cn>1$ depending only on $\bbb, \kap, c_0$.

$(b)$
There exists a constant $\eps>0$ depending only on $\bbb, \kap, d$ such that
if the operator $f(\delta \sqrt{L})$ is in $\cL(\bbb, \kap)$
and $\|f(\delta \sqrt{L})\|_\star < \eps$, then $\Id-f(\delta \sqrt{L})$ is invertible
and $[\Id-f(\delta \sqrt{L})]^{-1}- \Id$ belongs to $\cL(\bbb, \kap)$.

\end{theorem}


\noindent
{\bf Proof.}
Clearly, to prove Part (a) of the theorem it suffices to show that there exists a constant
$c_\natural>0$, depending only on $\bbb, \kap, c_0$, such that
\begin{align}\label{prod-est}
& \int_M |B(u, \delta)|^{-1} \exp\Big\{-\kap\Big(\frac{\rho(x, u)}{\delta}\Big)^{\bbb}
-\kap\Big(\frac{\rho(u, y)}{\delta}\Big)^{\bbb}\Big\} d\mu(u)\notag\\
&\qquad\qquad\qquad\qquad\qquad\qquad\qquad\qquad\qquad\qquad
\le c_\natural\exp\Big\{-\kap\Big(\frac{\rho(x, y)}{\delta}\Big)^{\bbb}\Big\}.
\end{align}
The proof of this relies on the following inequality:
For any $x, y, u\in M$
\begin{equation}\label{ineq}
\rho(x, u)^\bbb + \rho(y, u)^\bbb \ge \rho(x, y)^\bbb +(2-2^\bbb)\rho(x, u)^\bbb
\quad\hbox{if}\quad \rho(x, u) \le \rho(y, u).
\end{equation}
To prove this inequality, suppose $\rho(x, u) \le \rho(y, u)$ and
let $\rho(y, u)= t\rho(x, u)$, $t\ge 1$.
Then using that $0<\bbb<1$
\begin{align*}
\rho(x, u)^\bbb + \rho(y, u)^\bbb
&= (1+t^\bbb)\rho(x, u)^\bbb\\
&= [(1+t)\rho(x, u)]^\bbb + [1+t^\bbb-(1+t)^\bbb]\rho(x, u)^\bbb\\
&\ge [\rho(x, u)+\rho(y, u)]^\bbb + \min_{t\ge 1}[1+t^\bbb-(1+t)^\bbb]\rho(x, u)^\bbb\\
&\ge \rho(x, y)^\bbb + (2-2^\bbb)\rho(x, u)^\bbb,
\end{align*}
which confirms (\ref{ineq}).

Let $x, y\in M$, $x\ne y$.
We split $M$ into two:
$M':=\{u\in M: \rho(x, u) \le \rho(y, u)\}$ and $M'':= M\setminus M'$.
Denote $I':=\int_{M'} \cdots$ and $I'':= \int_{M''}\cdots$.
To estimate $I'$ we use inequality (\ref{ineq}) and obtain
\begin{align*}
I' &\le \exp\Big\{-\kap\Big(\frac{\rho(x, y)}{\delta}\Big)^{\bbb}\Big\}
\int_M |B(u, \delta)|^{-1} \exp\Big\{-\kap(2-2^\bbb)\Big(\frac{\rho(x, u)}{\delta}\Big)^{\bbb}\Big\} d\mu(u)\\
&\le c\exp\Big\{-\kap\Big(\frac{\rho(x, y)}{\delta}\Big)^{\bbb}\Big\}
\int_M |B(u, \delta)|^{-1} (1+\delta^{-1}\rho(x, u))^{-2d-1} d\mu(u)\\
&\le c\exp\Big\{-\kap\Big(\frac{\rho(x, y)}{\delta}\Big)^{\bbb}\Big\},
\end{align*}
where $c>0$ is a constant depending on $\bbb, \kap, c_0$.
Because of the symmetry the same estimate holds for $I''$ and the proof of (a) is complete.

Part (b) follows immediately from (a).
$\qed$

We shall also need a discrete version of inequality (\ref{prod-est}):


\begin{lemma}\label{lem:discr-prod}
Suppose $\cX$ is a maximal $\delta$-net on $M$
and $\{A_\xi\}_{\xi\in\cX}$ is a~companion disjoint partition of $M$ as in $\S\ref{sec:max-d-nets}$.
Let $\dst \ge \delta$.
Then
\begin{equation}\label{discr-prod}
\sum_{\xi\in \cX} \frac{|A_\xi|}{|B(\xi, \dst)|}
\exp\Big\{-\kappa\Big(\frac{\rho(x, \xi)}{\dst}\Big)^{\bbb}
-\kappa\Big(\frac{\rho(y, \xi)}{\dst}\Big)^{\bbb}\Big\}
\le \cs\exp\Big\{-\kappa\Big(\frac{\rho(x, y)}{\dst}\Big)^{\bbb}\Big\},
\end{equation}
where $\cs>1$ depends only on $\bbb, \kappa, c_0$.
\end{lemma}

\noindent
{\bf Proof.}
We proceed similarly as above.
Let $x, y\in M$, $x\ne y$.
We split $\cX$ into two sets:
$\cX':=\{\xi\in \cX: \rho(x, \xi) \le \rho(y, \xi)\}$ and $\cX'':= \cX\setminus \cX'$.
Set $\Sigma':=\sum_{\xi\in\cX'} \cdots$ and $\Sigma'':=\sum_{\xi\in\cX''} \cdots$.
Now, using inequality (\ref{ineq}) we get
\begin{align*}
\Sigma'
&\le \exp\Big\{-\kap\Big(\frac{\rho(x, y)}{\dst}\Big)^\bbb\Big\}
\sum_{\xi\in\cX'}\frac{|A_\xi|}{|B(\xi, \dst)|}
\exp\Big\{-\kap(2-2^\bbb)\Big(\frac{\rho(x, \xi)}{\dst}\Big)^\bbb\Big\}\\
&\le c\exp\Big\{-\kap\Big(\frac{\rho(x, y)}{\dst}\Big)^\bbb\Big\}
\sum_{\xi\in \cX} \frac{|A_\xi|}{|B(\xi, \dst)|}\big(1+\dst^{-1}\rho(x, \xi)\big)^{-2d-1}\\
&\le c\exp\Big\{-\kap\Big(\frac{\rho(x, y)}{\dst}\Big)^\bbb\Big\},
\end{align*}
where in the last inequality we used estimate (\ref{basic-est}).
By the same token, the same estimate holds for $\Sigma''$.
$\qed$

\subsection{Spectral spaces}\label{spectral-spaces}
As elsewhere we adhere to the setting described in the introduction.
We let $E_\lambda$, $\lambda\ge 0$, be the spectral resolution associated with
the self-adjoint positive operator $L$ on $\LL^2:=L^2(M, d\mu)$.
Further, we let $F_\lambda$, $\lambda\ge 0$, denote the spectral resolution associated with $\sqrt L$,
i.e. $F_\lambda=E_{\lambda^2}$.
As in \S\ref{sec:func-calculus} we are interested in operators of the form $f(\sqrt L)$.
Then $f(\sqrt L)= \int_0^\infty f(\lambda) d F_\lambda$
and the spectral projectors are defined by
$E_\lambda = \ONE_{[0, \lambda]}(L) := \int_0^\infty \ONE_{[0,\lambda]}(u) dE_u$
and
\begin{equation}\label{spect-projector}
F_\lambda = \ONE_{[0,  \lambda]}(\sqrt L)
:=\int_0^\infty  \ONE_{[0,\lambda]}(u)dF_u
=\int_0^\infty  \ONE_{[0,\lambda]}(\sqrt u)dE_u.
\end{equation}
Recall the definition of the spectral spaces $\Sigma^p_\lambda$, $1\le p \le \infty$, from \cite{CKP}:
$$
\Sigma^p_\lambda
:= \{f\in \LL^p: \theta (\sqrt L)f = f \hbox{ for all }
\theta \in C^\infty_0(\bR_+), \;  \theta \equiv 1 \hbox{ on } \;  [0, \lambda]\}
$$
and for any compact $K \subset [0, \infty)$
$$
\Sigma^p_K
:= \{f\in \LL^p: \theta (\sqrt L)f = f \hbox{ for all }
\theta \in C^\infty_0(\bR_+), \;  \theta \equiv 1 \hbox{ on } \;  K\}.
$$
We now extend this definition: Given a space $Y$ of measurable functions on $M$
$$
\Sigma_\lambda = \Sigma_\lambda(Y)
:= \{f\in Y: \theta (\sqrt L)f = f \hbox{ for all }
\theta \in C^\infty_0(\bR_+), \;  \theta \equiv 1 \hbox{ on } \;  [0, \lambda]\}.
$$
The space $Y$ usually will be obvious from the context and will not be mentioned explicitly.

We next relate different weighted $L^p$-norms of spectral functions.


\begin{proposition}\label{prop:Nikolski}
Let $0 < p \le q \le \infty$ and $\gg\in\bR$.
Then there exists a constant $c>0$ such that
\begin{equation}\label{norm-relation2}
\| |B(\cdot, \lambda^{-1})|^\gg g(\cdot)\|_q
\le c\| |B(\cdot, \lambda^{-1})|^{\gg+1/q-1/p} g(\cdot)\|_p
\quad\hbox{for}\quad g \in \Sigma_\lambda, \; \lambda \ge 1.
\end{equation}
Therefore, assuming in addition the non-collapsing condition $(\ref{non-collapsing})$ we have
$\Sigma_\lambda^p \subset \Sigma_\lambda^q$ and

\begin{equation}\label{norm-relation1}
\|g\|_q \le c\lambda^{d(1/p-1/q)}\|g\|_p,
\quad g \in \Sigma_\lambda^p, \; \lambda  \ge 1.
\end{equation}
\end{proposition}

\noindent
{\bf Proof.}
Let $g \in \Sigma_\lambda$, $\lambda\ge 1$, and set $\delta:= \lambda^{-1}$.
Let $\theta \in C^\infty_0(\bR_+)$ be so that $\theta \equiv 1$ on $[0,1]$.
Denote briefly $H(x, y):=\theta(\delta\sqrt L)(x, y)$
the kernel of the operator $\theta(\delta\sqrt L)$.
By Theorem \ref{thm:main-local-kernels} it obeys
\begin{equation}\label{est-ker-theta}
|H(x, y)| \le c_\sigma D_{\delta, \sigma+d/2}(x, y)
\le c_\sigma'|B(x, \delta)|^{-1}\Big(1+\frac{\rho(x, y)}{\delta}\Big)^{-\sigma}
\quad \forall \sigma >0.
\end{equation}

Suppose $1<p<\infty$.
Clearly,
$
g(x)=\theta(\delta\sqrt L)g(x)=\int_M H(x, y)g(y) d\mu(y)
$
and using (\ref{est-ker-theta}) with $\sigma\ge dp'(|\gg|+1/p)+d+1$ (here $1/p+1/p'=1$),
H\"{o}lder's inequality, and (\ref{D2})
we obtain
\begin{align*}
|g(x)|
&\le \| |B(\cdot, \delta)|^{\gg-1/p} g(\cdot)\|_p
\Big(
\int_M \Big(|H(x, y)| |B(y, \delta)|^{-\gg+1/p}\Big)^{p'} d\mu(y)
 \Big)^{1/p'}\\
 & \le c\| |B(\cdot, \delta)|^{\gg-1/p} g(\cdot)\|_p
\Big(
\int_M \frac{|B(x, \delta)|^{(-\gg+1/p-1)p'}}
{\big(1+\frac{\rho(x,y)}{\delta}\big)^{s}} d\mu(y)\Big)^{1/p'}\\
 &\le c\| |B(\cdot, \delta)|^{\gg-1/p} g(\cdot)\|_p
|B(x, \delta)|^{-\gg}.
\end{align*}
Here $s:=\sigma - dp'(|\gg|+1/p) \ge d+1$ and for the latter inequality we used (\ref{tech1}).
Therefore,
\begin{equation}\label{norm-rel-infty}
\| |B(\cdot, \delta)|^{\gg} g(\cdot)\|_\infty
\le c\| |B(\cdot, \delta)|^{\gg-1/p} g(\cdot)\|_p,
\quad 1<p\le \infty.
\end{equation}
Thus (\ref{norm-relation2}) holds in the case $q=\infty$.

Let now $0<p\le 1$. Then we use estimate (\ref{norm-rel-infty}) with $p=2$ to obtain
\begin{align*}
\| |B(\cdot, \delta)|^{\gg} g(\cdot)\|_\infty
&\le c\| |B(\cdot, \delta)|^{\gg-1/2} g(\cdot)\|_2\\
& =c\Big(\int_M \big[|B(x, \delta)|^{\gg} |g(x)|\big]^{2-p}
\big[|B(x, \delta)|^{\gg-1/p} |g(x)|\big]^p d\mu(x)\Big)^{1/2}\\
&\le c\| |B(\cdot, \delta)|^{\gg} g(\cdot)\|_\infty^{1-p/2}
\| |B(\cdot, \delta)|^{\gg-1/p} g(\cdot)\|_p^{p/2},
\end{align*}
which yields the validity of (\ref{norm-rel-infty}) for $0<p\le \infty$.

Finally, we derive  (\ref{norm-relation2}) in the case $0<p<q<\infty$ from (\ref{norm-rel-infty})
(with $\gg$ replaced by $\gg+1/q$)
as follows
\begin{align*}
\| |B(\cdot, \delta)|^{\gg} g(\cdot)\|_q
&= \Big(\int_M \big[|B(x, \delta)|^{\gg+\frac{1}{q}} |g(x)|\big]^{q-p}
\big[|B(x, \delta)|^{\gg+\frac{1}{q}-\frac{1}{p}} |g(x)|\big]^p d\mu(x)\Big)^{1/q}\\
&\le c\| |B(\cdot, \delta)|^{\gg+\frac{1}{q}} g(\cdot)\|_\infty^{1-\frac{p}{q}}
\Big(\int_M ||B(x, \delta)|^{\gg+\frac{1}{q}-\frac{1}{p}} g(x)|^p d\mu(x)\Big)^{\frac{1}{q}}\\
&\le c\| |B(\cdot, \delta)|^{\gg+1/q-1/p} g(\cdot)\|_p.
\end{align*}
The proof of (\ref{norm-relation2}) is complete.

The non-collapsing condition (\ref{non-collapsing}) yields (\ref{est-muB-unify1}),
which along with (\ref{norm-relation2}) leads to (\ref{norm-relation1}).
$\qed$

\subsection{Kernel norms}\label{sec:kernel-norms}

Bounds on the $\LL^p$-norms of the kernels of operators of the form
$\theta (\delta \sqrt L)$ are developed in \S 3.6 in \cite{CKP} and play an important role
in the development of frames.
We present them next in the form we need them.


\begin{theorem}\label{thm:norms}\cite{CKP}
Assume that the reverse doubling condition $(\ref{reverse-doubling})$ is valid,
and let $\theta\in C^\infty(\bR_+)$, $\theta \ge 0$,
$\supp \theta\subset [0, R]$ for some $R>1$,
and $\theta^{(2\nu+1)}(0)=0$, $\nu=0, 1, \dots$.
Suppose that either

$(i)$\;\;
$\theta(u)\ge 1$ for $u\in [0, 1]$, or

$(ii)$
$\theta(u)\ge 1$ for $u\in [1, b]$, where $b>1$ is a sufficiently large constant.

\noindent
Then for $0 < p \le \infty$, $0<\delta \le \min\{1, \frac{\diam M}{3}\}$, and $x \in M$ we have
\begin{equation}\label{est-norm-1}
c_1|B(x,\delta)|^{1/p-1}
\le \| \theta (\delta \sqrt L) (x,.)\|_p
\le c_2|B(x,\delta)|^{1/p-1},
\end{equation}
where $c_1, c_2>0$ are independent of $x, \delta$.
\end{theorem}

{\bf The constant \boldmath $b>1$} that appears in the above theorem will play a distinctive
role in what follows.

\section{Construction of frames}\label{sec:frames}
\setcounter{equation}{0}

Our goal here is to construct a pair of dual frames whose elements
are band limited and have sub-exponential space localization.
This is a major step forward compared with the frames from \cite{CKP},
where the elements of the second (dual) frame have limited space localization.
We shall utilize the main idea of the construction in \cite{CKP}
and also adopt most of the notation from \cite{CKP}.

We shall first provide the main ingredients for this construction and
then describe the two main steps of our scheme:
(i) Construction of Frame~\#~1, and
(ii) Construction of a nonstandard dual Frame \#~2.

\subsection{Sampling theorem and cubature formula}\label{sec:sampling-cubature}

The main vehicle in constructing frames is a sampling theorem for $\Sigma_\lambda^2$ and
a cubature formula for $\Sigma_\lambda^1$.
Their realization relies on the nearly exponential localization of operator kernels
induced by smooth cut-off functions $\varphi$ (Theorem~\ref{thm:main-local-kernels}):
If $\varphi \in C^\infty_0(\bR_+)$,
 $\supp\varphi \subset [0,b]$, $b>1$, $0\le \varphi \le 1$, and $\varphi = 1$ on $[0,1]$,
then there exists a constant $\alpha>0$ such that
for any $\delta >0$ and $x,y,x' \in M$
\begin{align}
|\varphi(\delta \sqrt{L})(x,y) | &\le \Css D_{\delta, \sigma}(x,y)
\quad\hbox{and}\label{local-Phi}\\
|\varphi(\delta \sqrt{L})(x,y)- \varphi(\delta \sqrt{L})(x',y) |
&\leq \Css \Big(\frac{\rho(x,x')}{\delta}\Big)^\alpha D_{\delta, \sigma}(x,y),
\; \rho(x,x')\le \delta. \label{Lip-Phi}
\end{align}
Here $\Css>1$ is a constant depending on $\varphi$, $\sigma$ and the other parameters,
but independent of $x, y, x'$ and $\delta$.

The above allows to establish a {\bf Marcinkiewicz-Zygmund inequality} for $\Sigma_\lambda^1$
\cite[Proposition 4.1]{CKP}:
{\em
Given $\lambda \ge 1$, let $\cX_\delta$ be a maximal $\delta-$net on $M$ with
$\delta :=\gamma \lambda^{-1}$,
where $0<\gamma < 1$,
and
suppose $\{A_\xi\}_{\xi\in\cX_\delta}$ is a companion disjoint partition of $M$
as described in \S\ref{sec:max-d-nets}.
Then for any $f \in \Sigma_\lambda^p$, $1\le p <\infty$,
\begin{equation}\label{Marcink1}
\sum_{\xi \in \cX_\delta} \int _{A_\xi}    |f(x) - f(\xi)|^p dx
\leq  [\CPhi  \gamma^{\alpha } \Cdiam]^p \| f \|_p^p,
\end{equation}
and a similar estimate holds when $p=\infty$.
Here $\CPhi$ is the constant from $(\ref{local-Phi})-(\ref{Lip-Phi})$
with $\sigma_*:=2d +1$ and
$\Cdiam>1$ depends only on $c_0, \Cstar, \cstar$ from $(\ref{doubling-0})-(\ref{lip})$.
}


The needed {\bf sampling theorem} takes the form \cite[Theorem 4.2]{CKP}:
{\em
Given a constant $0 < \eps <1$, let $0<\gamma<1$ be so that
$
\CPhi  \gamma^{\alpha }\Cdiam \le \eps/3.
$
Suppose $\cX_\delta$ is a~maximal $\delta-$net on $M$ and
$\{A_\xi\}_{\xi\in\cX_\delta}$ is a companion disjoint partition of $M$
with $\delta :=\gamma \lambda^{-1}$. %
Then for any $f\in \Sigma_\lambda^2$
\begin{equation}\label{samp3}
(1-\eps)\| f \|_2^2 \le \sum_{\xi \in \cX_\delta} |A_\xi|| f(\xi)|^2 \le (1+\eps)\| f \|_2^2.
\end{equation}
}

The Marcinkiewicz-Zygmund inequality (\ref{Marcink1}) is also used for the construction of
a {\bf cubature formula} \cite[Theorem 4.4]{CKP}:
{\em
Let $0< \gamma < 1$ be selected so that
$
\CPhi  \gamma^{\alpha }\Cdiam = \frac 14.
$
Given $\lambda \ge 1$, suppose $\cX_\delta$ is a maximal $\delta$-net on $M$
with $\delta :=\gamma \lambda^{-1}$.
Then there exist positive constants $($weights$)$
$\{\ww^\lambda_\xi\}_{\xi \in \cX_\delta}$ such that
\begin{equation}\label{quadrature}
\int_M f(x) d\mu(x) = \sum_{\xi \in \cX_\delta} \ww^\lambda_\xi f(\xi)
\quad \forall f \in \Sigma^1_\lambda,
\end{equation}
and
$(2/3)|B(\xi, \delta/2)| \le \ww^\lambda_\xi \le 2|B(\xi, \delta)|$,
$\xi \in \cX_\delta$.
}

\subsection{Construction of Frame \# 1}\label{natural-frame}

We begin with the construction of a well-localized frame based on the kernels of spectral operators
considered in \S\ref{sec:local-kernels}.

We use Theorem~\ref{thm:band-sub-exp} to construct a cut-off function $\Phi$ with the following properties:
$\Phi\in C^\infty (\bR_+)$,
$\Phi(u)=1$ for $u\in[0, 1]$,
$0\le \Phi \le 1$, and
$\supp \Phi \subset [0, b]$, where $b > 1$ is the constant from Theorem~\ref{thm:norms}.

Set $\Psi(u):=\Phi(u)-\Phi(bu)$.
Clearly, $0\le \Psi \le 1$ and $\supp \Psi\subset [b^{-1}, b]$.
We also assume that $\Phi$ is selected so that
$\Psi(u) \ge c>0$ for $u\in [b^{-3/4}, b^{3/4}]$.

From Theorem~\ref{thm:band-sub-exp} it follows that $\Phi(\delta\sqrt{L})$ and $\Psi(\delta\sqrt{L})$
are integral operators whose kernels
$\Phi(\delta\sqrt{L})(x, y)$ and $\Psi(\delta\sqrt{L})(x, y)$
have sub-exponential localization, namely,
\begin{equation}\label{sub-exp-Phi}
|\Phi(\delta \sqrt{L})(x, y)|, |\Psi(\delta\sqrt{L})(x, y)|
\le \cd E_{\delta, \kap}(x, y),
\quad x, y\in M,
\end{equation}
with
\begin{equation}\label{def-E}
E_{\delta, \kappa}(x, y)
:= \frac{\exp\big\{-\kappa\big(\frac{\rho(x, y)}{\delta}\big)^\bb\big\}}
{\big(|B(x, \delta)||B(y, \delta)|)^{1/2}}.
\end{equation}
Here $0<\bb<1$ is an arbitrary constant (as close to $1$ as we wish),
and $\kappa>0$ and $\cd >1$ are constants depending only on
$\bb, b$ and the constants $c_0, \Cstar, \cstar$ from $(\ref{doubling-0})-(\ref{lip})$.
Also, $\Phi(\delta \sqrt{L})(x, y)$ and $\Psi(\delta\sqrt{L})(x, y)$ are H\"{o}lder continuous,
namely,
\begin{equation}\label{holder-Phi}
|\Phi(\delta \sqrt{L})(x, y)-\Phi(\delta \sqrt{L})(x, y')|
\le c \big(\delta^{-1}\rho(y, y')\big)^\alpha E_{\delta, \kap}(x, y)
\;\hbox{ if }\; \rho(y, y')\le \delta,
\end{equation}
and the same holds for $\Psi(\delta\sqrt{L})(x, y)$.
Furthermore, for any $m\ge 1$
\begin{equation}\label{m-sub-exp-Phi}
|L^m\Phi(\delta \sqrt{L})(x, y)|, |L^m\Psi(\delta \sqrt{L})(x, y)|
\le c_m\delta^{-2m} E_{\delta, \kap}(x, y),
\quad x, y\in M.
\end{equation}
We shall regard $\bb$ and $\kap$ as {\em parameters} of our frames and they will be fixed from now on.

Set
\begin{equation}\label{def-Psi-j}
\Psi_0(u):=\Phi(u) \quad\hbox{ and }\quad \Psi_j(u):=\Psi(b^{-j}u),\;\; j\ge 1.
\end{equation}
Clearly, $\Psi_j\in C^\infty (\bR_+)$, $0\le \Psi_j \le 1$, $\supp \Psi_0 \subset [0, b]$,
$\supp \Psi_j \subset [b^{j-1}, b^{j+1}]$, $j\ge 1$, and
$\sum_{j\ge 0}\Psi_j(u) = 1$ for $u\in\bR_+$.
By Corollary~3.9 in \cite{CKP} (see also Proposition~\ref{prop:decomp-DD2} below)
we have the following Littlewood-Paley decomposition
\begin{equation}\label{repres-Psi-j}
f= \sum_{j\ge 0} \Psi_j(\sqrt L)f
\quad\hbox{for $f\in \LL^p$, $\;1\le p\le\infty$. \;$(\LL^\infty:=\UCB)$}
\end{equation}
From above it follows that
\begin{equation}\label{sum-Psi-j}
\frac 12 \le \sum_{j\ge 0} \Psi_j^2(u) \leq 1, \quad u\in \R_+.
\end{equation}
As
$
\|\Psi_j(\sqrt L) f\|_2^2
= \langle \Psi_j(\sqrt L) f, \Psi_j(\sqrt L) f \rangle
= \langle \Psi_j^2(\sqrt L) f,  f \rangle,
$
we obtain
$$
\sum_{j \ge 0} \|\Psi_j(\sqrt L) f\|_2^2
= \int_0^\infty  \sum_{j\ge 0} \Psi_j^2(u) d \langle F_u f,f\rangle,
$$
and using (\ref{sum-Psi-j}) we get
\begin{equation}\label{frame1}
\frac 12  \| f\|_2^2 \le \sum_{j \ge 0} \|\Psi_j(\sqrt L) f \|_2^2 \le \| f\|_2^2,
\quad f\in \LL^2.
\end{equation}

At this point we introduce a constant $0 <\eps <1$ by
\begin{equation}\label{pick-eps}
\eps:= (8\cs\cn^2\cd^2)^{-1},
\end{equation}
where the constant
$\cs>1$ is from Lemma~\ref{lem:discr-prod},
$\cn >1$ is from Theorem~\ref{thm:algebra}, and
$\cd>1$ is from (\ref{sub-exp-Phi}).
Pick $0< \gamma <1$ so that
\begin{equation}\label{pick-gamma}
\CPhi \gamma^{\alpha }\Cdiam = \eps/3,
\end{equation}
where $\CPhi$ is the constant from (\ref{local-Phi})-(\ref{Lip-Phi}) with $\sig_*:=2d+1$
and $\Cdiam>1$ is from (\ref{Marcink1}).

For any $j\ge 0$ let $\XX_j \subset M$ be a maximal $\ddj-$net on $M$
with
$\ddj:=\gamma b^{-j-2}$ and
suppose $\{A_\xi^j\}_{\xi\in\XX_j}$ is a companion disjoint partition of $M$
consisting of measurable sets such that
$B(\xi, \ddj/2) \subset  A_\xi^j \subset B(\xi,\ddj)$, $\xi \in \XX_j$,
as in \S\ref{sec:max-d-nets}.
By the sampling theorem (\S\ref{sec:sampling-cubature}) and the definition of $\Psi_j$
it follows that
\begin{equation}\label{sampling-L2}
(1-\eps)\| f \|_2^2 \le \sum_{\xi \in \XX_j} |A_\xi^j|| f(\xi)|^2 \le (1+\eps)\| f \|_2^2
\quad\hbox{for}\quad f\in \Sigma_{b^{j+2}}^2.
\end{equation}
From the definition of $\Psi_j$ we have
$\Psi_j(\sqrt L) f \in \Sigma^2_{b^{j+1}}$ for $f\in \LL^2$,
and hence (\ref{frame1}) and (\ref{sampling-L2}) yield
\begin{equation}\label{frame2}
\frac 14 \|f\|_2^2 \le
\sum_{j \geq 0} \sum_{\xi \in \XX_j} |A^j_\xi || \Psi_j(\sqrt L) f (\xi) |^2
\leq  2\| f\|_2^2,
\quad f\in L^2.
\end{equation}
Observe that
\begin{align*}
\Psi_j(\sqrt L) f (\xi)
&= \int_M f(u) \Psi_j(\sqrt L)(\xi,u)d\mu(u)\\
&= \int_M f(u) \overline{ \Psi_j(\sqrt L)(u,\xi )} d\mu(u)
= \big\langle f,  \Psi_j(\sqrt L)( . ,\xi  ) \big\rangle.
\end{align*}
We define the system $\{\psi_{\xi}\}$ by
\begin{equation}\label{def-frame}
\psi_{\xi} (x):= |A^{j}_\xi |^{1/2} \Psi_j(\sqrt L)(x,\xi),
\quad \xi \in \XX_j, j\ge 0.
\end{equation}
Write
$\cX:=\cup_{j\ge 0}\cX_j$,
where equal points from different sets $\cX_j$ will be regarded as distinct elements of $\cX$,
so $\cX$ can be used as an index set.
From the above observation and (\ref{frame2}) it follows that
$\{\psi_{\xi}\}_{\xi\in\cX}$
is a frame for $\LL^2$.

We next record the main properties of this system.


\begin{proposition}\label{prop:frame-prop}
$(a)$ {\rm Localization:} For any $0<\kaph<\kap$ there exist a constant $\ch>0$
such that for any $\xi\in\XX_j$, $j\ge 0$,
\begin{equation}\label{prop-psi-1}
|\psi_{\xi} (x)|
\le \ch |B(\xi, b^{-j})|^{-1/2} \exp\big\{-\kaph(b^j\rho(x, \xi))^\bbb\big\}
\end{equation}
and for any $m\ge 1$
\begin{equation}\label{prop-psi-11}
|L^m\psi_{\xi} (x)|
\le c_m|B(\xi, b^{-j})|^{-1/2}b^{2jm} \exp\big\{-\kaph(b^j\rho(x, \xi))^\bbb\big\}.
\end{equation}
Also, if $\rho(x, y) \le b^{-j}$
\begin{equation}\label{prop-psi-Lip}
|\psi_{\xi} (x)- \psi_{\xi} (y)|
\le \ch |B(\xi, b^{-j})|^{-1/2}(b^j\rho(x, y))^\alpha \exp\big\{-\kaph(b^j\rho(x, \xi))^\bbb\big\},
\quad \alpha>0.
\end{equation}

$(b)$ {\rm Norms:} If in addition the reverse doubling condition $(\ref{reverse-doubling})$ is valid, then
\begin{equation}\label{prop-psi-2}
\|\psi_{\xi}\|_p \sim |B(\xi, b^{-j})|^{\frac 1p-\frac 12},
\quad 0< p \leq \infty.
\end{equation}
$(c)$ {\rm Spectral localization:}
$\psi_{\xi}\in \Sigma_b^p$ if $\xi\in \XX_0$ and
$\psi_{\xi}\in \Sigma_{[b^{j-1}, b^{j+1}]}^p$
if $\xi\in \XX_j$, $j\ge 1$, $0<p\le\infty$.

$(d)$
The system $\{\psi_{\xi}\}$ is a frame for $\LL^2$, namely, 
\begin{equation}\label{frame3}
 4^{-1}\|f\|_2^2 \le
\sum_{j \geq 0} \sum_{\xi \in \XX_j} |\langle f, \psi_{\xi}\rangle|^2
\leq  2\| f\|_2^2,
\quad \forall f\in \LL^2.
\end{equation}

\end{proposition}

\noindent
{\bf Proof.}
From (\ref{sub-exp-Phi}) and the inequality
$|B(x, b^{-j})| \le c_0(1+b^j\rho(\xi, x))^d |B(\xi, b^{-j})|$, see (\ref{D2}),
we derive for $\xi\in\cX_j$
\begin{align*}
|\psi_{\xi} (x)|
&\le c|B(x, b^{-j})|^{-1/2}\exp\big\{-\kap\big(b^j\rho(x,y)\big)^\bb\big\}\\
&\le \ch|B(\xi, b^{-j})|^{-1/2}\exp\big\{-\kaph\big(b^j\rho(x,y)\big)^\bb\big\},
\end{align*}
which confirms (\ref{prop-psi-1}).
Estimate (\ref{prop-psi-11}) follows in the same way from (\ref{m-sub-exp-Phi})
and (\ref{prop-psi-Lip}) follows from (\ref{holder-Phi});
(\ref{prop-psi-2}) follows by Theorem~\ref{thm:norms}.
The spectral localization is obvious by the definition.
Estimates (\ref{frame3}) follow by (\ref{frame2}).
$\qed$

\subsection{Construction of Frame \# 2}\label{dual-frame}

Here the cardinal problem is to construct a dual frame to $\{\psi_{\xi}\}$ with similar
space and spectral localization.

The first step in this construction is to introduce two new cut-off functions
by dilating
$\Psi_0$ and $\Psi_1$ from \S\ref{natural-frame}:
\begin{equation}\label{def-LLam}
\LLam_0(u) := \Phi (b^{-1}u)\quad \hbox{and} \quad
\LLam_1(u) := \Phi (b^{-2}u)- \Phi(bu)
= \LLam_0(b^{-1}u) - \LLam_0(b^2u).
\end{equation}
Clearly,
$\supp\LLam_0 \subset [0, b^2]$, $\LLam_0(u)=1$ for $u\in [0, b]$, 
$\supp\LLam_1 \subset [b^{-1}, b^3]$, $\LLam_1(u)=1$ for $u\in [1, b^2]$,
$0\le \LLam_0, \LLam_1 \le 1$,
and
\begin{equation}\label{prop-LLam}
\LLam_0(u)\Psi_0(u) = \Psi_0(u), \quad
\LLam_1(u)\Psi_1(u) = \Psi_1(u).
\end{equation}
We shall also need the 
cut-off function
$\TTheta(u):=\Phi(b^{-3}u)$.
Note that
$\supp \TTheta \subset [0, b^4]$,
$\TTheta(u)=1$ for $u\in[0, b^3]$, and $\TTheta\ge 0$.
Hence,
$\TTheta(u) \LLam_j(u) = \LLam_j(u)$, $j=0,1$.


The kernels of the operators
$\Gamma_0(\delta \sqrt{L})$, $\Gamma_1(\delta \sqrt{L})$, 
and $\Theta(\delta \sqrt{L})$
inherit the localization and H\"{o}lder continuity
of $\Phi(\delta \sqrt{L})(x, y)$, see (\ref{sub-exp-Phi}) and (\ref{holder-Phi})-(\ref{m-sub-exp-Phi}).
More precisely,
if $f=\Gamma_0$ or $f=\Gamma_1$ or $f=\Theta$, then 
\begin{equation}\label{gen-local-1}
|f(\delta \sqrt{L})(x, y)| \le \cd E_{\delta, \kappa}(x, y),
\end{equation}
\begin{equation}\label{holder-f}
|f(\delta \sqrt{L})(x, y)-f(\delta \sqrt{L})(x, y')|
\le c \big(\delta^{-1}\rho(y, y')\big)^\alpha E_{\delta, \kap}(x, y)
\;\hbox{ if }\; \rho(y, y')\le \delta,
\end{equation}
and for any $m\in\bN$
\begin{equation}\label{gen-local-2}
|L^m f(\delta \sqrt{L})(x, y)| \le c_m \delta^{-2m}E_{\delta, \kappa}(x, y).
\end{equation}

The next lemma will be the main tool in constructing Frame \# 2.


\begin{lemma}\label{lem:instrument}

Given $\lambda\ge 1$, let $\XX_\dd$ be a~maximal $\dd-$net on $M$
with
$\dd:=\gamma\lambda^{-1}b^{-3}$ and
suppose $\{A_\xi\}_{\xi\in\XX_\dd}$ is a companion disjoint partition of $M$
consisting of measurable sets such that
$B(\xi, \dd/2) \subset  A_\xi \subset B(\xi,\dd)$,  $\xi \in \XX_\dd$ $(\S\ref{sec:max-d-nets})$.
%
Set $\omega_\xi := \frac{1}{1+\eps}|A_\xi |\sim  |B(\xi, \delta)|.$
Let $\LLam=\LLam_0$ or $\LLam=\Gamma_1$.
Then there exists an operator $\TT_\lambda: \LL^2\to\LL^2$
of the form $\TT_\lambda = \Id + \SSS_\lambda$ such that

$(a)$
$$
\|f \|_2\le \|\TT_\lambda f\|_2 \le \frac 1{1-2\varepsilon}\|f \|_2
\quad \forall f \in \LL^2.
$$

$(b)$
$\SSS_\lambda$ is an integral operator with kernel
$\SSS_\lambda(x, y)$ verifying
\begin{equation}\label{local-S}
|\SSS_\lambda(x,y)| \le  cE_{\lambda^{-1}, \kappa/2}(x,y),\quad x,y\in M.
\end{equation}

$(c)$
$\SSS_\lambda (\LL^2)\subset \Sigma_{\lambda b^2}^2$ if $\LLam = \LLam_0$ and
$\SSS_\lambda (\LL^2)\subset \Sigma_{[\lambda b^{-1}, \lambda b^3]}^2$ if $\LLam = \LLam_1$.

$(d)$ For any $f\in \LL^2$ such that $ \Gamma(\lambda^{-1} \sqrt L)f =f$ we have
\begin{equation}\label{instr-1}
f(x) =  \sum_{\xi \in \XX_\dd}  \omega_\xi  f(\xi)  \TT_\lambda [\LLam_\lambda(\cdot, \xi)](x),
\quad x\in M,
\end{equation}
where $\LLam_\lambda(\cdot, \cdot)$ is the kernel of the operator
$\LLam_\lambda:=\LLam(\lambda^{-1} \sqrt L)$.
\end{lemma}


\noindent
{\bf Proof.}
By the sampling theorem in \S\ref{sec:sampling-cubature} we have %
$$
(1-\eps)\|f\|^2_2
\le \sum_{\xi \in \XX_\dd} |A_\xi||f(\xi)|^2
\le (1+\eps)\| f \|^2_2
\quad\hbox{for}\; f\in \Sigma_{\lam b^3}^2,
$$
and with
$\omega_\xi := \frac{1}{1+\varepsilon}|A_\xi|$
we obtain
\begin{equation}\label{instr-2}
(1- 2\eps)\|f\|^2_2\le \sum_{\xi \in \XX_\dd}  \omega_\xi|f(\xi)|^2\le \| f \|^2_2
\quad\hbox{for}\; f\in \Sigma_{\lam b^3}^2.
\end{equation}
Write briefly $\Theta_\lambda:=\Theta(\lambda^{-1} \sqrt L)$ and
let $\Theta_\lambda(\cdot, \cdot)$ be the kernel of this operator.
Consider now the positive self-adjoint operator $U_\lambda$ with kernel
$$
U_\lambda(x,y)
= \sum_{\xi \in \XX_\dd} \omega_\xi  \Theta_\lambda(x,\xi)\Theta_\lambda(\xi,y).
$$
For $f\in\Sigma_{\lambda b^3}^2$ we have
$
\langle U_\lambda f, f  \rangle =   \sum_{\xi \in \XX_\dd}  \omega_\xi \  | f(\xi)|^2
$
and hence, using (\ref{instr-2}),
\begin{equation}\label{R1}
(1-2\varepsilon) \| f \|_2^2 \leq\langle U_\lambda f, f  \rangle  \le \|f \|_2^2
\quad \hbox{for } f\in\Sigma_{\lambda b^3}^2.
\end{equation}
Now, write $\LLam_\lambda:= \LLam(\lambda^{-1} \sqrt L)$
and let $\LLam_\lambda(x,y)$ be the kernel of this operator
(recall that $\LLam = \LLam_0$ or $\LLam = \LLam_1$).
We introduce one more self-adjoint kernel operator by
$$
R_\lambda := \LLam_\lambda(\Id - U_\lambda) \LLam_\lambda
= \LLam_\lambda^2 - \LLam_\lambda U_\lambda \LLam_\lambda.
$$
Set $V_\lambda:= \LLam_\lambda U_\lambda \LLam_\lambda$
and denote by $V_\lambda (x,y)$ its kernel.
Since
$ \Theta(u)\LLam(u) = \LLam(u)$, we have
\begin{align*}
V_\lambda (x, y)
&= \sum_{\xi \in \XX_\dd}  \omega_\xi
\int_M \int_M \LLam_\lambda(x,u)\Theta_\lambda(u,\xi)\Theta_\lambda(\xi,v)\LLam_\lambda(v,y) du dv\\
&= \sum_{\xi \in \XX_\dd}  \omega_\xi \LLam_\lambda(x,\xi)\LLam_\lambda(\xi,y).
\end{align*}
By (\ref{gen-local-1}) and Lemma~\ref{lem:discr-prod} we obtain
$$ 
|V_\lambda (x,y)|  \le \cs\cd^2 E_{\lambda^{-1}, \kappa}(x,y).
$$ 
%
Also, by (\ref{gen-local-1}) and Theorem~\ref{thm:algebra}
$$
|\Gamma^2(x, y)| \le \cn\cd^2 E_{\lambda^{-1}, \kappa}(x,y).
$$
These two estimates yield
\begin{align*}
|R_\lambda(x,y)|
\le (\cn\cd^2 + \cs\cd^2) E_{\lambda^{-1}, \kappa}(x,y)
\le 2\cn\cs\cd^2 E_{\lambda^{-1}, \kappa}(x,y).
\end{align*}
To simplify our notation we set
$\tc:=2\cn\cs\cd^2$.
Thus we have
\begin{equation}\label{est-Rxy}
|R_\lambda(x,y)| \le \tc E_{\lambda^{-1}, \kappa}(x,y).
\end{equation}

From the definition of $R_\lambda$ we derive
$$
\langle R_\lambda f, f  \rangle
= \|\LLam_\lambda f\|_2^2 -
\langle U_\lambda \LLam_\lambda f , \LLam_\lambda f \rangle
\quad \hbox{for $f\in\LL^2$.}
 $$
Since $ \LLam_\lambda (\LL^2) \subset \Sigma^2_{\lambda b^3}$, then
$\Theta_\lambda \LLam_\lambda f = \LLam_\lambda f$, and by (\ref{R1})
$$
(1-2\varepsilon) \| \LLam_\lambda f \|^2_2
\le \langle U_\lambda \LLam_\lambda f , \LLam_\lambda f \rangle
\le  \| \LLam_\lambda f \|^2_2,
\quad f\in\LL^2.
$$
Hence,
$$
0 \leq  \langle R_\lambda f, f  \rangle
\le 2\eps \| \LLam_\lambda f \|^2_2
\le 2\eps \| f \|_2^2,
\quad f\in\LL^2,
$$
where for the last inequality we used that $ \|\LLam\|_\infty \le 1$.
Therefore,
$$
\| R_\lambda\|_{2\rightarrow 2} \leq 2\eps <1
\quad \hbox{and}\quad
(1-2\eps)\|f\|_2 \le  \|( \Id-R_\lambda) f\|_2 \le \| f \|_2,
\quad f\in\LL^2.
$$

We now define
$
\TT_\lambda :=(\Id- R_\lambda)^{-1} = \Id + \sum_{k\ge 1} R_\lambda^k =: \Id + \SSS_\lambda.
$
Clearly,
\begin{equation}\label{B0}
\|f\|_2 \le  \|\TT_\lambda f \|_2 \le  \frac 1{1-2\eps}\|f \|_2
\quad  \forall f \in \LL^2.
\end{equation}
If  $\LLam_\lambda f = f$, then
$$
f= \TT_\lambda (f -R_\lambda f)
= \TT_\lambda \big( f- \LLam_\lambda f + V_\lambda f \big)
=  \TT_\lambda V_\lambda f.
$$
On the other hand, if $\LLam_\lambda f =f$, then
$
(V_\lambda f)(x) = \sum_{\xi \in \XX_\dd}  \omega_\xi  f(\xi)   \LLam_\lambda(x,\xi)
$
and hence %
\begin{equation}\label{B1}
f(x) =  \sum_{\xi \in \XX_\dd}  \omega_\xi  f(\xi)  \TT_\lambda[\LLam_\lambda(\cdot,\xi)](x).
\end{equation}
By construction %
\begin{equation}\label{B2}
\SSS_\lambda: \LL^2 \mapsto \Sigma^2_{\lambda b^3}
\quad\hbox{if $\LLam=\LLam_0$ and} \quad
\SSS_\lambda: \LL^2 \mapsto \Sigma^2_{[\lambda b^{-1}, \lambda b^3]}
\quad \hbox{if $\LLam=\LLam_1$.}
\end{equation}


It remains to establish the space localization of the kernel $\SSS_\lambda (x, y)$
of the operator~$\SSS_\lambda$.
Denoting by $R_\lambda^k(x,y)$ the kernel of $R_\lambda^k$,
we have
$$
|\SSS_\lambda (x,y) |  \le  \sum_{k\ge 1}  |R_\lambda^k(x,y)|.
$$
Evidently, $R_\lambda^k = \Theta_\lambda R_\lambda^k \Theta_\lambda$.
From this, (\ref{gen-local-1}) with $f=\Theta$, and the fact that
$\|R_\lambda\|_{2\rightarrow 2} \le 2\eps$
we obtain, applying Proposition~\ref{prop:prod-oper},
\begin{equation}\label{kernel-Rk-1}
|R_\lambda^k(x,y)|
\le  \frac{\ct\cd^2 \|R_\lambda\|^k_{2\rightarrow 2}}
{\big(|B(x, \lambda^{-1})| |B(y, \lambda^{-1})|\big)^{1/2}}
\le  \frac{(2\eps)^k\ct\cd^2}{\big(|B(x, \lambda^{-1})| |B(y, \lambda^{-1})|\big)^{1/2}}.
\end{equation}
On the other hand, applying repeatedly Theorem~\ref{thm:algebra} $k-1$ times using (\ref{est-Rxy})
we obtain
\begin{equation}\label{kernel-Rk-2}
|R_\lambda^k(x,y)|\le  \cn^{k-1}\tc^k  E_{\lambda^{-1}, \kappa}(x,y).
\end{equation}
Taking the geometric average of (\ref{kernel-Rk-1}) and (\ref{kernel-Rk-2})
($0\le a\le b$, $a\le c$ $\Longrightarrow$ $a\le \sqrt{bc}$)
we get
\begin{align*}
|R_\lambda^k(x,y)|
&\le \frac{(\ct\cd^2\cn^{-1})^{1/2} (2\eps \cn\tc)^{k/2}\exp\{-\frac{\kappa}{2}(\lambda \rho(x, y))^\beta\}}
{\big(|B(x, \lambda^{-1})| |B(y, \lambda^{-1})|\big)^{1/2}}\\
&\le  \sqrt{\ct}\cd 2^{-k/2} E_{\lambda^{-1}, \kappa/2}(x,y),
\end{align*}
where we used the notation from (\ref{def-E}) and the fact that
$2\eps \cn\tc=\frac{1}{2}$, which follows by the selection of $\eps$ in (\ref{pick-eps}).
Now, summing up we arrive at
$$
|S_\lambda(x,y)|
\le  \sqrt{\ct}\cd E_{\lambda^{-1}, \kappa/2}(x,y) \sum_{k\ge 1}2^{-k/2}
\le  3\sqrt{\ct}\cd E_{\lambda^{-1}, \kappa/2}(x,y).
$$
This completes the proof of the lemma.
$\qed$

\smallskip

We can now complete the construction of the dual frame.
We shall utilize the functions and operators introduced in \S\ref{natural-frame} and above.

Write briefly
$\LLam_{\lam_0}:= \LLam_0(\sqrt L)$ and
$\LLam_{\lam_j}:= \LLam_1(b^{-j+1}\sqrt L)$ for $j\ge 1$,
$\lam_j:=b^{-j+1}$.
Observe that since $\LLam_0(u)=1$ for $u \in [0, b]$ and $\LLam_1(u)=1$ for $u\in [1, b^2]$,
then $\LLam_{\lam_0}(\Sigma_b^2)=\Sigma_b^2$ and
$\LLam_{\lam_j}(\Sigma_{[b^{j-1}, b^{j+1}]}^2) = \Sigma_{[b^{j-1}, b^{j+1}]}^2$, $j\ge 1$.
On the other hand, clearly
$\Psi_0(\cdot, y) \in \Sigma_{b}^2$ and
$\Psi_j(\cdot, y) \in \Sigma_{[b^{j-1}, b^{j+1}]}^2$ if $j\ge 1$.
Therefore, we can apply Lemma~\ref{lem:instrument} with
$\XX_j$ and $\{A_\xi^j\}_{\xi\in\XX_j}$ from \S\ref{natural-frame}, and
$\lambda=\lambda_j=b^{j-1}$
to obtain
\begin{equation}\label{rep-Psi-j}
\Psi_j(\sqrt L)(x, y)=
\sum_{\xi \in \XX_j}  \omega_\xi  \Psi_j(\xi, y) \TT_{\lambda_j} [\LLam_{\lambda_j}(\cdot, \xi)](x),
\quad \omega_\xi=(1+\eps)^{-1}|A_\xi^j|.
\end{equation}
By (\ref{def-frame}) we have
$\psi_{\xi}(x)= |A_\xi^j|^{1/2}\Psi_j(\xi, x)$ for $\xi\in\cX_j$ and we now set
\begin{equation}\label{rep-psi-tpsi}
\tilde\psi_{\xi}(x):=c_\eps |A_\xi^j|^{1/2}\TT_{\lambda_j} [\LLam_{\lambda_j}(\cdot, \xi)](x),
\quad \xi\in\XX_j, \quad c_\eps:=(1+\eps)^{-1}.
\end{equation}
Thus $\{\tilde\psi_{\xi}\}_{\xi\in\XX}$ with $\cX:=\cup_{j\ge 0}\cX_j$, is the desired dual frame.
Note that (\ref{rep-Psi-j}) takes the form
\begin{equation}\label{rep2-Psi-j}
\Psi_j(\sqrt L)(x, y)=
\sum_{\xi \in \XX_j} \psi_{\xi}(y)\tilde\psi_{\xi}(x).
\end{equation}

We next record the main properties of the dual frame $\{\tilde\psi_{\xi}\}$.
They are similar to the properties of $\{\psi_{\xi}\}$.


\begin{theorem}\label{thm:dual-frame}
$(a)$ {\rm Representation:}
For any $f\in\LL^p$, $1\le p\le \infty$, with $\LL^\infty:=\UCB$ we have
\begin{equation}\label{rep-L2}
f = \sum_{\xi \in \XX} \langle f, \tilde\psi_{\xi}\rangle \psi_{\xi}
= \sum_{\xi \in \XX} \langle f, \psi_{\xi}\rangle \tilde\psi_{\xi}
\quad\hbox{in}\;\; \LL^p.
\end{equation}

$(b)$ {\rm Frame:}
The system $\{\tilde\psi_{\xi}\}$ as well as $\{\psi_{\xi}\}$
is a frame for $\LL^2$, namely, there exists a constant $c>0$ such that
\begin{equation}\label{frame-tpsi}
 c^{-1}\|f\|_2^2 \le
\sum_{\xi \in \XX} |\langle f, \tilde\psi_{\xi}\rangle|^2
\leq  c\| f\|_2^2,
\quad \forall f\in \LL^2.
\end{equation}

$(c)$ {\rm Space localization:} For any $0<\kaph<\kap/2$, $m \ge 0$, and any $\xi\in\XX_j$, $j\ge 0$,
\begin{equation}\label{prop-tpsi-1}
|L^m\tilde\psi_{\xi} (x)|
\le c_m b^{2jm}|B(\xi, b^{-j})|^{-1/2} \exp\big\{-\kaph(b^j\rho(x, \xi))^\bbb\big\},
\end{equation}
and if $\rho(x, y) \le b^{-j}$
\begin{equation}\label{prop-tpsi-Lip}
|\tilde\psi_{\xi} (x)- \tilde\psi_{\xi} (y)|
\le \ch|B(\xi, b^{-j})|^{-1/2}(b^j\rho(x, y))^\alpha \exp\big\{-\kaph(b^j\rho(x, \xi))^\bbb\big\}.
\end{equation}

$(d)$ {\rm Spectral localization:}
$\tilde\psi_{\xi}\in \Sigma_b^p$ if $\xi\in \XX_0$ and
$\tilde\psi_{\xi}\in \Sigma_{[b^{j-2}, b^{j+2}]}^p$ if $\xi\in \XX_j$, $j\ge 1$,
$0<p\le\infty$.

$(e)$ {\rm Norms:} If in addition the reverse doubling condition $(\ref{reverse-doubling})$ is valid, then
\begin{equation}\label{prop-tpsi-2}
\|\tilde\psi_{\xi}\|_p \sim |B(\xi, b^{-j})|^{\frac 1p-\frac 12}
\quad\hbox{for} \;\; 0< p \le \infty.
\end{equation}
\end{theorem}

\noindent
{\bf Proof.}
By the definition of $\tilde\psi_{\xi}$ in (\ref{rep-psi-tpsi}) and Lemma~\ref{lem:instrument}
we have
\begin{equation}\label{rep-tpsi}
\tilde\psi_{\xi}(x)
= c_\eps |A_\xi^j|^{1/2}\big[\LLam_{\lambda_j}(x, \xi)
+\SSS_{\lambda_j} [\LLam_{\lambda_j}(\cdot, \xi)](x)\big],
\quad \xi\in \cX_j.
\end{equation}
By the proof of Lemma~\ref{lem:instrument}
$\SSS_{\lambda} = \sum_{k\ge 1} R_\lambda = R_\lambda(\Id + \SSS_\lambda)$.
Hence, for $\xi\in \cX_j$
$$
L^m\tilde\psi_{\xi}(x)
= c_\eps |A_\xi^j|^{1/2}\Big(L^m\LLam_{\lambda_j}(x, \xi)
+ L^m R_{\lambda_j}[\LLam_{\lambda_j}(\cdot, \xi)](x)
+ L^m R_{\lambda_j}\SSS_{\lambda_j}[\LLam_{\lambda_j}(\cdot, \xi)](x)\Big).
$$
Clearly, the kernel $L^m R_{\lambda_j}(x, y)$ of the operator $L^m R_{\lambda_j}$ is given  by
$$
L^m R_{\lambda_j}(x, y) = L^m\LLam_{\lambda_j}^2(x, y)
+ \sum_{\xi\in\cX_j} w_\xi L^m\LLam_{\lambda_j}(x, \xi)\LLam_{\lambda_j}(\xi, y)
$$
and just as in the proof of Lemma~\ref{lem:instrument} we get
$
|L^mR_{\lambda_j}(x,y)| \le c b^{2jm}E_{\lambda_j^{-1}, \kappa}(x,y).
$
Now, this implies (\ref{prop-tpsi-1}) using (\ref{gen-local-1}), (\ref{gen-local-2}), (\ref{local-S}),
Theorem~\ref{thm:algebra}, and (\ref{D2}).
The~H\"{o}lder continuity estimate (\ref{prop-tpsi-Lip}) follows by (\ref{rep-tpsi}) using
(\ref{holder-f}) and (\ref{local-S}).
The other claims of the theorem are as in Theorem~5.3 in \cite{CKP}.
$\qed$

\subsection{Frames in the case when \boldmath $\{\Sigma_\lambda^2\}$  possess the polynomial property}
\label{sec:polyn-prop}

The~construction of frames with the desired excellent space and spectral localization is
simple and elegant in the case when the spectral spaces $\Sigma_\lambda^2$
have the polynomial property under multiplication:
{\em
Let $\{F_\lambda, \lambda\ge 0\}$ be the spectral resolution associated with the operator $\sqrt{L}$.
We say that the associated spectral spaces
$$
\Sigma_\lambda^2= \{ f \in \LL^2: F_\lambda f =f\}
$$
have the polynomial property if there exists a constant $\kk>1$ such that
\begin{equation}\label{polyp}
\Sigma_\lambda^2\cdot \Sigma_\lambda^2 \subset \Sigma^1_{\kk\lambda},
\quad \hbox{i.e.}\quad f, g\in \Sigma_\lambda^2 \Longrightarrow fg\in \Sigma_{\kk\lambda}^1.
\end{equation}
}

The construction begins with the introduction of a pair of cut-off functions
$\Psi_0, \Psi \in C^\infty(\bR_+)$ with the following properties:
\begin{align*}
&\supp \Psi_0\subset [0, b],
\quad \supp \Psi \subset [b^{-1}, b],
\quad
0\le \Psi_0, \Psi \le 1,\\
&\Psi_0(u) \ge c>0,\;\; u\in [0, b^{3/4}], \quad
\Psi(u) \ge c>0,\;\; u\in [b^{-3/4}, b^{3/4}],\\
& \Psi_0(u)=1 \quad u\in [0,1], %
\Psi_0^2(u)+\sum_{j\ge 1} \Psi^2(b^{-j}u)=1, \quad u\in \bR_+,
\end{align*}
and the kernels of the operators $\Psi_0(\delta\sqrt{L})$ and $\Psi(\delta\sqrt{L})$
have sub-exponential localization and H\"{o}lder continuity as in (\ref{sub-exp-Phi})-(\ref{m-sub-exp-Phi}).
Above $b>1$ is the constant from Theorem~\ref{thm:norms}.
The existence of functions like these follows by Theorem~\ref{thm:band-sub-exp}.

Set $\Psi_j(u):=\Psi(b^{-j}u)$.
Then
$\sum_{j\ge 0} \Psi_j^2(u)=1$, $u\in \bR_+$,
which leads to the following Calder\'{o}n type decomposition (see Proposition~\ref{prop:decomp-DD2} below)
\begin{equation}\label{calderon}
f=\sum_{j\ge 0} \Psi_j^2(\sqrt{L})f, \quad f\in \LL^p, \; 1\le p \le\infty,
\quad (L^\infty:=\UCB).
\end{equation}

The key observation is that the polynomial property (\ref{polyp}) of the spectral spaces allows to
discretize the above expansion and as a result to obtain the desired frame.
To be more specific, by construction
$\Psi_j(\sqrt{L})$ is a kernel operator whose kernels has
sub-exponential localization and
$\Psi_j(\sqrt{L})(x, \cdot)\in \Sigma_{b^{j+1}}$.
Now, choosing $\cX_j$ ($j\ge 0$) to be a maximal $\delta$-net on $M$
with $\delta:=\gamma \kk^{-1}b^{-j-1}\sim b^{-j}$
we get from (\ref{quadrature}) a~cubature formula of the form
$$
\int_M f(x) d\mu(x) = \sum_{\xi\in\cX_j} \ww_{j\xi} f(\xi)
\quad\hbox{for}\quad f\in \Sigma_{\kk b^{j+1}}^1,
$$
where $\frac 23|B(\xi, \delta/2)| \le \ww_{j\xi} \le 2|B(\xi, \delta)|$.
Since $\Psi_j(\sqrt{L})(x,\cdot)\Psi_j(\sqrt{L})(\cdot,y) \in \Sigma_{\kk b^{j+1}}^1$
by (\ref{polyp}), we can use the cubature formula from above to obtain
\begin{align}\label{Psi-Psi}
\Psi_j(\sqrt{L})\Psi_j(\sqrt{L})(x,y)
&= \int_M \Psi_j(\sqrt{L})(x,u)\Psi_j(\sqrt{L})(u,y) d\mu(u)\\
&=\sum_{\xi\in\cX_j} \ww_{j\xi} \Psi_j(\sqrt{L})(x,\xi)\Psi_j(\sqrt{L})(\xi,y).\notag
\end{align}
Now, the frame elements are defined by
\begin{equation}\label{def-psi-xi}
\psi_{\xi}(x):=\sqrt{\ww_{j\xi}}\Psi_j(\sqrt{L})(x,\xi),
\;\; \xi\in\cX_j, \; j\ge 0.
\end{equation}
As in \S\ref{natural-frame}, set $\cX:=\cup_{j\ge 0} \cX_j$.
It will be convenient to use $\cX$ as an index set and for this
equal points from different $\cX_j$'s will be regarded as distinct element of $\cX$.

Observe that
$\{\psi_{\xi}\}_{\xi\in\cX}$ is a tight frame for $\LL^2$.
More precisely, for any $f\in\LL^p$, $1\le p\le\infty$, $(\LL^\infty:=\UCB)$ we have
\begin{equation}\label{frame2-rep}
f = \sum_{\xi \in \cX} \langle f, \psi_{\xi}\rangle \psi_{\xi}
\quad\hbox{in}\;\LL^p\quad\hbox{and}\quad
\|f\|_2^2 = \sum_{\xi \in \cX}|\langle f, \psi_{\xi}\rangle|^2
\quad\hbox{for $f\in \LL^2$.}
\end{equation}
The convergence in (\ref{frame2-rep}) for test functions and distributions is given
in Proposition~\ref{prop:decomp-DD2} below.
Furthermore, the frame elements $\psi_{\xi}$ have all other properties of the elements constructed
in \S\ref{natural-frame} (see Proposition~\ref{prop:frame-prop}).

\section{Distributions}\label{distributions}\label{sec:distributions}
\setcounter{equation}{0}

The Besov and Triebel-Lizorkin spaces that will be developed are in general
spaces of distributions. There are some distinctions, however, between the tests functions
and distributions that we shall use depending on whether $\mu(M) <\infty$ or $\mu(M)=\infty$.
We shall clarify them in this section.

\subsection{Distributions in the case \boldmath $\mu(M) <\infty$}\label{sec:distr-compact}

To introduce distributions we shall use as test functions the class $\cD$ of all functions
$\phi\in  \cap_m D(L^m)$ with topology induced by
\begin{equation}\label{norm1_D}
\cP_m (\phi)
:= \|L^m\phi\|_2,
\quad m\ge 0,
\end{equation}
or equivalently by
\begin{equation}\label{norm11-D}
\cP_m^*(\phi) := \max_{0\le r \le m}\|L^r\phi\|_2,
\quad m\ge 0.
\end{equation}
The norms $\cP^*_m (\phi)$, $m=0, 1, \dots$, are usually more convenient since
they form a~directed family of norms.
Another alternative is to use the norms
\begin{equation}\label{norm2_D}
\cP^{**}_m (\phi) := \sup_{\lambda \ge 0}\, (1+\lambda)^m \|(\Id-E_\lambda) \phi\|_2,
\quad m=0, 1, \dots,
\end{equation}
where as before $E_\lambda$, $\lambda\ge 0$, is the spectral resolution associated with the operator~$L$.
The equivalence of the norms $\{\cP_m^* (\cdot)\}_{m\ge 0}$ and $\{\cP_m^{**}(\cdot)\}_{m\ge 0}$
follows by the identity
$$ %
\|L^m\phi\|_2^2
= \int_0^\infty \lambda^{2m} d(E_\lambda \phi, \phi)
= \int_0^\infty \lambda^{2m} d \|E_\lambda \phi\|_2^2.
$$ 
Indeed, clearly
$$
\lambda^{2m}\|( I_d - E_\lambda)\phi\|_2^2
=\lambda^{2m} \int_\lambda^\infty d\| E_t \phi \|_2^2
\le \int_0^\infty t^{2m} d\| E_t \phi \|_2^2 =\|L^m\phi\|_2^2,
$$
and hence
$\cP_m^{**}(\phi) \le c \cP_m^*(\phi)$.
On the other hand,
\begin{align*}
\|L^m\phi\|_2^2
&= \int_0^1 \lambda^{2m} d\|E_\lambda \phi \|_2^2
+ \sum_{j\ge 0}\int_{2^j}^{2^{j+1}} \lambda^{2m} d\|E_\lambda \phi \|_2^2\\
&\le \|\phi\|_2^2 + \sum_{j\ge 0}2^{(j+1)2m}\|(\Id-E_{2^j})\phi\|_2^2,
\end{align*}
implying
$\cP_m^{*}(\phi) \le c \cP_{m+1}^{**}(\phi)$.

In the next proposition we collect some simple facts about test functions.


\begin{proposition}\label{prop:char-D}
$(a)$
$\cD$ is a Fr\'{e}chet space.

$(b)$
$\Sigma_\lambda \subset \cD$, $\lambda \ge 0$, and
for every $\phi\in \cD$,
$
\phi=\lim_{\lambda\to\infty} E_\lambda \phi
$
in the topology of $\cD$.

$(c)$
If $\varphi$ is in the Schwartz class $\cS(\R)$ of $C^\infty$ rapidly decaying
$($with all their derivatives$)$ functions on $\bR$,
$\varphi$ is real-valued,
and $\varphi^{(2\nu+1)}(0)=0$ for $\nu=0, 1, \dots$, then the kernel
$\varphi(\sqrt L)(x, y)$ of the operator $\varphi(\sqrt L)$ belongs
to $\cD$ as a function of $x$ and as a function of $y$.
\end{proposition}

\noindent
{\bf Proof.}
Part (a) follows by the completeness of $\LL^2$ and the fact that $L$ being a self-adjoint operator is closed
(see the proof of Proposition~\ref{prop:char-D-infty} below).

Part (b) is also easy to prove, indeed, for $t>1$ we have for any $m\ge 1$
\begin{align*}
\cP^{**}_m (\phi- E_t\phi)
&= \sup_{u\ge 0} (1+u)^m \|(Id-E_u) (Id- E_t)\phi)\|_2\\
&= \sup_{u\ge t} (1+u)^m \|(Id- E_u)\phi)\|_2
\le \sup_{u\ge t} c_{m+1}u^{-1}
\le ct^{-1}
\end{align*}
and the claimed convergence follows.
Part (c) follows by Corollary~\ref{thm:LS-local-kernels}.
$\qed$

\smallskip


The space $\cD'$ of distributions on $M$ is defined
as the space of all continuous linear functionals on $\cD$.
The pairing of $f\in \cD'$ and $\phi\in\cD$ will be denoted by
$\langle f, \phi \rangle:= f(\overline{\phi})$; this will be consistent with the inner product
$
\langle f, g \rangle :=\int_M f\overline{g} d\mu
$
in $\LL^2$.

We shall be dealing with integral operators $\cH$ of the form
$\cH f(x):=\int_M \cH (x, \cdot)fd\mu$,
where $\cH(x, \cdot)\in \cD$ for all $x\in\cD$.
We set
\begin{equation}\label{def-Hf}
\cH f (x) := \langle f,  \overline{\cH(x,\cdot)} \rangle
\quad\mbox{for}\quad f\in\cD',
\end{equation}
where on the right $f$ acts on $\cH(x, y)$ as a function of $y$.

As is shown in \cite{CKP}, \S 3.7, in the case  $\mu(M) <\infty$ the spectrum
of $L$ is discrete and hence the spectrum
of the operator $\sqrt L$ is discrete as well.
Furthermore, the spectrum of $\sqrt L$ is of the form
$\spec \sqrt L=\{\lambda_1, \lambda_2, \dots\}$, where $0\le \lambda_1 < \lambda_2<\dots$ and
$\lambda_n \to\infty$.
Also, the eigenspace $\cE_\lambda$ associated with each $\lambda\in\spec \sqrt L$
is of finite dimension, say, $N_\lambda$.
Let $\{e_{\lambda m}: m=1, 2, \dots, N_\lambda\}$ be an orthonormal basis for~$\cE_\lambda$.
Then
$
E_t(x, y)=\sum_{0\le \lambda\le t}\sum_{m=1}^{N_\lambda}
e_{\lambda m}(x) \overline{e_{\lambda m}(y)}
$
is the kernel of the projector~$E_t$.
Therefore, for any distribution $f\in\cD'$
\begin{equation}\label{Et-f}
E_tf=\langle f, \overline{E_t(x, \cdot)}\rangle
=\sum_{0\le \lambda\le t}\sum_{m=1}^{N_\lambda}
\langle f, e_{\lambda m}\rangle e_{\lambda m}(x).
\end{equation}
Consequently, for any $f\in\cD'$ we have $E_tf\in \Sigma_t=\bigoplus_{\lambda\le t}\cE_\lambda$.

We collect this and some other simple fact about the distributions we introduced above
in the following


\begin{proposition}\label{prop:dec-D1}
$(a)$
A linear functional $f$ belongs to $\cD'$
if and only if there exist $m\ge 0$ and $c_m>0$ such that
\begin{equation}\label{Dis1}
|\langle f, \phi\rangle|\le c_m\cP_m^*(\phi)
\quad \mbox{for all } \; \phi \in \cD.
\end{equation}
Hence, for any $f \in \cD'$ there exist $m\ge 0$ and $c_m>0$ such that
\begin{equation}\label{Dis2}
\|(\Id- E_\lambda) f\|_2
\le c_m(1+\lambda)^{-m}, \quad \forall \lambda\ge 1.
\end{equation}

$(b)$
For any $f\in\cD'$ we have $E_\lambda f\in \Sigma_\lambda$ and also
$
 \langle E_\lambda f, \phi \rangle= \langle f, E_\lambda \phi \rangle
$
for all $\phi\in\cD$.

$(c)$
For any $f\in\cD'$ we have
$
f=\lim_{\lambda\to\infty} E_\lambda f
$
in distributional sense, i.e.
\begin{equation}\label{Dis3}
\langle f, \phi\rangle
= \lim_{\lambda\to\infty}\langle E_\lambda f, \phi  \rangle
= \lim_{t\to\infty}\langle E_\lambda f, E_\lambda \phi  \rangle
\quad\mbox{for all}
\quad \phi\in\cD.
\end{equation}

$(d)$
If $\varphi\in C^\infty_0(\bR_+)$, $\varphi^{(2\nu+1)}(0)=0$ for $\nu=0, 1, \dots$,
and $\supp \varphi \subset [0, R]$, then for any $\delta >0$ and $f\in \cD'$
we have $\varphi (\delta\sqrt L) f\in \Sigma_{R/\delta}$.
\end{proposition}

\noindent
{\bf Proof.}
Part (a) follows at once by the fact that the
topology in $\cD$ can be defined by the norms $\cP_m^*(\cdot)$ from (\ref{norm2_D}).
Part (b) follows from (\ref{Et-f}).
For Part (c) we use Proposition~\ref{prop:char-D} (b) and (b) from above to obtain
$$
\langle f, \phi \rangle = \lim_{t\to\infty}\langle f, E_t \phi \rangle
= \lim_{t\to\infty}\langle E_tf, \phi \rangle
= \lim_{t\to\infty}\langle E_tf, E_t\phi \rangle,
$$
which completes the proof.
The proof of Part (d) is similar taking into account that the integral is actually
a discrete sum.
$\qed$

Basic convergence results for distributions will be given in the next subsection.

\subsection{Distributions in the case \boldmath $\mu(M)=\infty$}\label{sec:distr-noncompact}

In this case the class of test functions $\cD$ is defined as the set of all functions
$\phi\in  \cap_m D(L^m)$ such that
\begin{equation}\label{norm-D-infty}
\cP_{m, \ell}(\phi)
:= \sup_{x\in M}
(1+\rho(x, x_0))^\ell |L^m\phi(x)| < \infty
\quad \forall m, \ell\ge 0.
\end{equation}
Here $x_0\in M$ is selected arbitrarily and fixed once and for all.
Clearly, the~particular selection of $x_0$ in the above definition is not important,
since if  $\cP_{m, \ell} (\phi) <\infty$ for one $x_0\in M$, then
$\cP_{m, \ell} (\phi) <\infty$ for any other selection of $x_0\in M$.

It is often more convenient to have a directed family of norms.
For this reason we introduce the following norms on $\cD$:
\begin{equation}\label{norm-D-infty2}
\cP_{m, \ell}^*(\phi)
:= \max_{0\le r \le m, 0\le l\le \ell} \cP_{r, l}(\phi).
\end{equation}

Note that unlike in the case $\mu(M) <\infty$, in general,
$\Sigma_t^p \not\subset \cD$.
However, there are still sufficiently many test functions.
This becomes clear from the following


\begin{proposition}\label{prop:char-D-infty}
$(a)$
$\cD$ is a Fr\'{e}chet space
and $\cD\subset {\rm UCB}$.

$(b)$
If $\varphi$ is in the Schwartz class $\cS(\R)$, $\varphi$ is real-valued 
and $\varphi^{(2\nu+1)}(0)=0$ for $\nu=0, 1, \dots$, then the kernel
$\varphi(\sqrt L)(x, y)$ of the operator $\varphi(\sqrt L)$ belongs
to $\cD$ as a function of $x$ $($and as a function of $y$$)$.
Moreover,
$\varphi(\sqrt L) \phi \in \cD$ for any $\phi\in\cD$.
Also, $e^{-tL}(x,\cdot) \in \cD$ and $e^{-tL}(\cdot, y) \in \cD$, $t>0$.
\end{proposition}

\noindent
{\bf Proof.}
To prove that $\cD$ is a Fr\'{e}chet space we only have to establish the completeness of $\cD$.
Let $\{\phi_j\}_{j\ge 1}$ be a Cauchy sequence in $\cD$, i.e.
$\cP_{m, \ell}(\phi_j -\phi_n) \to 0$ as $j,n \to \infty$ for all $m, \ell\ge 0$.
Choose $\ell \in \bN$ so that $\ell \ge (d+1)/2$. Then clearly for any $m\ge 0$
\begin{align*}
\|L^m\phi_j - L^m\phi_n\|_2
&\le \cP_{m, \ell}(\phi_j -\phi_n)\int_M(1+\rho(x, x_0)^{-d-1}d\mu(x)\\
&\le c|B(x_0, 1)|\cP_{m, \ell}(\phi_j -\phi_n),
\end{align*}
where we used (\ref{tech1}).
Therefore,
$\|L^m\phi_j - L^m\phi_n\|_2 \to 0$ as $j,n \to \infty$
and by the completeness of $L^2$ there exists $\Psi_m\in L^2$ such that
$\|L^m\phi_j - \Psi_m\|_2 \to 0$ as $j \to \infty$.
Write $\phi:=\Psi_0$.
From $\|\phi_j - \phi\|_2 \to 0$, $\|L\phi_j - \Psi_1\|_2 \to 0$,
and the fact that $L$ being a self-adjoint operator is closed \cite{RS}
it follows that $\phi\in D(L)$ and
$\|L\phi_j - L\phi\|_2 \to 0$.
Using the same argument inductively we conclude that $\phi\in \cap_m D(L^m)$ and
\begin{equation}\label{L2-conv}
\|L^m\phi_j - L^m\phi\|_2 \to 0
\quad\hbox{as $j\to\infty$\; for all\; $m\ge 0$.}
\end{equation}
On the other hand,
$\|L^m\phi_j-L^m\phi_n\|_\infty = \cP_{m, 0}(\phi_j-\phi_n) \to 0$ as $j, n\to\infty$
and from the completeness of $L^\infty$ the sequence
$\{L^m\phi_j\}_{j\ge 0}$ converges in $L^\infty$.
This and (\ref{L2-conv}) yield
\begin{equation}\label{Linfty-conv}
\|L^m\phi_j - L^m\phi\|_\infty \to 0
\quad\hbox{as $j\to\infty$ \; for all \; $m\ge 0$.}
\end{equation}
In turn, this along with $\cP_{m, \ell}(\phi_j -\phi_n) \to 0$ as $j,n \to \infty$ implies
$\cP_{m, \ell}(\phi_j -\phi) \to 0$ as $j\to \infty$ for all $m,\ell\ge 0$,
which confirms the completeness of $\cD$.


In Proposition~\ref{prop:decomp-DD2} (a) below it will be shown that any $\phi\in\cD$
can be approximated in $L^\infty$ by H\"{o}lder continuous functions,
which implies that $\phi$ is uniformly continuous and hance $\cD\subset {\rm UCB}$.


For the proof of Part (b), we note that
if $\varphi \in \cS(\R)$ 
and $\varphi^{(2\nu+1)}(0)=0$ for $\nu=0, 1, \dots$, then by Theorem~\ref{thm:LS-local-kernels}
$L^m\varphi(\sqrt L)$ is an integral operator whose kernel obeys
\begin{equation}\label{test-func-D-infty}
|L^m\varphi(\sqrt L)(x, y)|
\le c_{\sigma, m} |B(x, 1)|^{-1}\big(1+\rho(x, y)\big)^{-\sigma}
\quad \hbox{for all}\quad \sigma>0,\; m\ge 0.
\end{equation}
Therefore, $\varphi(\sqrt L)(x, \cdot) \in \cD$ with $x$ fixed,
and $\varphi(\sqrt L)(\cdot, y) \in \cD$ with $y$ fixed.
These follow by (\ref{test-func-D-infty}) and the identity
\begin{equation}\label{Lm-phi}
L^m\big[\varphi(\sqrt L)(x, \cdot)\big] = L^m\varphi(\sqrt L)(x, \cdot)
\quad\hbox{for any fixed $x\in M$.}
\end{equation}
To prove this, suppose first that $\varphi \in C^\infty_0(\R) $.
Then $h := \varphi(\sqrt L)(x, \cdot) \in \cup_\lambda \Sigma_\lambda$ and hence
$L^m h \in  \cup_\lambda \Sigma_\lambda$, which implies $h, L^mh  \in \cap_k D(L^k)$.
For $ \theta \in \cup_\lambda \Sigma_\lambda$
$$
\int_M L^m h(u) \theta(u) d\mu(u) = \int_M  \varphi(\sqrt L)(x, u)  L^m\theta(u) d\mu(u)
= \varphi(\sqrt L)(L^m\theta)(x)
$$
$$= [ \varphi(\sqrt L)L^m ]\theta(x) = [L^m\varphi(\sqrt L)]\theta(x)
= \int_M [L^m \varphi(\sqrt L)] (x,u) \theta(u) d\mu(u).
$$
Here we used that $\overline{L^m\theta} = L^m\overline{\theta}$.
Now, we derive (\ref{Lm-phi}) for $\varphi \in \cS(\R)$ by a limiting argument.

In going further, from above it readily follows that $\varphi(\sqrt L) g \in \cD$ for any $g\in\cD$.
Also, Corollary~\ref{thm:LS-local-kernels} yields
$e^{-tL}(x,\cdot) \in \cD$ and $e^{-tL}(\cdot, y) \in \cD$, $t>0$.
$\qed$

\smallskip

As usual the space $\cD'$ of distributions on $M$ is defined
as the set of all continuous linear functionals on $\cD$ and
the pairing of $f\in \cD'$ and $\phi\in\cD$ will be denoted by
$\langle f, \phi \rangle:= f(\overline{\phi})$.

We next record some basic properties of distributions in the case $\mu(M)=\infty$.


\begin{proposition}\label{prop:dec-D1-inft}
$(a)$
A linear functional $f$ belongs to $\cD'$
if and only if there exist $m, \ell \ge 1$ and a constant $c>0$ such that
\begin{equation}\label{Dis-221}
|\langle f, \phi\rangle|\le c\cP_{m, \ell}^* (\phi)
\quad \mbox{for all } \quad \phi \in \cD.
\end{equation}

$(b)$ If $\varphi\in C^\infty_0(\bR_+)$, $\varphi^{(2\nu+1)}(0)=0$ for $\nu=0, 1, \dots$,
and $\supp \varphi \subset [0, R]$, then for any $\delta >0$ and $f\in \cD'$
we have $\varphi (\delta\sqrt L) f\in \Sigma_{R/\delta}^p$ for $0< p \le \infty$.
\end{proposition}

\noindent
{\bf Proof.}
Part (a) is immediate from the definition of distributions
and (b) follows by the fact that the kernel $\varphi (\delta\sqrt L)(x, y)$ of the operator
$\varphi (\delta\sqrt L)$ belongs to $\cD\cap \Sigma_{R/\delta}$
as a function of $x$ and as a function of $y$.
$\qed$

\smallskip

We now give our main convergence result for distributions and in $\LL^p$.


\begin{proposition}\label{prop:decomp-DD2}
$(a)$
Let $\varphi \in C^\infty(\bR_+)$, $\varphi$ be real-valued, $\supp \varphi \subset [0, R]$, $R>0$,
$\varphi(0)=1$, and $\varphi^{(\nu)}(0)=0$ for $\nu\ge 1$.
Then for any $\phi\in \cD$
\begin{equation}\label{decomp-dist-1}
\phi=\lim_{\delta\to 0} \varphi (\delta\sqrt L) \phi
\quad \mbox{ in }\; \cD,
\end{equation}
and for any $f\in \cD'$
\begin{equation}\label{decomp-dist-2}
f=\lim_{\delta\to 0} \varphi (\delta\sqrt L) f
\quad \mbox{ in }\; \cD'.
\end{equation}

$(b)$
Let $\varphi_0, \varphi \in C^\infty(\R_+)$, $\varphi_0, \varphi$ be real-valued,
$\supp \varphi_0 \subset [0, b]$ and $\supp \varphi \subset [b^{-1}, b]$ for some $b>1$,
and
$
\varphi_0(\lambda)+\sum_{j\ge 1}\varphi(b^{-j}\lambda) =1
$
for $\lambda\in \R_+$.
Set $\varphi_j(\lambda):= \varphi(b^{-j}\lambda)$, $j\ge 1$;
hence
$\sum_{j\ge 0} \varphi_j(\lambda)=1$ on $\bR_+$.
Then for any $f\in \cD'$
\begin{equation}\label{decomp-D12}
f=\sum_{j\ge 0}\varphi_j(\sqrt L)f \;\;\mbox{ in }\; \cD'.
\end{equation}

$(c)$
Let $\{\psi_\xi\}_{\xi\in\cX}$, $\{\tilde\psi_\xi\}_{\xi\in\cX}$ be the pair of frames
from $\S\S\ref{natural-frame}-\ref{dual-frame}$.
Then for any $f\in \cD'$
\begin{equation}\label{decomp-D3}
f=\sum_{\xi\in\cX}\langle f, \tilde\psi_\xi \rangle \psi_\xi
=\sum_{\xi\in\cX}\langle f, \psi_\xi \rangle \tilde\psi_\xi
\;\;\mbox{ in }\; \cD'.
\end{equation}
Furthermore, $(\ref{decomp-dist-2})-(\ref{decomp-D3})$ hold in $\LL^p$ for any $f\in\LL^p$,
$1\le p\le \infty$ $(\LL^\infty:=\UCB)$.
\end{proposition}

\noindent
{\bf Proof.}
We shall only consider the case $\mu(M)=\infty$.
The case $\mu(M)<\infty$ is easier.
For the proof of Part (a) it suffices to prove only (\ref{decomp-dist-1}),
since then (\ref{decomp-dist-2}) follows by duality.
To prove (\ref{decomp-dist-1}) we have to show that for any $m, \ell\ge 0$
$$
\lim_{\delta\to 0} \cP_{m, \ell}\big(\phi-\varphi(\delta\sqrt{L})\phi\big)
=\lim_{\delta\to 0} \sup_{x\in M}(1+\rho(x, x_0))^\ell|L^m[\phi-\varphi(\delta\sqrt{L})\phi](x)| =0.
$$
Given $m, \ell\ge 0$,
pick the smallest $k, r\in \bN$ so that
$k\ge \ell+ 5d/2$
and
$2r\ge k+d+1$.
Set $\omega(\lambda):=\lambda^{-2r}(1-\varphi(\lambda))$.
Then $1-\varphi(\delta\sqrt{\lambda}) = \delta^{2r}\omega(\delta\sqrt{\lambda})\lambda^{r}$
and hence
\begin{align*}
L^m[\phi-\varphi(\delta\sqrt{L})\phi](x)
&= \delta^{2r}\omega(\delta\sqrt{L})L^{m+r}\phi(x)\\
&= \delta^{2r}\int_M \omega(\delta\sqrt{L})(x, y)L^{m+r}\phi(y) d\mu(y).
\end{align*}
From the definition of $\omega$ we have $\omega\in C^\infty(\bR_+)$,
$\omega^{(\nu)}(0)=0$ for $\nu\ge 0$, and
$$
|\omega^{(\nu)}(\lambda)| \le c_\nu (1+\lambda)^{-2r},
\quad \lambda\in \bR_+, \;\nu\ge 0.
$$
Now, we apply Theorem~\ref{thm:S-local-kernels}, taking into account that $k\ge d+1$ and $2r\ge k+d+1$,
to conclude that the kernel $\omega(\delta\sqrt{L})(x, y)$ of the operator
$\omega(\delta\sqrt{L})$ obeys
$$ 
|\omega(\delta\sqrt{L})(x, y)| \le c_k D_{\delta, k}(x, y)
\le \frac{c}{|B(x, \delta)|(1+\delta^{-1}\rho(x, y))^{k-d/2}}.
$$ 
By (\ref{doubling}), (\ref{D2}) it readily follows that for $0<\delta <1$
$$
|B(x_0, 1)| \le c_0(1+\rho(x, x_0))^d|B(x, 1)| \le c_0^2\delta^{-d}(1+\rho(x, x_0))^d|B(x, \delta)|.
$$
Also, since $\phi\in \cD$ we have
$
|L^{m+r}\phi(x)| \le c(1+\rho(x, x_0))^{-k}.
$
Putting all of the above together we get
\begin{align*}
& (1+\rho(x, x_0))^\ell \big|L^m[\phi-\varphi(\delta\sqrt{L})\phi](x)\big|\\
& \qquad\qquad\qquad
\le c\frac{\delta^{2r-d}}{|B(x_0, 1)|}\int_M
\frac{(1+\rho(x, x_0))^{\ell+d}}{(1+\rho(x, y))^{k-d/2}(1+\rho(y, x_0))^{k-d/2}}
d\mu(y)\\
& \qquad\qquad\qquad
\le c\frac{\delta^{2r-d}(1+\rho(x, x_0))^{\ell+d}}{(1+\rho(x, x_0))^{k-3d/2}}
\le c \delta \to 0
\quad \hbox{as} \quad \delta\to 0.
\end{align*}
Here for the latter estimate we used that $k\ge \ell+5d/2$ and $2r > d+1$,
and for the former we used (\ref{tech2}).
This completes the proof of Part (a).

To show Part (b), set $\theta(\lambda):= \varphi_0(\lambda) + \varphi(b^{-1}\lambda)$
and note that
$\sum_{k=0}^j \varphi_k(\lambda)= \theta(b^{-j}\lambda)$ for $j\ge 1$.
Then the result follows readily by Part (a).


For the proof of Part (c) it suffices to show that
\begin{equation}\label{decomp-Dk-left}
\phi =\sum_{\xi\in\cX}\langle \phi, \psi_\xi \rangle \tilde\psi_\xi
\quad\hbox{and}\quad
\phi=\sum_{\xi\in\cX}\langle \phi, \tilde\psi_\xi \rangle \psi_\xi
\;\;\mbox{ in }\; \cD
\;\;\mbox{ for all }\; \phi\in \cD.
\end{equation}
We shall only prove the left-hand side identity in (\ref{decomp-Dk-left});
the proof of the right-hand side identity is similar.
Let $\{\Psi_j\}_{j\ge 0}$ be from the definition of $\{\psi_\xi\}$ in \S\ref{natural-frame}.
Then
$\sum_{j\ge 0}\Psi_j(u) = 1$ for $u\in\bR_+$
and by Part (b)
$\phi=\sum_{j\ge 0}\Psi_j(\sqrt{L})\phi$ in $\cD$ for all $\phi\in\cD$.
Therefore, to prove the left-hand side identity in (\ref{decomp-Dk-left}) it suffices to show that for each $j\ge 0$
\begin{equation}\label{decomp-Psij-left}
\Psi_j(\sqrt{L})\phi =\sum_{\xi\in\cX_j}\langle \phi, \psi_\xi \rangle \tilde\psi_\xi
\;\;\mbox{ in }\; \cD \;\;\forall \phi\in\cD.
\end{equation}
By (\ref{rep2-Psi-j})
$$
\Psi_j(\sqrt L)(x, y)=
\sum_{\xi \in \cX_j} \psi_{\xi}(y)\tilde\psi_{\xi}(x),\quad x, y\in M.
$$
From this and the sub-exponential space localization of
$\psi_{\xi}(x)$ and $L^m\tilde\psi_{\xi}(x)$, $m\ge 0$, given in (\ref{prop-psi-1}) and (\ref{prop-tpsi-1})
(see also (\ref{est-psi-51})-(\ref{est-psi-52}))
it readily follows that
$$
L^m\Psi_j(\sqrt L)(x, y)
=\sum_{\xi \in \cX_j} \psi_{\xi}(y)L^m\tilde\psi_{\xi}(x),
\quad x, y\in M,
$$
and hence
$$
L^m\Psi_j(\sqrt L)\phi
=\sum_{\xi \in \cX_j} \langle \psi_{\xi}, \overline\phi \rangle L^m\tilde\psi_{\xi},
\quad \forall \phi\in\cD.
$$
Clearly, to prove (\ref{decomp-Psij-left}) it suffices to
show that for any $\ell, m\ge 0$ and $\phi\in\cD$
\begin{equation}\label{conv-1}
\lim_{K\to\infty}\sup_{x\in M} (1+\rho(x, x_0))^\ell
\sum_{\xi \in \cX_j: \, \rho(\xi, x_0)\ge K} \int_M |\psi_{\xi}(y)\phi(y)|d\mu(y)
|L^m\tilde\psi_{\xi}(x)| =0.
\end{equation}
Given $\ell, m\ge 0$, choose $\sigma \ge \ell +3d+1$.
From (\ref{prop-psi-1}) and (\ref{prop-tpsi-1}) it follows that
\begin{align}
|\psi_{\xi}(x)|
&\le c_\sigma |B(\xi, b^{-j})|^{-1/2}(1+b^j\rho(x, \xi))^{-\sigma},\label{est-psi-51}\\
|L^m\tilde\psi_{\xi}(x)|
&\le c_{\sigma, m} b^{2jm}|B(\xi, b^{-j})|^{-1/2}(1+b^j\rho(x, \xi))^{-\sigma}, \xi\in\cX_j.\label{est-psi-52}
\end{align}
On the other hand, since
$\phi\in\cD$ we have $|\phi(x)| \le c(1+\rho(x, x_0))^{-\sigma}$.
Therefore,
\begin{align*}
&\int_M |\psi_{\xi}(y)\phi(y)|d\mu(y)
\le c\int_M\frac{d\mu(y)}
{|B(\xi, b^{-j})|^{1/2}(1+b^j\rho(y, \xi))^{\sigma}(1+\rho(y, x_0))^{\sigma}}\\
&\qquad\qquad
\le c|B(\xi, b^{-j})|^{1/2}\int_M\frac{d\mu(y)}
{|B(y, b^{-j})|(1+b^j\rho(y, \xi))^{\sigma-d}(1+\rho(y, x_0))^{\sigma}}\\
&\qquad\qquad
\le \frac{c|B(\xi, b^{-j})|^{1/2}}{(1+\rho(\xi, x_0))^{\sigma-d}},
\end{align*}
where for the second inequality we used (\ref{D2}) and for the last inequality
(\ref{tech5}).
From above and (\ref{est-psi-52})
\begin{align*}
&(1+\rho(x, x_0))^\ell
\sum_{\xi \in \cX_j: \, \rho(\xi, x_0)\ge K} \int_M |\psi_{\xi}(y)||\phi(y)|d\mu(y)
|L^m\tilde\psi_{\xi}(x)|\\
&\qquad
\le \sum_{\xi \in \cX_j: \, \rho(\xi, x_0)\ge K}
\frac{cb^{2jm}(1+\rho(x, x_0))^\ell}{(1+\rho(\xi, x_0))^{\sigma-d}(1+b^j\rho(\xi, x))^\sigma}\\
&\qquad
\le \sum_{\xi \in \cX_j: \, \rho(\xi, x_0)\ge K}
\frac{cb^{2jm}}{(1+\rho(\xi, x_0))^{\sigma-\ell-d}(1+b^j\rho(\xi, x))^{\sigma-\ell}}\\
&\qquad
\le \frac{cb^{2jm}}{(1+K)^{\sigma-\ell-d}}\sum_{\xi \in \cX_j}
\frac{1}{(1+b^j\rho(\xi, x))^{\sigma-\ell}}
\le \frac{cb^{2jm}}{1+K} \to 0
\quad\mbox{as}\quad K\to\infty.
\end{align*}
Here for the second inequality we used that
$1+\rho(x, x_0) \le (1+\rho(\xi, x_0))(1+\rho(\xi, x))$
and for the last inequality we used (\ref{discr-tech11}).
The above implies (\ref{conv-1}) and the proof of (\ref{decomp-D3}) is complete.


The convergence of $(\ref{decomp-dist-2})-(\ref{decomp-D3})$ in $\LL^p$ for $f\in \LL^p$
follows by a standard argument,
see also Theorem~5.3 in~\cite{CKP}.
$\qed$

\subsection{Distributions on \boldmath $\bR^d$ and $\bT^d$ induced by $L= -\Delta$}

The purpose of this subsection is to show that in the cases of $M=\bT^d$ and $M=\bR^d$
with $L= -\Delta$ ($\Delta$ being the Laplace operator)
the distributions defined as in \S\S\ref{sec:distr-compact}-\ref{sec:distr-noncompact}
are just the classical distributions on the torus $\bT^d$ and the tempered distributions on $\bR^d$.

The case of $M=\bT^d=\bR^d/\bZ^d$ and $L= -\Delta$ is quite obvious. The eigenfunctions of $-\Delta$
are $e^{2\pi ik \cdot x}$, $k\in \bZ^d$.
Clearly, in this case the class of test functions $\cD$ defined in \S\ref{sec:distr-compact}
consists of all functions $\phi\in L^2(\bT^d)$ whose Fourier coefficients $\hat\phi(k)$
obey $|\hat\phi(k)| \le c_N(1+|k|)^{-N}$ for each $N>0$.
It is easy to see that this is necessary and sufficient for $\phi\in C^\infty(\bT^d)$.
Therefore, $\cD=C^\infty(\bT^d)$ as in the classical case.
For more details, see e.g. \cite{Edwards}.
Observe that the situation with distributions on the unit sphere $\bS^{d-1}$ in $\bR^d$ is quite similar,
see e.g. \cite{NPW2}.

The case of $M=\bR^d$ and $L= -\Delta$ is not so obvious and since we do not find the argument in the literature
we shall consider it in more detail. Note first that in this case the class of test functions $\cD$ defined
in \S\ref{sec:distr-noncompact} consists of all functions
\begin{equation}\label{D-on-Rd}
\phi\in C^\infty(\bR^d)\quad\hbox{s.t.}\quad
\cP_{m,\ell}(\phi):= \sup_{x\in \bR^d}(1+|x|)^\ell|\Delta^m \phi(x)| <\infty, \quad \forall m, \ell \ge 0.
\end{equation}
Recall that the Schwartz class $\cS$ on $\bR^d$ consists of all functions $\phi\in C^\infty(\bR^d)$ such that
$\|\phi\|_{\alpha, \beta, \infty}:=\sup_x|x^\alpha \partial^\beta\phi(x)| <\infty$ for all multi-indices $\alpha, \beta$.
We shall also need the semi-norms $\|\phi\|_{\alpha, \beta, 2}:=\|x^\alpha \partial^\beta\phi\|_{L^2}$.
It is well known %
that on $\cS$ the semi-norms $\{\|\phi\|_{\alpha, \beta, \infty}\}$
are equivalent to the semi-norms $\{\|\phi\|_{\alpha, \beta, 2}\}$,
see e.g. \cite{RS}, Lemma 1, p.~141.


\begin{proposition}\label{prop:distr-Td-Rd}
The classes $\cD$ $($defined in \S\ref{sec:distr-noncompact}$)$
and $\cS$ on $\bR^d$ are the same with the same topology.
\end{proposition}

\noindent
{\bf Proof.}
We only have to prove that $\cD \subset \cS$, since obviously $\cS \subset \cD$.

Assume $\phi\in \cD$, i.e. $\phi\in C^\infty(\bR^d)$ and
$\cP_{m,\ell}(\phi)<\infty$, $\forall m, \ell \ge 0$.
This readily implies
$\|x^\alpha\Delta^m \phi\|_2 <\infty$ for all multi-indices $\alpha$ and $m\ge 0$.
Denoting by $\hat\phi(\xi)$ the Fourier transform of $\phi$ we infer using Plancherel's identity
\begin{equation}\label{alpha-beta-1}
\|\partial^\alpha (|\xi|^{2m}\hat\phi)\|_2 <\infty, \quad \forall\alpha \hbox{ and } m\ge 0.
\end{equation}
We claim that this yields
\begin{equation}\label{alpha-beta-2}
\|\xi^\alpha\partial^\beta\hat\phi\|_2 <\infty, \quad \forall\alpha, \beta.
\end{equation}
We shall carry out the proof by induction in $|\beta|$.
Indeed, (\ref{alpha-beta-2}) when $|\beta|=0$ is immediate from (\ref{alpha-beta-1}) with $|\alpha| =0$.
Clearly,
\begin{equation}\label{partial-j}
\partial_j (|\xi|^{2m}\hat\phi(\xi))= 2m\xi_j|\xi|^{2m-2}\hat\phi(\xi) + |\xi|^{2m}\partial_j\hat\phi(\xi)
\end{equation}
and hence
$$
\||\xi|^{2m}\partial_j\hat\phi\|_2
\le \|\partial_j (|\xi|^{2m}\hat\phi)\|_2 + 2m\||\xi|^{2m-1}\hat\phi\|_2 <\infty,
$$
where we used (\ref{alpha-beta-1}) and the already established (\ref{alpha-beta-2}) when $|\beta|=0$.
The above yields (\ref{alpha-beta-2}) for $|\beta|=1$ and all multi-indices $\alpha$.

We differentiate (\ref{partial-j}) and use (\ref{alpha-beta-1}) and that (\ref{alpha-beta-2}) holds for $|\beta|=0, 1$
and all $\alpha$'s just as above to show that (\ref{alpha-beta-2}) holds when $|\beta|=2$
and for all multi-indices~$\alpha$.
We complete the proof of (\ref{alpha-beta-2}) by induction.

Applying the inverse Fourier transform we obtain from (\ref{alpha-beta-2})
$$
\|\partial^\alpha(x^\beta\phi)\|_2 <\infty, \quad \forall\alpha, \beta.
$$
In turn, just as above this leads to
$\|\phi\|_{\alpha, \beta, 2}=\|x^\alpha \partial^\beta\phi\|_2 <\infty$, $\forall\alpha, \beta$.
As was mentioned, the semi-norms $\{\|\phi\|_{\alpha, \beta, 2}\}$
are equivalent to the semi-norms $\{\|\phi\|_{\alpha, \beta, \infty}\}$.
Therefore,
$\|\phi\|_{\alpha, \beta, \infty}=\|x^\alpha \partial^\beta\phi\|_\infty <\infty$,
$\forall\alpha, \beta$,
and hence $\cD \subset \cS$.
Clearly, the equivalence of the semi-norms  $\{\cP_{m,\ell}(\phi)\}$ and
$\{\|\phi\|_{\alpha, \beta, \infty}\}$ follows from the above considerations.
$\qed$

\section{Besov spaces}\label{besov-spaces}
\setcounter{equation}{0}

We shall use the well known general idea %
\cite{Peetre, Triebel1, Triebel2}
of employing spectral decompositions induced by a self-adjoint positive operator
to introduce (inhomogeneous) Besov spaces in the general set-up of this paper.
A~new point in our development is that we allow the smoothmess to be negative
and $p<1$.
To better deal with possible anisotropic geometries we introduce two types of Besov spaces:
(i)  classical Besov spaces  $B_{pq}^{s}=B_{pq}^{s}(L)$,
which for $s>0$ and $p\ge 1$ can be identified as approximation spaces of
linear approximation from $\Sigma_t^p$ in $\LL^p$,
and
(ii)  nonclassical Besov spaces $\tB_{pq}^{s}=\tB_{pq}^{s}(L)$, which for certain indices appear in
nonlinear approximation.
We shall utilize real-valued functions
$\varphi_0, \varphi\in C^\infty(\bR_+)$ such that
\begin{align}
&\supp \varphi_0 \subset   [0, 2] ,\;
\varphi_0^{(2\nu+1)}(0) = 0 \hbox{ for } \nu\ge 0,\;
|\varphi_0(\lambda)| \ge c>0 \;\hbox{ for } \lambda\in [0, 2^{3/4}], \label{cond_phi}\\
&\supp \varphi \subset   [1/2, 2], \;\;
|\varphi(\lambda)| \ge c>0 \;\hbox{ for } \lambda\in [2^{-3/4}, 2^{3/4}]. \label{cond_psi}
\end{align}
Then
$|\varphi_0(\lambda)| +\sum_{j\ge 1} |\varphi(2^{-j}\lambda)| \ge c >0$,
$\lambda \in  \R_+$.
Set $\varphi_j(\lambda):= \varphi(2^{-j}\lambda)$ for $j\ge 1$.


\begin{definition}\label{def-B-spaces}
Let $s \in \R$ and $0<p,q \le \infty$.

$(i)$ The Besov space  $B_{pq}^{s}=B_{pq}^{s}(L)$
is defined as the set of all $f \in \cD'$ such that
\begin{equation}\label{def-Besov-space1}
\|f\|_{B_{pq}^{s}} :=
\Big(\sum_{j\ge 0} \Big(2^{s j}
\|\varphi_j(\sqrt{L}) f(\cdot)\|_{\LL^p}
\Big)^q\Big)^{1/q} <\infty.
\end{equation}

$(ii)$ The Besov space  $\tB_{pq}^{s}= \tB_{pq}^{s}(L)$ is defined as the set
of all $f \in \cD'$ such that
\begin{equation}\label{def-Besov-space2}
\|f\|_{\tB_{pq}^{s}} :=
\Big(\sum_{j\ge 0} \Big(
\| |B(\cdot, 2^{-j})|^{-s/d}
\varphi_j(\sqrt{L}) f(\cdot)\|_{\LL^p}
\Big)^q\Big)^{1/q} <\infty.
\end{equation}
Above the $\ell^q$-norm is replaced by the sup-norm if $q=\infty$.
\end{definition}


\begin{remark}\label{rem:dimension}
{\rm
A word of caution concerning the smoothness parameter $s$ is in order.
The spaces $\tB_{pq}^{s}$ are completely independent of $d$,
but for convenience in the definition of $\|f\|_{\tB_{pq}^{s}}$ in $(\ref{def-Besov-space2})$
the smoothness parameter $s$ is normalized as if $\dim M=d$ which, in general, is not the case.
However, if $|B(x, r)|\sim r^d$ uniformly in $x\in M$, like in the classical case on $\bR^d$, then
$\|f\|_{B_{pq}^{s}} \sim \|f\|_{\tB_{pq}^{s}}$.
}
\end{remark}

It will be convenient to introduce (quasi-)norms on
$B_{pq}^{s}$ and $\tB_{pq}^{s}$,
where in the spectral decomposition $2^j$ is replaced by $b^j$ with $b>1$ the constant
from the definition of frames in \S\ref{sec:frames} (see Theorem~\ref{thm:norms}).
Let the real-valued functions
$\Phi_0, \Phi\in C^\infty(\bR_+)$ obey the conditions
\begin{align}
&\supp \Phi_0 \subset   [0, b] ,\;\;
\Phi_0(\lambda) = 1 \hbox{ for }    \lambda \in [0, 1],\;\;
\Phi_0(\lambda) \ge c>0 \;\hbox{ for } \lambda\in [0, b^{3/4}], \label{con_Phi-0}\\
&\supp \Phi \subset   [b^{-1}, b], \;\;
\Phi(\lambda) \ge c>0 \;\hbox{ for } \lambda\in [b^{-3/4}, b^{3/4}],\;\;
\hbox{and}\quad \Phi_0, \Phi \ge 0. \label{con_Phi}
\end{align}
Set $\Phi_j(\lambda):=\Phi(b^{-j}\lambda)$ for $j\ge 1$.
We define
\begin{equation}\label{def-Besov-norm3}
\|f\|_{B_{pq}^{s}(\Phi)} :=
\Big(\sum_{j\ge 0} \Big(b^{s j}
\|\Phi_j(\sqrt{L}) f(\cdot)\|_{\LL^p}
\Big)^q\Big)^{1/q}
\end{equation}
and
\begin{equation}\label{def-Besov-norm4}
\|f\|_{\tB_{pq}^{s}(\Phi)} :=
\Big(\sum_{j\ge 0} \Big(
\||B(\cdot, b^{-j})|^{-s/d}
\Phi_j(\sqrt{L}) f(\cdot)\|_{\LL^p}
\Big)^q\Big)^{1/q}
\end{equation}
with the usual modification when $q=\infty$.


\begin{proposition}\label{prop:Bspq-independ}
For all admissible indices
$\|\cdot\|_{B_{pq}^{s}}$ and $\|\cdot\|_{B_{pq}^{s}(\Phi)}$
are equivalent quasi-norms in $B_{pq}^{s}$,
and $\|\cdot\|_{\tB_{pq}^{s}}$ and $\|\cdot\|_{\tB_{pq}^{s}(\Phi)}$
are equivalent quasi-norms in $\tB_{pq}^{s}$.
Consequently, the definitions of $B_{pq}^{s}$ and $\tB_{pq}^{s}$ are independent
of the particular selection of the functions $\varphi_0, \varphi$ satisfying
$(\ref{cond_phi})$--$(\ref{cond_psi})$.
\end{proposition}

For the proof of this theorem and in the sequel we shall need an~analogue of Peetre's inequality
which involves the maximal operator from (\ref{def:max-op}).


\begin{lemma}\label{lem:Peetre}
Let $t, r>0$ and $\gamma\in\R$. Then there exists a constant $c>0$ such that
for any $g\in\Sigma_t$
\begin{equation}\label{est-Peetre}
\sup_{y\in M}\frac{|B(y, t^{-1})|^{\gamma}|g(y)|}{(1+t\rho(x, y))^{d/r}}
\le c\cM_r\big(|B(\cdot, t^{-1})|^{\gamma}g\big)(x),
\quad x\in M.
\end{equation}
\end{lemma}

\noindent
{\bf Proof.}
Let $g\in\Sigma_t$.
As before, let $\theta \in C^\infty_0(\bR_+)$ and $\theta(\lambda)= 1$ for $\lambda \in [0,1]$.
Denote briefly $\cH_\delta:=\theta(\delta\sqrt L)$ with $\delta=t^{-1}$
and let $\cH_\delta(x, y)$ be its kernel.
Evidently, $\cH_\delta g=g$ and hence
$g(y)=\int_M\cH_\delta(y, z)g(z) d\mu(z)$.
For the kernel $\cH_\delta(\cdot, \cdot)$
we know from Theorem~\ref{thm:main-local-kernels}
that for any $\sigma>0$
\begin{equation}\label{kernel-in-Lip}
|\cH_\delta(y, z)-\cH_\delta(u, z)|
\le c_\sigma \frac{(t\rho(y,u))^\alpha}{|B(y, t^{-1})|(1+t\rho(y, z))^\sigma}
\quad \hbox{if } \rho(y, u)\le t^{-1}.
\end{equation}
Fix $0<\eps<1$. Then for $y\in M$
$$
|g(y)| \le \inf_{u\in B(y, \eps t^{-1})}|g(u)|
+ \sup_{u\in B(y, \eps t^{-1})}|g(y)-g(u)|
$$
and hence
\begin{align*}
G(x) &:= \sup_{y\in M}\frac{|B(y, t^{-1})|^{\gamma}|g(y)|}{(1+t\rho(x, y))^{d/r}}
\le \sup_{y\in M}\frac{|B(y, t^{-1})|^{\gamma}\inf_{u\in B(y, \eps t^{-1})}|g(u)|}
{(1+t\rho(x, y))^{d/r}}\\
& + \sup_{y\in M}\frac{|B(y, t^{-1})|^{\gamma}\sup_{u\in B(y, \eps t^{-1})}|g(y)-g(u)|}
{(1+t\rho(x, y))^{d/r}}
=: G_1(x)+G_2(x).
\end{align*}

To estimate $G_1(x)$ we first observe that
$$
\inf_{u\in B(y, \eps t^{-1})}|g(u)|
\le \Big(
\frac{1}{|B(y, \eps t^{-1})|}\int_{B(y, \eps t^{-1})}|g(u)|^r d\mu(u)
\Big)^{1/r},
$$
which implies
\begin{align}\label{est-G1}
G_1(x) &\le
\left(\frac{|B(x, \rho(x, y)+\eps t^{-1})|}
{|B(y, \eps t^{-1})|(1+t\rho(x, y))^{d}}\right)^{1/r}\notag\\
& \times \left(
\frac{1}{|B(x, \rho(x, y)+\eps t^{-1})|}
\int_{B(y, \eps t^{-1})}\big(|B(y, t^{-1})|^{\gamma}|g(u)|\big)^r d\mu(u)
\right)^{1/r}.
\end{align}
Note that if $u\in B(y, \eps t^{-1})$, then
$B(y, t^{-1})\subset B(u, 2t^{-1})$
and
$B(u, t^{-1})\subset B(y, 2t^{-1})$.
Therefore, the doubling condition (\ref{doubling-0}) yields
$$
c_0^{-1}|B(u, t^{-1})| \le |B(y, t^{-1})| \le c_0|B(u, t^{-1})|,
\quad u\in B(y, \eps t^{-1}).
$$
Also, since $B\big(x, \rho(x, y)+\eps t^{-1}\big) \subset B\big(y, 2\rho(x, y)+\eps t^{-1}\big)$,
then using (\ref{D2})
\begin{align*}
\big|B\big(x, \rho(x, y)+\eps t^{-1}\big)\big|
&\le  \big|B\big(y, 2\rho(x, y)+\eps t^{-1}\big)\big|\\
&\le c_0(1+\eps^{-1}t[2\rho(x, y)+\eps t^{-1}])^d|B(y, \eps t^{-1})|\\
&\le c\eps^{-d}(1+ t\rho(x, y))^d|B(y, \eps t^{-1})|.
\end{align*}
We use the above in (\ref{est-G1}) and enlarge the set of integration
in (\ref{est-G1}) from $B(y, \eps t^{-1})$
to  $B(x, \rho(x, y)+\eps t^{-1})$ to bound $G_1(x)$ by
\begin{align*}
& c \eps^{-d/r}\sup_{y\in M}
\left(\frac{1}{|B(x, \rho(x, y)+\eps t^{-1})|}
\int_{B(x, \rho(x, y)+\eps t^{-1})} \big(|B(u, t^{-1})|^{\gamma}|g(u)|\big)^r d\mu(u)
\right)^{1/r}\\
& \le c \eps^{-d/r} \cM_r\big(|B(\cdot, t^{-1})|^{\gamma}g(\cdot)\big)(x).
\end{align*}
Thus
\begin{equation}\label{est-G1-2}
G_1(x) \le c \eps^{-d/r} \cM_r\big(|B(\cdot, t^{-1})|^{\gamma}g(\cdot)\big)(x).
\end{equation}


We next estimate $G_2(x)$.
Using (\ref{kernel-in-Lip}) we obtain
\begin{align*}
\sup_{u\in B(y, \eps t^{-1})}|g(y)-g(u)|
&\le \sup_{u\in B(y, \eps t^{-1})}
\int_M |\cH_\delta(y, z)- \cH_\delta(u, z)||g(z)| d\mu(z)\\
& \le c\sup_{u\in B(y, \eps t^{-1})} |B(y, t^{-1})|^{-1}
\int_M \frac{(t\rho(y, u))^\alpha |g(z)|}{(1+t\rho(y, z))^\sigma} d\mu(z)\\
&\le c\eps^\alpha |B(y, t^{-1})|^{-1}
\int_M \frac{|g(z)|}{(1+t\rho(y, z))^\sigma} d\mu(z)
\end{align*}
and choosing $\sigma = d/r+d|\gamma|+d+1$ we get
$$ %
G_2(x) \le c \eps^\alpha \sup_{y\in M}
\frac{1}{|B(y, t^{-1})|}
\int_M \frac{|B(y, t^{-1})|^{\gamma}|g(z)|}
{(1+t\rho(y, x))^{\frac{d}{r}}(1+t\rho(y, z))^{\frac{d}{r}+d|\gamma|+d+1}} d\mu(z).
$$ 
Clearly,
$$
(1+t\rho(x, z)) \le (1+t\rho(y, x))(1+t\rho(y, z))
$$
and by (\ref{D2})
$$
c_0^{-1}(1+t\rho(y, z))^{-d}|B(z, t^{-1})|
\le |B(y, t^{-1})| \le c_0(1+t\rho(y, z))^d|B(z, t^{-1})|.
$$
We use these in the above estimate of $G_2(x)$ to obtain
\begin{align*}
G_2(x) &\le c \eps^\alpha \sup_{y\in M}
\frac{1}{|B(y, t^{-1})|}
\int_M \frac{|B(z, t^{-1})|^{\gamma}|g(z)|}
{(1+t\rho(x, z))^{\frac{d}{r}}(1+t\rho(y, z))^{d+1}} d\mu(z)\\
&\le c' \eps^\alpha
\sup_{z\in M}
\frac{|B(z, t^{-1})|^{\gamma}|g(z)|}
{(1+t\rho(x, z))^{\frac{d}{r}}}
\sup_{y\in M}
\frac{1}{|B(y, t^{-1})|}
\int_M \frac{1}
{(1+t\rho(y, z))^{d+1}} d\mu(z)\\
&\le c'' \eps^\alpha G(x),
\end{align*}
where for the last inequality we used (\ref{tech1}).
From this and (\ref{est-G1-2}) we infer
$$
G(x) \le c \eps^{-d/r} \cM_r\big(|B(\cdot, t^{-1})|^{\gamma}g(\cdot)\big)(x) + c'' \eps^\alpha G(x).
$$
Here the constants $c$ and $c''$ are independent of $\eps$.
Consequently, choosing $\eps$ so that  $c'' \eps^\alpha \le 1/2$ we arrive at estimate
(\ref{est-Peetre}).
$\qed$

\medskip

\noindent
{\bf Proof of Proposition~\ref{prop:Bspq-independ}.}
We shall only prove the equivalence of  $\|\cdot\|_{\tB_{pq}^{s}}$ and $\|\cdot\|_{\tB_{pq}^{s}(\Phi)}$.
The proof of the equivalence of $\|\cdot\|_{B_{pq}^{s}}$ and $\|\cdot\|_{B_{pq}^{s}(\Phi)}$
is similar.

It is easy to see (e.g. \cite{F-J-W}) that there exist functions
$\tilde\Phi_0, \tilde\Phi \in C^\infty_0(\bR_+)$
with the properties of $\Phi_0$, $\Phi$ from (\ref{con_Phi-0})-(\ref{con_Phi})
such that
$$
\tilde\Phi_0(\lambda)\Phi_0(\lambda) + \sum_{j\ge 1}\tilde\Phi(b^{-j}\lambda)\Phi(b^{-j}\lambda)=1,
\quad \lambda\in \bR_+.
$$
Set $\tilde\Phi_j(\lambda):=\tilde\Phi(b^{-j}\lambda)$, $j\ge 1$.
Then
$
\sum_{j\ge 0}\tilde\Phi_j(\lambda)\Phi_j(\lambda)=1.
$
By Proposition~\ref{prop:decomp-DD2} it follows that for any $f\in\cD'$
$$
f=\sum_{j\ge 0}\tilde\Phi_j(\sqrt{L})\Phi_j(\sqrt{L})f
\quad\hbox{in}\quad \cD'.
$$
Assume $1<b<2$ (the case $b\ge 2$ is similar) and let $j\ge 1$.
Evidently, there exist $\ell >1$ (depending only on $b$) and $m\ge 1$ such that
$[2^{j-1}, 2^{j+1}] \subset [b^{m-1}, b^{m+\ell+1}]$.
Then $2^j\sim b^m$.
Using the above we have
\begin{align*}
\varphi_j(\sqrt L)f(x)
&= \sum_{\nu=m}^{m+\ell} \varphi_j(\sqrt L)\tilde\Phi_\nu(\sqrt{L})\Phi_\nu(\sqrt{L})f(x)\\
& = \sum_{\nu=m}^{m+\ell}\int_M K_{j\nu}(x, y)\Phi_\nu(\sqrt{L})f(y),
\end{align*}
where $K_{j\nu}(\cdot, \cdot)$ is the kernel of the operator
$\varphi_j(\sqrt L)\tilde\Phi_\nu(\sqrt{L})$.

Choose $0<r<p$ and $\sigma\ge |s|+d/r+3d/2+1$.
By Theorem~\ref{thm:main-local-kernels} we have the following bounds on the kernels of the operators
$\varphi_j(\sqrt L)$ and $\tilde\Phi_\nu(\sqrt{L})$:
$$
|\varphi_j(\sqrt L)(x, y)|\le cD_{2^{-j}, \sigma}(x, y)\le cD_{b^{-\nu}, \sigma}(x, y),
\quad 
|\tilde\Phi_\nu(\sqrt{L})(x, y)|\le cD_{b^{-\nu}, \sigma}(x, y),
$$
and applying (\ref{Comp}) we obtain
$$
|K_{j\nu}(x, y)| \le c'D_{b^{-\nu}, \sigma}(x, y)
\le c|B(x, b^{-\nu})|^{-1}(1+b^\nu\rho(x, y))^{-\sigma+d/2},
\;\; m\le \nu\le m+\ell,
$$
which implies
$$
|\varphi_j(\sqrt L)\tilde\Phi_\nu(\sqrt{L})\Phi_\nu(\sqrt{L})f(x)|
\le \frac{c}{|B(x, b^{-\nu})|}
\int_M\frac{|\Phi_\nu(\sqrt{L})f(y)|}{(1+b^\nu\rho(x, y))^{\sigma-d/2}} d\mu(y).
$$
Observe that $\supp \Phi_\nu \subset [0, b^{\nu+1}]$ and, therefore,
by Proposition~\ref{prop:dec-D1} and Proposition~\ref{prop:dec-D1-inft},
$\Phi_\nu(\sqrt{L})f \in \Sigma_{b^{\nu+1}}$.
Now, using this, (\ref{doubling}) and (\ref{D2}) we get
\begin{align*}
&|B(x, 2^{-j})|^{-s/d}|\varphi_j(\sqrt{L})f(x)|\\
&  \le c\sum_{\nu=m}^{m+\ell}
\frac{1}{|B(x, b^{-\nu})|}
\int_M\frac{|B(x, b^{-\nu})|^{-s/d}|\Phi_\nu(\sqrt{L})f(y)|}
{(1+b^\nu\rho(x, y))^{\sigma-d/2}} d\mu(y)\\
&  \le c\sum_{\nu=m}^{m+\ell}
\frac{1}{|B(x, b^{-\nu})|}
\int_M\frac{|B(y, b^{-\nu})|^{-s/d}|\Phi_\nu(\sqrt{L})f(y)|}
{(1+b^\nu\rho(x, y))^{\sigma-|s|-d/2}} d\mu(y)\\
&  \le c\sum_{\nu=m}^{m+\ell}
\sup_{y\in M}\frac{|B(y, b^{-\nu-1})|^{-s/d}|\Phi_\nu(\sqrt{L})f(y)|}
{(1+b^{\nu+1}\rho(x, y))^{d/r}}
\int_M\frac{|B(x, b^{-\nu})|^{-1}}
{(1+b^\nu\rho(x, y))^{d+1}} d\mu(y)\\
&\le c\sum_{\nu=m}^{m+\ell}
\cM_r\Big( |B(\cdot, b^{-\nu-1})|^{-s/d}\Phi_\nu(\sqrt{L})f(\cdot)\Big)(x).
\end{align*}
Here for the last inequality we used Lemma~\ref{lem:Peetre} and (\ref{tech1}).
Finally, applying the maximal inequality (\ref{max-ineq}) for individual functions ($0<r<p$)
we get
\begin{align*}
\||B(\cdot, 2^{-j})|^{-s/d}\varphi_j(\sqrt{L})f(\cdot)\|_p
&\le c\sum_{\nu=m}^{m+\ell}
\Big\|\cM_r\Big( |B(\cdot, b^{-\nu-1})|^{-s/d}\Phi_\nu(\sqrt{L})f(\cdot)\Big)\Big\|_p\\
&\le c\sum_{\nu=m}^{m+\ell} \||B(\cdot, b^{-\nu})|^{-s/d}\Phi_\nu(\sqrt{L})f(\cdot)\|_p,
\quad j\ge 1.
\end{align*}

Just as above a~similar estimate is proved for $j=0$.
Taking into account that $\ell$ is a constant the~above estimates imply
$\|f\|_{\tB_{pq}^s} \le c\|f\|_{\tB_{pq}^s(\Phi)}$.
In the same manner one proves the~estimate
$\|f\|_{\tB_{pq}^s(\Phi)} \le c\|f\|_{\tB_{pq}^s}$.
$\qed$


\begin{proposition}\label{prop:Bspq-embedding}
The Besov spaces $B_{pq}^s$ and $\tB_{pq}^s$ are quasi-Banach spaces which
are continuously embedded in $\cD'$.
More precisely,
for all admissible indices $s, p, q$,
we have:

$(a)$
If $\mu(M)<\infty$, then
\begin{equation}\label{embed-Bspq-1}
|\langle f, \phi \rangle|
\le c \|f\|_{B_{pq}^s}\cP_{m}^*(\phi),
\quad
f\in B_{pq}^s, \;\; \phi\in\cD,
\end{equation}
provided
$2m > d\big(\frac{1}{\min\{p, 1\}}-1\big)-s$,
and
\begin{equation}\label{embed-tBspq-1}
|\langle f, \phi \rangle|
\le c \|f\|_{\tB_{pq}^s}\cP_{m}^*(\phi),
\quad
f\in \tB_{pq}^s, \;\; \phi\in\cD,
\end{equation}
provided
$2m > \max\big\{0, d\big(\frac{1}{\min\{p, 1\}}-1\big)-s\big\}$.

$(b)$
If $\mu(M)=\infty$, then
\begin{equation}\label{embed-Bspq-2}
|\langle f, \phi \rangle|
\le c \|f\|_{B_{pq}^s}\cP_{m,\ell}^*(\phi),
\quad
f\in B_{pq}^s, \;\; \phi\in\cD,
\end{equation}
provided
$2m > d\big(\frac{1}{\min\{p, 1\}}-1\big)-s$
and $\ell>2d$, and
\begin{equation}\label{embed-tBspq-2}
|\langle f, \phi \rangle|
\le c \|f\|_{\tB_{pq}^s}\cP_{m,\ell}^*(\phi),
\quad
f\in \tB_{pq}^s, \;\; \phi\in\cD,
\end{equation}
provided
$2m > \max\big\{0, d\big(\frac{1}{\min\{p, 1\}}-1\big)-s\big\}$
and $\ell>\max\big\{2d, \big|d\big(\frac{1}{p}-1\big)\big|+|s|\big\}$.
\end{proposition}

\noindent
{\bf Proof.}
Observe first that the~completeness of $B_{pq}^s$ and $\tB_{pq}^s$ follows readily by
the continuous embedding of $B_{pq}^s$ and $\tB_{pq}^s$ in $\cD'$.
We shall only prove the continuous embedding of $\tB_{pq}^s$ in $\cD'_{m,\ell}$
in the case when $\mu(M)=\infty$. All other cases are easier and we skip the details.

Choose real-valued functions $\varphi_0, \varphi\in C^\infty_0(\bR_+)$ so that
$\supp \varphi_0\subset [0, 2]$, $\varphi_0(\lambda)=1$ for $\lambda\in[0, 1]$,
$\supp \varphi\subset [2^{-1}, 2]$, and
$\varphi_0^2(\lambda)+\sum_{j\ge 1} \varphi^2(2^{-j}\lambda)=1$ for $\lambda\in \R_+$.
Set $\varphi_j(\lambda):=\varphi(2^{-j}\lambda)$ for $j\ge 1$.
Then $\sum_{j\ge 0} \varphi_j^2(\lambda)=1$ for $\lambda\in \bR_+$
and hence, using Proposition~\ref{prop:decomp-DD2}, for any $f\in\cD'$
\begin{equation}\label{rep-f-B}
f=\sum_{j\ge 0} \varphi_j^2(\sqrt L)f
\quad \hbox{in $\cD'$.}
\end{equation}
Also, observe that $\{\varphi_j\}_{j\ge 0}$ are just like the functions
in the definition of $\tB_{pq}^s$
and can be used to define an equivalent norm on $\tB_{pq}^s$ as in (\ref{def-Besov-norm4}).
From (\ref{rep-f-B}) we get
\begin{equation}\label{embedB-0}
\langle f, \phi \rangle
= \sum_{j\ge 0} \langle \varphi_j^2(\sqrt L)f, \phi \rangle
= \sum_{j\ge 0} \big\langle \varphi_j(\sqrt L)f, \varphi_j(\sqrt L)\phi \big\rangle,
\quad \phi\in \cD.
\end{equation}

We next estimate $|\varphi_{j}(\sqrt L)\phi(x)|$, $j\ge 1$.
To this end we set
$\omega(\lambda):= \lambda^{-2m}\varphi(\lambda)$.
Then $\varphi_{j}(\sqrt\lambda)= 2^{-2mj}\omega(2^{-j}\sqrt \lambda)\lambda^{2m}$ and hence
$$
\varphi_{j}(\sqrt L)\phi(x)
= 2^{-2mj}\int_M \omega(2^{-j}\sqrt L)(x, y)L^m\phi(y)d\mu(y).
$$
Clearly, $\omega\in C^\infty_0(\bR_+)$ and $\supp \omega\subset [1/2, 2]$.
Therefore, by Theorem~\ref{thm:main-local-kernels}
$$
|\omega(2^{-j}\sqrt L)(x, y)|
\le c|B(y, 2^{-j})|^{-1}(1+2^j\rho(x, y))^{-\ell},
$$
where $\ell > 2d$ is from the assumption in (b).
On the other hand, since $\phi\in \cD$
we have
$|L^m\phi(x)| \le c(1+\rho(x, x_0))^{-\ell}\cP_{m, \ell}(\phi)$.
From the above we obtain
\begin{align}\label{embed-B3}
|\varphi_{j}(\sqrt L)\phi(x)|
&\le c2^{-2mj}\cP_{m, \ell}(\phi)
\int_M \frac{d\mu(y)}
{|B(y,2^{-j})|(1+2^j\rho(x, y))^{\ell}(1+\rho(y, x_0))^{\ell}} \notag\\
&\le c2^{-2mj}\cP_{m, \ell}(\phi)(1+\rho(x, x_0))^{-\ell}.
\end{align}
Here for the last inequality we used (\ref{tech5}) and that $\ell > 2d$.


We are now prepared to estimate the inner products in (\ref{embedB-0}).
We consider two cases:

Case 1: $1\le p \le \infty$. Then applying H\"{o}lder's inequality ($1/p+1/p'=1$) we get

\begin{equation}\label{inner-prod-B1}
\begin{aligned}
& \big|\big\langle \varphi_j(\sqrt L)f, \varphi_j(\sqrt L)\phi \big\rangle\big|\\
& \qquad
\le \int_M |B(x, 2^{-j})|^{-s/d} |\varphi_j(\sqrt L)f(x)|
|B(x, 2^{-j})|^{s/d}|\varphi_j(\sqrt L)\phi(x)| d\mu(x)\\
& \qquad
\le \||B(x, 2^{-j})|^{-s/d} \varphi_j(\sqrt L)f\|_p
\||B(x, 2^{-j})|^{s/d} \varphi_j(\sqrt L)\phi\|_{p'}\\
& \qquad
\le c\|f\|_{\tB_{pq}^s}
\||B(x, 2^{-j})|^{s/d} \varphi_j(\sqrt L)\phi\|_{p'},
\quad j\ge 0.
\end{aligned}
\end{equation}
Using (\ref{embed-B3}) we obtain for $j\ge 1$
$$
Q:=\||B(x, 2^{-j})|^{s/d} \varphi_j(\sqrt L)\phi\|_{p'}^{p'}
\le c2^{-2mjp'}\cP_{m, \ell}(\phi)^{p'}
\int_M\frac{|B(x, 2^{-j})|^{sp'/d}}{(1+\rho(x, x_0))^{\ell p'}}d\mu(x).
$$
Two cases present themselves here depending on whether $s\ge 0$ or $s<0$.
If $s\ge 0$, then by (\ref{D2}) we have
$|B(x, 2^{-j})| \le |B(x, 1)| \le c_0(1+\rho(x, x_0))^d|B(x_0, 1)|$
and hence
\begin{align*}
Q
&\le c2^{-2mjp'}\cP_{m, \ell}(\phi)^{p'}|B(x_0, 1)|^{sp'/d}
\int_M\frac{d\mu(x)}{(1+\rho(x, x_0))^{(\ell-s)p'}}\\
&\le c2^{-2mjp'}\cP_{m, \ell}(\phi)^{p'}|B(x_0, 1)|^{sp'/d+1},
\end{align*}
where for the last inequality we used (\ref{tech1}) and that $(\ell-s)p'>d$,
which follows from $\ell>|d(1/p-1)|+|s|$.
In the case $s<0$ we use that
$$
|B(x_0,1)| \le c_0(1+\rho(x, x_0))^d |B(x, 1)| \le c_0^22^{jd}(1+\rho(x, x_0))^d|B(x, 2^{-j})|,
$$
which is immediate from (\ref{doubling}),(\ref{D2}), to obtain
\begin{align*}
Q &\le c2^{-j(2m+s)p'}\cP_{m, \ell}(\phi)^{p'}|B(x_0, 1)|^{sp'/d}
\int_M\frac{d\mu(x)}{(1+\rho(x, x_0))^{(\ell+s)p'}}\\
&\le c2^{-j(2m+s)p'}\cP_{m, \ell}(\phi)^{p'}|B(x_0, 1)|^{sp'/d+1}.
\end{align*}
Here we again used (\ref{tech1}) and that $(\ell+s)p'>d$ due to $\ell>|d(1/p-1)|+|s|$.
From the above estimates on $Q$ and (\ref{inner-prod-B1}) we get for $j\ge 1$
\begin{equation}\label{embed-B5}
\big|\big\langle \varphi_{j}(\sqrt L)f, \varphi_{j}(\sqrt L)\phi \big\rangle\big|
\le c2^{-j(2m+\min\{s, 0\})}|B(x_0, 1)|^{s/d+1-1/p}\|f\|_{\tB_{pq}^s}\cP_{m,\ell}^*(\phi).
\end{equation}


It remains to consider the easier case when $j=0$.
Applying Theorem~\ref{thm:main-local-kernels} to $\varphi_0$ and since $\phi\in\cD$,
we obtain
$$
|\varphi_0(\sqrt L)(x, y)|
\le c|B(y,1)|^{-1}(1+\rho(x, y))^{-\ell}
\quad\hbox{and}\quad
|\phi(x)| \le c(1+\rho(x, x_0))^{-\ell}\cP_{0, \ell}(\phi),
$$
which as in (\ref{embed-B3}) imply
$
|\varphi_0(\sqrt L)\phi(x)|
\le c(1+\rho(x, x_0))^{-\ell}\cP_{0, \ell}(\phi).
$
As above this leads to
$$
\||B(\cdot, 1)|^{s/d} \varphi_j(\sqrt L)\phi\|_{p'}
\le c|B(x_0, 1)|^{s/d+1-1/p}\cP_{0,\ell}(\phi).
$$
In turn, this and (\ref{inner-prod-B1}) yield (\ref{embed-B5}) with $m=0$ for $j=0$.

Summing up estimates (\ref{embed-B5}),
taking into account (\ref{embedB-0}) and that $2m >\max\{0, -s\}$,
we arrive at (\ref{embed-tBspq-2}).


\smallskip

Case 2: $0<p<1$.
Setting $\gamma:= s/d-1/p+1$, we have for $\phi\in \cD$ and $j\ge 1$
$$
\big|\big\langle \varphi_j(\sqrt L)f, \varphi_j(\sqrt L)\phi \big\rangle\big|
\le \||B(x, 2^{-j})|^{-\gamma} \varphi_j(\sqrt L)f\|_1
\||B(x, 2^{-j})|^{\gamma} \varphi_j(\sqrt L)\phi\|_\infty.
$$
Since $\varphi_j(\sqrt L)f\in \Sigma_{2^{j+1}}$,
Proposition~\ref{prop:Nikolski} yields
\begin{align*} 
\||B(\cdot, 2^{-j})|^{-\gamma} \varphi_j(\sqrt L)f\|_1
&\le c \||B(\cdot, 2^{-j})|^{-\gamma+1-1/p} \varphi_j(\sqrt L)f\|_p\notag\\
&=c \||B(\cdot, 2^{-j})|^{-s/d} \varphi_j(\sqrt L)f\|_p
\le c\|f\|_{\tB_{pq}^s}.
\end{align*}
On the other hand by (\ref{embed-B3})
$$
R:=\||B(\cdot, 2^{-j})|^{\gamma} \varphi_j(\sqrt L)\phi\|_\infty
\le c2^{-2mj}\cP_{m, \ell}(\phi)
\sup_{x\in M}\frac{|B(x, 2^{-j})|^{\gamma}}{(1+\rho(x, x_0))^{\ell}}.
$$
As in the estimation of $Q$ above we obtain
$R \le c2^{-2mj}\cP_{m, \ell}(\phi)|B(x_0,1)|^{\gamma}$ if $\gamma\ge 0$ and
$R \le c2^{-j(2m+d\gamma)}\cP_{m, \ell}(\phi)|B(x_0,1)|^{\gamma}$ if $\gamma<0$.
Here we used that $\ell >d|\gamma|$ due to
$\ell>|d(1/p-1)|+|s|$.
Therefore, for $ j\ge 1$
$$
\big|\big\langle \varphi_{j}(\sqrt L)f, \varphi_{j}(\sqrt L)\phi \big\rangle\big|
\le c2^{-j(2m+\min\{0, d\gamma\})}|B(x_0, 1)|^{s/d+1-1/p}\|f\|_{\tB_{pq}^s}\cP_{m,\ell}(\phi).
$$
Now we complete the proof of (\ref{embed-tBspq-2}) just as in the case $1\le p\le \infty$,
taking into account that $2m > -\min\{0, d\gamma\}=\max\{0, d(1/p-1)-s\}$.
$\qed$

\subsection{Heat kernel characterization of Besov spaces}\label{sec:heat-dec-B-spaces}

We shall show that the Besov spaces $B_{pq}^{s}$ and $\tB_{pq}^{s}$ can be equivalently
defined using directly the Heat kernel when $p$ is restricted to $1\le p\le \infty$.


\begin{definition}\label{def:BH-norms}
Given $s \in \R$, let $m$ be the smallest $m\in \bZ_+$ such that $m>s$.
We define
\begin{align*}
\|f\|_{B_{pq}^{s}(H)}
& := \|e^{-L}f\|_p
+ \Big(\int_0^1\Big[t^{-s/2}\|(tL)^{m/2} e^{-tL}f\|_p\Big]^q \frac{dt}{t}
\Big)^{1/q},\\ 
\|f\|_{\tB_{pq}^{s}(H)}
& := \||B(\cdot, 1)|^{-s/d}e^{-L}f\|_p
+ \Big(\int_0^1\big\||B(\cdot, t^{1/2})|^{-s/d}(tL)^{m/2} e^{-tL}f\big\|_p^q \frac{dt}{t}
\Big)^{1/q} 
\end{align*}
with the usual modification when $q=\infty$.
\end{definition}


\begin{theorem}\label{thm:Heat-B-charac}
Suppose $s \in \R$, $1\le p\le \infty$, $0< q \le \infty$,
and $m>s$, $m\in \bZ_+$, as in the above definition.
Let $f\in \cD'$. Then we have:

$(a)$
$f\in B_{pq}^{s}$ if and only if $e^{-L}f\in \LL^p$
and $\|f\|_{B_{pq}^{s}(H)} < \infty$.
Moreover, if $f\in B_{pq}^{s}$, then
$\|f\|_{B_{pq}^{s}}\sim \|f\|_{B_{pq}^{s}(H)}$.

$(b)$
$f\in \tB_{pq}^{s}$ if and only if $|B(\cdot, 1)|^{s/d}e^{-L}f\in \LL^p$
and $\|f\|_{\tB_{pq}^{s}(H)} < \infty$.
Moreover, if $f\in \tB_{pq}^{s}$, then
$\|f\|_{\tB_{pq}^{s}}\sim \|f\|_{\tB_{pq}^{s}(H)}$.

\end{theorem}

\noindent
{\bf Proof.}
We shall only prove Part (b); the proof of Part (a) is easier and will be omitted.
Let $\varphi_0$, $\varphi$, and $\varphi_j$, $j\ge 1$, be precisely as in the proof of
Proposition~\ref{prop:Bspq-embedding}.
Then for any $f\in \cD'$ we have $f=\sum_{j\ge 0} \varphi_j^2(\sqrt L)f$ and hence
$$
|B(\cdot, t^{1/2})|^{-s/d}(tL)^{m/2} e^{-tL}f
= \sum_{j\ge 0} |B(\cdot, t^{1/2})|^{-s/d}(tL)^{m/2} e^{-tL}\varphi_j^2(\sqrt L)f
=: \sum_{j\ge 0}F_j.
$$
It is readily seen that for $j\ge 1$
\begin{align*}
F_{j}
&= |B(\cdot, t^{1/2})|^{-s/d}(tL)^{m/2} e^{-tL}\varphi(2^{-j}\sqrt L)\varphi(2^{-j}\sqrt L)f\\
&= |B(\cdot, t^{1/2})|^{-s/d}\omega(2^{-j}\sqrt L)\varphi(2^{-j}\sqrt L)f,
\end{align*}
where
$\omega(\lambda):=(t\lambda^24^{j})^{m/2}e^{-t\lambda^24^{j}}\varphi(\lambda)$.
As
$\varphi\in C^\infty$, $\supp \varphi \subset [\frac 12, 2]$, and $0\le \varphi \le 1$
we have, by Theorem~\ref{thm:main-local-kernels},
\begin{equation}\label{omega-bound}
|\omega(2^{-j}\sqrt L)(x, y)|
\le c_\sigma (t4^{j})^{m/2}e^{-t4^{j}} |B(x, 2^{-j})|^{-1}(1+2^j\rho(x, y))^{-\sigma},
\;\forall \sigma>0.
\end{equation}
On the other hand, (\ref{doubling}) and (\ref{D2}) easily imply
\begin{equation}\label{est-from-doubl}
|B(x, t^{1/2})|^{-s/d}
\le c\big(1+ (t4^j)^{-s/2}\big)(1+2^j\rho(x, y))^{|s|}
|B(y, 2^{-j})|^{-s/d}.
\end{equation}
To this end one has to consider four cases, depending on whether $t^{1/2} \le 2^{-j}$ or $t^{1/2} > 2^{-j}$
and whether $s\ge 0$ or $s<0$.
Combining the above with (\ref{omega-bound}) we obtain
\begin{equation}\label{est-Fj-B}
|F_{j}(x)| \le c\big(1+ (t4^j)^{-\frac{s}{2}}\big)(t4^{j})^{\frac{m}{2}}e^{-t4^{j}}
\int_M \frac{|B(y, 2^{-j})|^{-\frac{s}{d}}|\varphi(2^{-j}\sqrt L)f(y)|}
{|B(x, 2^{-j})|(1+2^j\rho(x, y))^{\sigma-|s|}} d\mu(y).
\end{equation}
We now choose $\sigma \ge |s|+d+1$ and invoke Proposition~\ref{prop:young} to obtain
$$
\|F_{j}\|_p \le c \big(1+ (t4^j)^{-s/2}\big)(t4^{j})^{m/2}e^{-t4^{j}}
\||B(\cdot, 2^{-j})|^{-\frac{s}{d}}\varphi_{j}(\sqrt L)f\|_p,
\quad j\ge 1.
$$
One similarly obtains the estimate
$$
\|F_0\|_p \le c\||B(\cdot, 1)|^{-\frac{s}{d}}\varphi_0(\sqrt L)f\|_p.
$$
Putting the above estimates together we obtain for $0<t\le 1$
\begin{align*}
&\||B(\cdot, t^{1/2})|^{-s/d}(tL)^{m/2} e^{-tL}f\|_p\\
&\qquad\qquad\qquad
\le c\sum_{j\ge 0} \big[(t4^{j})^{m/2}+ (t4^j)^{(m-s)/2}\big]e^{-t4^{j}}
\||B(\cdot, 2^{-j})|^{-\frac{s}{d}}\varphi_{j}(\sqrt L)f\|_p.
\end{align*}
Let
$h_j(t):= \big[(t4^{j})^{m/2}+ (t4^j)^{(m-s)/2}\big]e^{-t4^{j}}$
and $b_j:= \||B(\cdot, 2^{-j})|^{-\frac{s}{d}}\varphi_{j}(\sqrt L)f\|_p$.
Then from above
\begin{align*}
&\Big(\int_0^1\big\||B(\cdot, t^{1/2})|^{-s/d}(tL)^{m/2} e^{-tL}f\big\|_p^q \frac{dt}{t}\Big)^{1/q}
\le c\Big(\int_0^1
\Big(\sum_{j\ge 0}h_j(t) b_j\Big)^q
\frac{dt}{t}\Big)^{1/q}\\
&= c\Big(\sum_{\nu\ge 0}\int_{4^{-\nu-1}}^{4^{-\nu}}
\Big(\sum_{j\ge 0}h_j(t) b_j\Big)^q \frac{dt}{t}\Big)^{1/q}
\le c\Big(\sum_{\nu\ge 0}
\Big(\sum_{j\ge 0}a_{j-\nu} b_j\Big)^q\Big)^{1/q}.
\end{align*}
Here
$$
a_{j-\nu} = \max \{h_j(t): t\in [4^{-\nu-1}, 4^{-\nu}]\}
\le (4^{(j-\nu)m/2}+ 4^{(j-\nu)(m-s)/2})e^{-4^{j-\nu-1}}
$$
and we set
$a_\nu := (4^{\nu m/2}+ 4^{\nu(m-s)/2})e^{-4^{\nu-1}}$, $\nu\in \bZ$.

Three cases present themselves here, depending on whether $q=\infty$, $1<q<\infty$ or $0<q<1$.
The case when $q=\infty$ is obvious.
If~$1<q<\infty$, we apply Young's inequality to the convolution of the above sequences to obtain
$$
\Big(\sum_{\nu\ge 0}\Big(\sum_{j\ge 0}a_{j-\nu} b_j\Big)^q\Big)^{1/q}
\le \sum_{\nu\in\bZ} a_\nu \Big(\sum_{j\ge 0} b_j^q\Big)^{1/q}
\le c\Big(\sum_{j\ge 0} b_j^q\Big)^{1/q},
$$
where we used that
$\sum_{\nu\in\bZ} a_\nu \le c$ due to $m>s$.
If $0<q\le 1$, we apply the~$q$-triangle inequality and obtain
$$
\sum_{\nu\ge 0}\Big(\sum_{j\ge 0}a_{j-\nu} b_j\Big)^q
\le \sum_{\nu\ge 0}\sum_{j\ge 0}a_{j-\nu}^q b_j^q
\le \sum_{\nu\in\bZ} a_\nu^q \sum_{j\ge 0} b_j^q
\le c\sum_{j\ge 0} b_j^q.
$$
Here we used that
$\sum_{\nu\in\bZ} a_\nu^q \le c$.
In both cases, we get
$$
\Big(\int_0^1\big\||B(\cdot, t^{1/2})|^{-s/d}(tL)^{m/2} e^{-tL}f\big\|_p^q \frac{dt}{t}\Big)^{1/q}
\le c\Big(\sum_{j\ge 0} b_j^q\Big)^{1/q}
\le c\|f\|_{\tB_{pq}^{s}}.
$$
It is easier to show that
$
\||B(\cdot, 1)|^{-s/d}e^{-L}f\|_p \le c\|f\|_{\tB_{pq}^{s}}.
$
The proof follows in the footsteps of the above proof and will be omitted.
Combining the above two estimates we get
$\|f\|_{\tB_{pq}^{s}(H)} \le c\|f\|_{\tB_{pq}^{s}}$.


We next prove an estimate in the opposite direction.
We only consider the case when $0<q<\infty$; the case $q=\infty$ is easier.
Assume that $\varphi_0$, $\varphi$, and $\varphi_j$, $j\ge 1$, are as in the definition
of $\tB^s_{pq}$ (Definition~\ref{def-B-spaces}).
We can write
\begin{align*}
&|B(x, 2^{-j})|^{-s/d}\varphi_{j}(\sqrt L)f(x)\\
& \qquad\qquad\qquad
= |B(x, 2^{-j})|^{-s/d}(tL)^{-m/2} e^{tL}\varphi(2^{-j}\sqrt L)(tL)^{m/2} e^{-tL}f(x)\\
&  \qquad\qquad\qquad
= |B(x, 2^{-j})|^{-s/d}\omega(2^{-j}\sqrt L)(tL)^{m/2} e^{-tL}f(x),
\end{align*}
where $\omega(\lambda):= (t\lambda^24^j)^{-m/2}e^{-t\lambda^24^j}\varphi(\lambda)$
and $t\in [4^{-j}, 4^{-j+1}]$.
Since $\supp \varphi \subset [1/2, 2]$ we have $\|\omega\|_\infty \le c$ and by
Theorem~\ref{thm:main-local-kernels}
$$
|\omega(2^{-j}\sqrt L)(x, y)| \le c_\sigma |B(x, 2^{-j})|^{-1}\big(1+2^j\rho(x, y)\big)^{-\sigma}.
$$
From (\ref{D2})
$
|B(x, 2^{-j})|^{-s/d} \le c\big(1+2^j\rho(x, y)\big)^{|s|}|B(y, t^{1/2})|^{-s/d},
$
$t\in [4^{-j}, 4^{-j+1}]$.
Therefore,
\begin{equation}\label{est-Bphij}
|B(x, 2^{-j})|^{-s/d}|\varphi_{j}(\sqrt L)f(x)|
\le c\int_M\frac{|B(y, t^{1/2})|^{-s/d}|(tL)^{m/2} e^{-tL}f(y)|}
{|B(x, 2^{-j})|\big(1+2^j\rho(x, y)\big)^{\sigma-|s|}} d\mu(y).
\end{equation}
Choosing $\sigma \ge |s|+d+1$ and applying Proposition~\ref{prop:young}, we get
$$
\big\||B(\cdot, 2^{-j})|^{-s/d}\varphi_{j}(\sqrt L)f\big\|_p
\le c \big\||B(\cdot, t^{1/2})|^{-s/d}(tL)^{m/2} e^{-tL}f\big\|_p
$$
for $t\in [4^{-j}, 4^{-j+1}]$ and hence for $j\ge 1$
$$
\big\||B(\cdot, 2^{-j})|^{-s/d}\varphi_{j}(\sqrt L)f\big\|_p^q
\le c \int_{4^{-j}}^{4^{-j+1}}\big\||B(\cdot, t^{1/2})|^{-s/d}(tL)^{m/2} e^{-tL}f\big\|_p^q\frac{dt}{t}.
$$
Also, one easily obtains
$$
\||B(\cdot, 1)|^{-s/d}\varphi_0(\sqrt L)f\|_p
\le c \||B(\cdot, 1)|^{-s/d}e^{-L}f\|_p.
$$
Summing up the former estimates for $j=1, 2, \dots$ and using the result
and the latter estimate in the definition of $\|f\|_{\tB_{pq}^{s}}$
we get
$\|f\|_{\tB_{pq}^{s}} \le c\|f\|_{\tB_{pq}^{s}(H)}$.
$\qed$


\begin{remark}\label{rem:heat-kernel-Besov}
{\rm
From the above proof it easily follows that whenever $f$ is a function the terms
$\|e^{-L}f\|_p$ and
$\||B(\cdot, 1)|^{-s/d}e^{-L}f\|_p$
in Definition~\ref{def:BH-norms} can be replaced by
$\|f\|_p$ and
$\||B(\cdot, 1)|^{-s/d}f\|_p$, respectively.

Also, observe that in the case when $s>0$ Theorem~\ref{thm:Heat-B-charac} (a) follows readily from
the characterization of $B^s_{pq}$ by means of linear approximation from $\Sigma_t^p$,
see \cite{CKP}, $\S 6.1$.
}
\end{remark}

\subsection{Frame decomposition of Besov spaces}\label{sec:frame-dec-B-spaces}

Our primary goal here is to show that the Besov spaces introduced by Definition~\ref{def-B-spaces}
can be characterized in terms of respective sequence norms of the frame coefficients
of distributions, using the frames from \S\ref{sec:frames}.

Everywhere in this subsection $\{\psi_\xi\}_{\xi\in\cX}$, $\{\tilde\psi_\xi\}_{\xi\in\cX}$
will be the pair of dual frames from \S\S\ref{natural-frame}-\ref{dual-frame},
$\cX:= \cup_{j\ge 0} \cX_j$ will denote the sets of the centers of the frame elements and
$\{A_\xi\}_{\xi\in\cX_j}$ will be the associated partitions of $M$.

Our first order of business is to introduce the sequence spaces  $b_{pq}^s$ and $\tb_{pq}^s$.


\begin{definition}\label{def:b-spaces}
Let $s\in \R$ and $0<p,q \le \infty$.

$(a)$
$b_{pq}^s$
is defined as the space of all complex-valued sequences
$a:=\{a_{\xi}\}_{\xi\in \cX}$ such that
\begin{equation}\label{def-tilde-berpq}
\|a\|_{b_{pq}^s}
:=\Bigl(\sum_{j\ge 0}b^{jsq}
\Bigl[\sum_{\xi\in \cX_j}\Big(|B(\xi, b^{-j})|^{1/p-1/2}|a_\xi|\Big)^p
\Bigr]^{q/p}\Bigr)^{1/q} <\infty.
\end{equation}

$(b)$
$\tb_{pq}^s$
is defined as the space of all complex-valued sequences
$a:=\{a_{\xi}\}_{\xi\in \cX}$ such that
\begin{equation}\label{def-berpq}
\|a\|_{\tb_{pq}^s}
:=\Bigl(\sum_{j\ge 0}
\Bigl[\sum_{\xi\in \cX_j}\Big(|B(\xi, b^{-j})|^{-s/d+1/p-1/2}|a_\xi|\Big)^p
\Bigr]^{q/p}\Bigr)^{1/q} <\infty.
\end{equation}
Above as usual the $\ell^p$ or $\ell^q$ norm is replaced by the $\sup$-norm if
$p=\infty$ or $q=\infty$.
\end{definition}

In our further analysis we shall use the ``analysis" and ``synthesis" operators
defined by
\begin{equation}\label{anal_synth_oprts2}
S_{\tilde\psi}: f\rightarrow \{\langle f, \tilde\psi_\xi\rangle\}_{\xi \in \cX}
\quad\text{and}\quad
T_{\psi}: \{a_\xi\}_{\xi \in \cX}\rightarrow \sum_{\xi\in \cX}a_\xi\psi_\xi.
\end{equation}
Here the roles of $\{\psi_\xi\}$ and $\{\tilde\psi_\xi\}$ can be interchanged.


\begin{theorem}\label{thm:B-character}
Let $s \in \R$ and  $0< p,q\le \infty$.
%
$(a)$
The operators
$S_{\tilde\psi}: B_{pq}^s \rightarrow  b_{pq}^s$ and
$T_{\psi}: b_{pq}^s \rightarrow B_{pq}^s$
are bounded and $T_{\psi}\circ S_{\tilde\psi}=Id$ on $B_{pq}^s$.
Consequently,  for $f\in \cD'$ we have $f\in B_{pq}^s$ if and
only if $\{\langle f, \tilde\psi_\xi\rangle\}_{\xi \in \cX}\in b_{pq}^s$.
Moreover, if $f\in B_{pq}^s$, then
$
\|f\|_{B_{pq}^s} \sim  \|\{\langle f,\tilde\psi_\xi\rangle\}\|_{b_{pq}^s}
$
and under the reverse doubling condition $(\ref{reverse-doubling})$
\begin{equation}\label{Bnorm-equivalence-1}
\|f\|_{B_{pq}^s} %
\sim \Big(\sum_{j\ge 0} b^{jsq}\Bigl[\sum_{\xi\in \cX_j}
\|\langle f,\tilde\psi_{\xi}\rangle\psi_{\xi}\|_p^p\Bigr]^{q/p}\Bigr)^{1/q}.
\end{equation}

\noindent
$(b)$
The operators
$S_{\tilde\psi}: \tB_{pq}^s \rightarrow  \tb_{pq}^s$ and
$T_{\psi}: \tb_{pq}^s \rightarrow \tB_{pq}^s$
are bounded and $T_{\psi}\circ S_{\tilde\psi}=Id$ on $\tB_{pq}^s$.
Hence,
$f\in \tB_{pq}^s \Longleftrightarrow \{\langle f, \tilde\psi_\xi\rangle\}_{\xi \in \cX}\in \tb_{pq}^s$.
Furthermore, if $f\in \tB_{pq}^s$, then
$\|f\|_{B_{pq}^s} \sim  \|\{\langle f,\tilde\psi_\xi\rangle\}\|_{b_{pq}^s}
$
and under the reverse doubling condition $(\ref{reverse-doubling})$
\begin{equation}\label{tBnorm-equiv-1}
\|f\|_{\tB_{pq}^s} 
\sim \Big(\sum_{j\ge 0} \Bigl[\sum_{\xi\in \cX_j}
\Big(|B(\xi, b^{-j})|^{-s/d}
\|\langle f,\tilde\psi_{\xi}\rangle\psi_{\xi}\|_{p}\Big)^p\Bigr]^{q/p}\Bigr)^{1/q}.
\end{equation}
Above in $(a)$ and $(b)$ the roles of $\{\psi_\xi\}$ and $\{\tilde\psi_\xi\}$ can be interchanged.
\end{theorem}

To prove this theorem we need some technical results which will be presented next.


\begin{definition}\label{def-h-star}
For any set of complex numbers
$\{a_{\xi}\}_{\xi\in \cX_j}$ $(j\ge 0)$
we define
\begin{equation}\label{def.h-star}
a^{\ast}_\xi:=\sum_{\eta\in\cX_j}\frac{|a_\eta|}{(1+b^j\rho(\eta,\xi))^\tau}
\quad\hbox{for}\quad \xi\in \cX_j,
\end{equation}
where $\tau >1$ is a sufficiently large constant that will be selected later on.
\end{definition}


\begin{lemma}\label{lem:disc-max}
Let $0<r<1$ and assume that $\tau$ in the definition $(\ref{def.h-star})$ of $a_\xi^*$ obeys
$\tau > d/r$.
Then for any set of complex numbers
$\{a_\xi\}_{\xi\in \cX_j}$ $(j\ge 0)$ we have
\begin{equation}\label{disc-max}
\sum_{\xi\in\cX_j} a_\xi^{\ast}\ONE_{A_\xi}(x)
\le c \cM_r\Big(\sum_{\eta\in \cX_j}|a_\eta|\ONE_{A_\eta}\Big)(x),
\quad x\in M.
\end{equation}
\end{lemma}

\noindent
{\bf Proof.}
Fix $\xi\in\cX_j$ and set
$
S_0:=\{\eta\in\cX_j: \rho(\xi, \eta) \le c^\diamond b^{-j}\}
$
and
$$
S_m:=\{\eta\in\cX_j: c^\diamond b^{-j+m-1} < \rho(\xi, \eta) \le c^\diamond b^{-j+m}\},
$$
where $c^\diamond:=\gamma b^{-1}$ with $\gamma$ being the constant from the construction of the frames
in \S\ref{natural-frame}.
Let
$
B_m:=B(\xi, c^\diamond(b^m+1)b^{-j}).
$
Note that
$A_\eta\subset B_m$ if $\eta\in S_\ell$, $0\le \ell\le m$,
and hence, using (\ref{doubling}),
\begin{equation}\label{est-B-A}
\frac{|B_m|}{|A_\eta|}
\le \frac{|B(\eta, 2c^\diamond(b^m+1)b^{-j})|}{|B(\eta, c^\diamond 2^{-1}b^{-j-1})|}
\le cb^{dm}.
\end{equation}
We have
\begin{align*}
a_\xi^* \le c\sum_{m\ge 0}b^{-m\tau} \sum_{\eta\in S_m} |a_\eta|
\le c\sum_{m\ge 0}b^{-m\tau} \Big(\sum_{\eta\in S_m} |a_\eta|^r\Big)^{1/r}
\end{align*}
and using (\ref{est-B-A})
\begin{align*}
a_\xi^*
& \le c\sum_{m\ge 0} b^{-m\tau}
\Big(\int_M\Big[\sum_{\eta\in S_m} |a_\eta| |A_\eta|^{-1/r}\ONE_{A_\eta}\Big]^r d\mu(y)\Big)^{1/r}\\
&\le c\sum_{m\ge 0} b^{-m\tau}\Big(\frac{1}{|B_m|}\int_{B_m}
\Big[\sum_{\eta\in S_m} \Big(\frac{|B_m|}{|A_\eta|}\Big)^{1/r}|a_\eta|
\ONE_{A_\eta}\Big]^r d\mu(y)\Big)^{1/r}\\
&\le c\sum_{m\ge 0} b^{-m(\tau-d/r)}\Big(\frac{1}{|B_m|}\int_{B_m}
\Big[\sum_{\eta\in S_m}|a_\eta|
\ONE_{A_\eta}\Big]^r d\mu(y)\Big)^{1/r}\\
& \le c\cM_r\Big(\sum_{\eta\in \cX_j}|a_\eta|\ONE_{A_\eta}\Big)(x)
\sum_{m\ge 0} b^{-m(\tau-d/r)}
\le c\cM_r\Big(\sum_{\eta\in \cX_j}|a_\eta|\ONE_{A_\eta}\Big)(x),
\end{align*}
for $x\in A_\xi$,
which confirms (\ref{disc-max}).
$\qed$

\medskip


\begin{lemma}\label{lem:half_shannon}
Let $0<p < \infty$ and $\gamma\in \R$.
Then for any $g\in \Sigma_{b^{j+2}}$, $j\ge 0$,
\begin{equation}\label{half_shannon}
 \Big(\sum_{\xi\in \cX_j}|B(\xi, b^{-j})|^{\gamma p}
 \sup_{x\in A_\xi}|g(x)|^p |A_\xi|\Big)^{1/p}
 \le c\||B(\cdot, b^{-j})|^\gamma g(\cdot)\|_{\LL^p}.
\end{equation}
\end{lemma}

\noindent
{\bf Proof.}
Let $0<r<p$. We have
\begin{align*}
&\sum_{\xi\in \cX_j}|B(\xi, b^{-j})|^{\gamma p}
\sup_{x\in A_\xi}|g(x)|^p |A_\xi|\\
& \qquad\qquad\qquad \le c \int_M \sum_{\xi\in \cX_j}
\sup_{x\in A_\xi}\left(\frac{|B(x, b^{-j})|^{\gamma}|g(x)|}{(1+b^j\rho(x, y))^{d/r}}\right)^p
\ONE_{A_\xi}(y) d\mu(y)\\
& \qquad\qquad\qquad \le c \int_M
\left(\sup_{x\in M}\frac{|B(x, b^{-j})|^{\gamma}|g(x)|}{(1+b^j\rho(x, y))^{d/r}}\right)^p
d\mu(y)\\
& \qquad\qquad\qquad \le c\int_M \Big[\cM_r\big(|B(\cdot, b^{-j})|^{\gamma}g\big)(y)\Big]^p d\mu(y)\\
& \qquad\qquad\qquad \le c\int_M \Big[|B(y, b^{-j})|^{\gamma}|g(y)|\Big]^p d\mu(y),
\end{align*}
which confirms (\ref{half_shannon}).
Here for the first inequality we used that
$A_\xi \subset B(\xi, cb^{-j})$ and $|B(\xi, b^{-j})| \sim |B(x, b^{-j})|$ if $x\in A_\xi$,
for the third we used Lemma~\ref{lem:Peetre},
and for the last inequality we used the boundedness of the maximal operator
$\cM_r$ on $\LL^p$ when $r<p$.
$\qed$

\medskip

\noindent
{\bf Proof of Theorem~\ref{thm:B-character}.}
We shall only carry out the proof for the spaces $\tB_{pq}^s$.
Also, we only consider the case when $p, q<\infty$. The other cases are similar.

We first prove the boundedness of the synthesis operator
$T_{\psi}: \tb_{pq}^s \rightarrow \tB_{pq}^s$.
To~this end we shall first prove it for finitely supported sequences and
then extend it to the general case.
Let
$a=\{a_\xi\}_{\xi\in\cX}$ be a finitely supported sequence
and set
$f:= T_{\psi}a=\sum_{\xi\in\cX} a_\xi\psi_\xi$.
We shall use the norm on $\tB_{pq}^s$ defined in (\ref{def-Besov-norm4})
(see Proposition~\ref{prop:Bspq-independ}).
We have
$$
\Phi_j(\sqrt{L})f
= \sum_{m=j-1}^{j+1}\sum_{\xi\in\cX_m} a_\xi\Phi_j(\sqrt{L}) \psi_\xi
\quad \hbox{with $\cX_{-1}:=\emptyset$.}
$$

Choose $r$ and $\sigma$ so that $0<r<p$ and $\sigma \ge |s|+d/r+5d/2+1$.
By Theorem~\ref{thm:main-local-kernels} we have the following bound on the kernel
$\Phi_j(\sqrt{L})(x, y)$ of the operator $\Phi_j(\sqrt{L})$:
\begin{align*}
|\Phi_j(\sqrt{L})(x, y)|
&\le cD_{b^{-j}, \sigma}(x, y) \le cD_{b^{-m}, \sigma}(x, y)\\
&\le c|B(y, b^{-m})|^{-1}(1+b^m\rho(x, y))^{-\sigma+d/2},
\quad j-1\le m \le j+1.
\end{align*}
On the other hand, by (\ref{prop-psi-1}) it follows that
$$
|\psi_{\xi} (x)|
\le c |B(\xi, b^{-m})|^{-1/2} \big(1+b^m\rho(x, \xi)\big)^{-\sigma},
\quad \xi\in\cX_m.
$$
Therefore, for $\xi\in\cX_m$, $j-1\le m \le j+1$,
\begin{align*}
&|\Phi_j(\sqrt{L}) \psi_\xi(x)|
= \Big|\int_M |\Phi_j(\sqrt{L})(x, y) \psi_\xi(y) d\mu(y)\Big|\\
& \qquad \le \frac{c}{|B(\xi, b^{-m})|^{1/2}}
\int_M \frac{ d\mu(y)}{|B(y, b^{-m})|(1+b^m\rho(x, y))^{\sigma-d/2}(1+b^m\rho(y, \xi))^\sigma}\\
& \qquad \le \frac{c}{|B(\xi, b^{-m})|^{1/2}(1+b^m\rho(x, \xi))^{\sigma-d/2}}.
\end{align*}
Here for the last inequality we used (\ref{tech3}) and that $\sigma >5d/2$.
From the above we infer
\begin{align*}
|B(x, b^{-j})|^{-s/d}|\Phi_j(\sqrt{L})f(x)|
&\le c\sum_{m=j-1}^{j+1}\sum_{\xi\in\cX_m}
\frac{|a_\xi||B(x, b^{-j})|^{-s/d}}{|B(\xi, b^{-m})|^{1/2}(1+b^m\rho(x, \xi))^{\sigma-d/2}}\\
&\le c\sum_{m=j-1}^{j+1}\sum_{\xi\in\cX_m}
\frac{|a_\xi||B(\xi, b^{-j})|^{-s/d-1/2}}{(1+b^m\rho(x, \xi))^{\sigma-|s|-d/2}},
\quad x\in M.
\end{align*}
Let $\cX_{m, x}:=\{\eta\in \cX_m: x\in A_\eta\}$
and set $Q_\xi:= |a_\xi||B(\xi, b^{-m})|^{-s/d-1/2}$.
Then the~above yields
\begin{align*}
|B(x, b^{-j})|^{-s/d}|\Phi_j(\sqrt{L})f(x)|
&\le c\sum_{m=j-1}^{j+1}\sum_{\eta\in\cX_{m, x}}\sum_{\xi\in\cX_m}
\frac{|a_\xi||B(\xi, b^{-j})|^{-s/d-1/2}}{(1+b^m\rho(\eta, \xi))^{d/r+1}}\\
&= c\sum_{m=j-1}^{j+1}\sum_{\eta\in\cX_{m, x}} Q_\eta^*\ONE_{A_\eta}(x)
= c\sum_{m=j-1}^{j+1}\sum_{\eta\in\cX_m} Q_\eta^*\ONE_{A_\eta}(x)\\
&\le c\sum_{m=j-1}^{j+1} \cM_r\Big(\sum_{\eta\in\cX_m} Q_\eta\ONE_{A_\eta}\Big)(x),
\end{align*}
where we used that $\sigma \ge |s|+d/r+d/2+1$ and for the last inequality
we applied Lemma~\ref{lem:disc-max} with $\tau=d/r+1$.
Therefore,
\begin{align*}
\||B(\cdot, b^{-j})|^{-s/d}|\Phi_j(\sqrt{L})f(\cdot)|\|_p
&\le c\sum_{m=j-1}^{j+1} \Big\|\cM_r\Big(\sum_{\eta\in\cX_m} Q_\eta\ONE_{A_\eta}\Big)(\cdot)\Big\|_p\\
&\le c\sum_{m=j-1}^{j+1} \Big\|\sum_{\eta\in\cX_m} Q_\eta\ONE_{A_\eta}(\cdot)\Big\|_p\\
& \le c\sum_{m=j-1}^{j+1}\Big(
\sum_{\eta\in\cX_m} \big[|B(\eta, b^{-m})|^{-s/d+1/p-1/2}|a_\eta|\big]^{p}
\Big)^{1/p}.
\end{align*}
Here for the second inequality we used the maximal inequality (\ref{max-ineq})
and for the last inequality
that $|A_\eta|\sim |B(\eta, b^{-m})|$ for $\eta\in\cX_m$.
We insert the above in the Besov norm from (\ref{def-Besov-norm4}) to obtain
$\|f\|_{\tB_{pq}^s(\Phi)} \le c\|\{a_\xi\}\|_{\tb_{pq}^s}$.
Thus
\begin{equation}\label{finite-embed-B}
\hbox{
$\|T_{\psi}a\|_{\tB_{pq}^s(\Phi)} \le c\|a\|_{\tb_{pq}^s}$
for any finitely supported sequence $a=\{a_\xi\}$.
}
\end{equation}
Let now $a=\{a_\xi\}_{\xi\in\cX}$ be an arbitrary sequence in $\tb_{pq}^s$.
We order arbitrarily the elements of $\{a_\xi\}_{\xi\in\cX}$ in a sequence with indices $1, 2, \dots$
and denote by $\cX^j \subset \cX$ the indices in $\cX$ of the first $j$ elements of the sequence.
Since $\|\{a_\xi\}\|_{\tb_{pq}^s}<\infty$ it readily follows that
$\{a_\xi\}_{\xi\in\cX^j} \to \{a_\xi\}_{\xi\in\cX}$ in $\tb_{pq}^s$ as $j\to\infty$.
This and (\ref{finite-embed-B}) implies that the series
$\sum_{\xi\in\cX} a_\xi\psi_\xi$ converges in the norm of $\tB_{pq}^s$
and by the continuous embedding of $\tB_{pq}^s$ into $\cD'$
(Proposition~\ref{prop:Bspq-embedding})
it converges in $\cD'$ as well.
Therefore,
$T_{\psi}a=\sum_{\xi\in\cX} a_\xi\psi_\xi$
is well defined for $a=\{a_\xi\} \in \tb_{pq}^s$.
The boundedness of the operator
$T_{\psi}: \tb_{pq}^s \rightarrow \tB_{pq}^s$
follows by a simple limiting argument from (\ref{finite-embed-B}).

\smallskip


We now turn to the proof of the boundedness of the operator
$S_{\tilde\psi}:\tB_{pq}^s\rightarrow \tb_{pq}^s$.
Let $f\in \tB_{pq}^s$.
From (\ref{rep-psi-tpsi}) or (\ref{rep-tpsi}) it follows that
$$
\langle f, \tilde\psi_\xi \rangle
= c_\eps |A_\xi|^{1/2} \big[\Gamma_{\lambda_j}f(\xi)
+ S_{\lambda_j}\Gamma_{\lambda_j}f(\xi)\big]
$$
and hence
\begin{align*}
\sum_{\xi\in\cX_j}\Big(|B(\xi, b^{-j})|^{-s/d+1/p-1/2} |\langle f, \tilde\psi_\xi \rangle|\Big)^p
&\le c \sum_{\xi\in\cX_j}|B(\xi, b^{-j})|^{-sp/d}|\Gamma_{\lambda_j}f(\xi)|^p|A_\xi|\\
&+ c\sum_{\xi\in\cX_j}|B(\xi, b^{-j})|^{-sp/d}|S_{\lambda_j}\Gamma_{\lambda_j}f(\xi)|^p|A_\xi|.
\end{align*}
Since $\Gamma_{\lambda_j}f\in \Sigma_{b^{j+2}}$ we can apply Lemma~\ref{lem:half_shannon} to obtain
\begin{equation}\label{est-Gamma}
\sum_{\xi\in\cX_j}|B(\xi, b^{-j})|^{-sp/d}|\Gamma_{\lambda_j}f(\xi)|^p|A_\xi|
\le c\||B(\cdot, b^{-j})|^{-s/d}\Gamma_{\lambda_j}f\|_p^p.
\end{equation}
To estimate 
the second sum above we denote $g_j(x):=|B(x, b^{-j})|^{-s/d}\Gamma_{\lambda_j}f(x)$
and choose $r$ and $\sigma$ so that $0<r<p$ and $\sigma \ge |s|+ d/r+3d/2+1$.
Observe that by Lemma~\ref{lem:instrument} (b) it follows that
$$
|S_{\lambda_j}(x, y)| \le cD_{b^j, \sigma}(x, y)
\le c|B(x, b^{-j})|^{-1}(1+b^j\rho(x, y))^{-\sigma+d/2}
$$
and hence
\begin{align}\label{est-SG}
&|B(\xi, b^{-j})|^{-s/d}|S_{\lambda_j}\Gamma_{\lambda_j}f(\xi)|\notag\\
&\qquad\qquad\qquad
\le c\int_M \frac{|\Gamma_{\lambda_j}f(y)||B(\xi, b^{-j})|^{-s/d-1}}
{(1+b^j\rho(\xi, y))^{\sigma-d/2}} d\mu(y)\notag\\
&\qquad\qquad\qquad
\le c\sup_{y\in M}\frac{|g_j(y)|}{(1+b^j\rho(\xi, y))^{\sigma-|s|-3d/2-1}}
\int_M \frac{|B(\xi, b^{-j})|^{-1}}{(1+b^j\rho(\xi, y))^{d+1}} d\mu(y)\\
& \qquad\qquad\qquad
\le c\sup_{y\in M}\frac{|g_j(y)|}{(1+b^j\rho(z, y))^{d/r}}
\le c\cM_r(g_j)(z), \quad z\in A_\xi.\notag
\end{align}
Here for the second inequality we used that
$$
|B(\xi, b^{-j})|^{-s/d} \le c (1+b^j\rho(\xi, y))^{|s|} |B(y, b^{-j})|^{-s/d},
$$
which follows by (\ref{D2}),
for the third inequality we used
$\sigma \ge |s|+d/r+3d/2+1$ and (\ref{tech1}),
and for the last inequality we applied Lemma~\ref{lem:Peetre}.
Thus, applying the maximal inequality
\begin{align*}
&\sum_{\xi\in\cX_j}|B(\xi, b^{-j})|^{-sp/d}|S_{\lambda_j}\Gamma_{\lambda_j}f(\xi)|^p|A_\xi|
\le c \sum_{\xi\in\cX_j}\int_{A_\xi} [\cM_r(g_j)(z)]^pd\mu(z)\\
&\qquad\qquad =c \int_M \big[\cM_r(|B(\cdot, b^{-j})|^{-s/d}\Gamma_{\lambda_j}f)(z)\big]^pd\mu(z)
\le c\||B(\cdot, b^{-j})|^{-s/d}\Gamma_{\lambda_j}f\|_p^p.
\end{align*}
From this and (\ref{est-Gamma}) we infer
$$
\sum_{\xi\in\cX_j}\Big(|B(\xi, b^{-j})|^{-s/d + 1/p - 1/2} |\langle f, \tilde\psi_\xi \rangle|\Big)^p
\le c\||B(\cdot, b^{-j})|^{-s/d}\Gamma_{\lambda_j}f\|_p^p,
\quad j\ge 0.
$$
Inserting this in the $\tb^s_{pq}$-norm we get
$
\|\{\langle f, \tilde\psi_\xi \rangle\}\|_{\tb^s_{pq}}
\le c \|f\|_{\tB^s_{pq}(\Gamma)} \le c \|f\|_{\tB^s_{pq}},
$
where we used that the functions $\Gamma_j$, $j\ge 0$, can be used to define
an equivalent norm in $B^s_{pq}$ (see Proposition~\ref{prop:Bspq-independ}).
Thus the boundedness of the operator $S_\psi$ is established.

The equality $T_{\psi}\circ S_{\tilde\psi}=Id$ on $\tB_{pq}^s$ follows by
Proposition~\ref{prop:decomp-DD2} (c).

Assuming the reverse doubling condition (\ref{reverse-doubling}), we have
$\|\psi_\xi\|_p \sim |B(\xi, b^{-j})|^{1/p-1/2}$ from (\ref{prop-psi-2}),
which leads to (\ref{tBnorm-equiv-1}).
$\qed$

\subsection{Embedding of Besov spaces}\label{sec:B-embedding}

Here we show that the Besov spaces $B^s_{pq}$ and $\tB^s_{pq}$ embed ``correctly".


\begin{proposition}\label{B-embedding}
Let $0<p\le p_1<\infty$, $0<q\le q_1\le \infty$, $-\infty<s_1\le s<\infty$.
Then we have the continuous embeddings
\begin{equation}\label{B-embed}
B_{pq}^{s} \subset B_{p_1q_1}^{s_1}
\quad\hbox{and}\quad
\tB_{pq}^{s} \subset \tB_{p_1q_1}^{s_1}
 \quad\mbox{if}\quad
s/d-1/p=s_1/d-1/p_1.
\end{equation}
Here for the left-hand side embedding we assume in addition the non-collapsing condition $(\ref{non-collapsing})$.
\end{proposition}

\noindent
{\bf Proof.} This assertion follows easily by Proposition \ref{prop:Nikolski}.
Let $\{\varphi_j\}_{j\ge 0}$ be the functions from the definition of Besov spaces
(Definition~\ref{def-B-spaces}).
Given $f\in \tB_{p_1q_1}^{s_1}$
we evidently have $\varphi_j(\sqrt{L}) f\in \Sigma_{2^{j+1}}$ and
using  (\ref{norm-relation2})
\begin{align*}
\| |B(\cdot, 2^{-j})|^{-s_1/d} \varphi_j(\sqrt{L}) f(\cdot)\|_{p_1}
&\le c\| |B(\cdot, 2^{-j-1})|^{-s_1/d} \varphi_j(\sqrt{L}) f(\cdot)\|_{p_1}\\
&\le c\| |B(\cdot, 2^{-j-1})|^{-s_1/d+1/p_1-1/p} \varphi_j(\sqrt{L}) f(\cdot)\|_{p}\\
&\le c\| |B(\cdot, 2^{-j})|^{-s/d} \varphi_j(\sqrt{L}) f(\cdot)\|_{p},
\end{align*}
which readily implies
$\|f\|_{\tB_{p_1q_1}^{s_1}} \le c\|f\|_{\tB_{pq}^{s}}$
and hence the right-hand embedding in (\ref{B-embed}) holds.
The left-hand side imbedding in (\ref{B-embed}) follows in the same manner using (\ref{norm-relation1}).
$\qed$

\subsection{Characterization of Besov spaces via linear approximation from \boldmath $\Sigma_t^p$}
\label{sec:char-Besov}

It is natural and easy to characterize the Besov spaces $B_{pq}^s$ with $s>0$ and $p\ge 1$
by means of linear approximation from $\Sigma_t^p$, $t\ge 1$.
In fact, in this case the Besov space $B_{pq}^s$ is the same as the respective approximation space
$A_{pq}^s$ associated with linear approximation from $\Sigma_t^p$.
We refer the reader to \cite{CKP}, \S 3.5 and \S 6.1, for a detailed account of this relationship and
more.

\subsection{Application of Besov spaces to nonlinear approximation}\label{Nonlin-app}

Our aim here is to deploy the Besov spaces to nonlinear approximation.
We shall consider nonlinear $n$-term approximation for the frame
$\{\psi_\eta\}_{\eta\in \cX}$ defined in \S\ref{natural-frame}
with dual frame $\{\tilde\psi_\eta\}_{\eta\in \cX}$ from \S\ref{dual-frame}
or the tight frame $\{\psi_\eta\}_{\eta\in \cX}$ from \S\ref{sec:polyn-prop}.

In this part, we make the additional assumption that
the {\em reverse doubling condition} (\ref{reverse-doubling}) is valid,
and hence (\ref{D3}) holds.

Denote by $\Omega_n$ the nonlinear set of all functions $g$ of the form
$$
g = \sum_{\xi \in \Lambda_n} a_\xi \psi_\xi,
$$
where $\Lambda_n \subset \cX$, $\#\Lambda_n \le n$,
and $\Lambda_n$ may vary with $g$.
We let $\sigma_n(f)_p$ denote the error of best $\LL^p$-approximation to
$f \in \LL^p(M, d\mu)$ from $\Omega_n$, i.e.
$$
 \sigma_n(f)_p := \inf_{g \in \Omega_n} \|f - g\|_p.
$$
The approximation will take place in $\LL^p$, $1\le p<\infty$.
Suppose $s>0$ and let $1/\tau := s/d+1/p$.
The Besov space
$$
\tB_\tau^{s}:=\tB^s_{\tau\tau}
$$
will play a prominent role.

We shall utilize the representation of functions in $\LL^p$ via $\{\psi_\eta\}_{\eta\in \cX}$,
given in Theorem~\ref{thm:dual-frame} \& Proposition~\ref{prop:decomp-DD2}:
For any $f\in\LL^p$, $1\le p < \infty$,
\begin{equation}\label{repres-f}
f=\sum_{\xi\in\cX} \langle f, \tilde\psi_\xi\rangle \psi_\xi
\quad \mbox{in } \LL^p.
\end{equation}
It is readily seen that Theorem~\ref{thm:B-character} and
$\|\psi_\xi\|_p\sim |B(\xi, b^{-j})|^{1/p-1/2}$
for $\xi\in \cX_j$, $j \ge 0$, $0<p\le \infty$, (see (\ref{prop-psi-2}))
imply the following representation of the norm in $\tB_\tau^s$:
\begin{equation}\label{Btau-norm}
\|f\|_{\tB_\tau^{s}}\sim
\Big(\sum_{\xi\in\cX} \|\langle f, \tilde\psi_\xi \rangle \psi_\xi\|_p^\tau\Big)^{1/\tau}
=: \cN(f).
\end{equation}

The next embedding result shows the importance of the spaces $\tB_\tau^{s}$ for
nonlinear $n$-term approximation from $\{\psi_\eta\}_{\eta\in \cX}$ .


\begin{proposition}\label{prop:embed-Lp}
If $f \in \tB_\tau^{s}$, then $f$ can be identified as a function
$f\in \LL^p$ and
\begin{equation}\label{embedding}
\|f\|_p \le
\Big\|\sum_{\xi\in\cX}|\langle f, \tilde\psi_\xi\rangle\psi_\xi(\cdot)|\Big\|_p
\le c\|f\|_{\tB_\tau^{s}}.
\end{equation}
\end{proposition}
We can now give the main result in this section (Jackson estimate):


\begin{theorem}\label{thm:jackson}
If $f \in \tB_\tau^{s}$, then
\begin{equation}\label{jackson}
\sigma_n(f)_p \le c n^{-s/d}\|f\|_{\tB_\tau^{s}},
\quad n\ge 1.
\end{equation}
\end{theorem}

The proofs of Proposition~\ref{prop:embed-Lp} and Theorem~\ref{thm:jackson}
rely on the following lemma.


\begin{lemma}\label{lem:nonlinear}
Let $g=\sum_{\xi\in\cY_n} a_\xi\psi_\xi$, where $\cY_n\subset \cX$ and $\# \cY_n \le n$.
Suppose $\|a_\xi\psi_\xi\|_p\le K$ for $\xi\in\cY_n$, where $0<p<\infty$.
Then
$
\|g\|_p\le cKn^{1/p}.
$
\end{lemma}

\noindent
{\bf Proof.}
This lemma is trivial when $0<p\le 1$.
Suppose $1<p<\infty$.
As in the definition of $\{\psi_\eta\}_{\eta\in \cX}$ in \S\ref{natural-frame},
assume that $\{A_\xi\}_{\xi\in\cX_j}$ $(j\ge 0)$ is a companion to $\cX_j$ disjoint partition of $M$
such that
$B(\xi, \delta_j/2) \subset  A_\xi \subset B(\xi,\delta_j)$, $\xi \in \cX_j$,
with $\delta_j=\gamma b^{-j-2}$.
Fix $0<t<1$, e.g. $t=1/2$.
By the excellent space localization of $\psi_\xi$, given in (\ref{prop-psi-1}), and (\ref{max-ineq-2})
it follows that
$$
|\psi_\xi(x)|\le c (\cM_t\tONE_{A_\xi})(x),
\quad x\in M, \; \xi\in\cX,
$$
and applying the maximal inequality (\ref{max-ineq}) we obtain
$$
\|g\|_p\le c \big\|\sum_{\xi\in\cY_n} \cM_t(a_\xi \tONE_{A_\xi})\big\|_p
\le c\big\|\sum_{\xi\in\cY_n} |a_\xi|\tONE_{A_\xi}\big\|_p.
$$
On the other hand, from $\|a_\xi\psi_\xi\|_p\le K$ and
$\|\psi_\xi\|_p\sim |B(\xi, b^{-j})|^{\frac 1p - \frac 12}$ 
it follows that
$
|a_\xi|
\le cK|A_\xi|^{\frac 12 - \frac 1p}
$
and hence
\begin{equation}\label{est-norm-F}
\|g\|_p\le c K\big\|\sum_{\xi\in\cY_n} |A_\xi|^{-1/p} \ONE_{A_\xi}\big\|_p.
\end{equation}

For any $\xi\in\cX$ we denote by $\cX_\xi$ the set of all $\eta\in\cX$ such that
$A_\eta\cap A_\xi\ne \emptyset$
and $\ell(\eta) \le \ell(\xi)$,
where $\ell(\eta)$, $\ell(\xi)$ are the levels of $\eta$, $\xi$ in $\cX$
(e.g. $\ell(\xi)=j$ if $\xi\in\cX_j$).

Suppose $\xi\in\cX_j$ and let $\eta\in\cX_\xi\cap \cX_\nu$  for some $\nu\le j$.
Since $A_\eta\cap A_\xi\ne \emptyset$ then
$\rho(\xi, \eta)\le cb^{-\nu}$.
Applying (\ref{D3}) we get
$
|B(\xi, \gamma b^{-\nu-2}/2)| \ge cb^{(j-\nu)\bet}|B(\xi, \gamma b^{-j-2})|
$
and also using (\ref{D2}) we obtain
\begin{align*}
|A_\xi|
&\le |B(\xi, \gamma b^{-j-2})|
\le cb^{-(j-\nu)\bet}|B(\xi, \gamma b^{-\nu-2}/2)|\\
&\le cb^{-(j-\nu)\bet}\big[1+ 2\gamma^{-1}b^{\nu+2} \rho(\xi, \eta)\big]^{d}|B(\eta, \gamma b^{-\nu-2}/2)|
\le cb^{-(j-\nu)\bet}|A_\eta|.
\end{align*}
Hence
$
|A_\xi|/|A_\eta|\le c b^{-(j-\nu)\bet}
$
and therefore
\begin{equation}\label{est-AA}
\sum_{\eta\in\cX_\xi} (|A_\xi|/|A_\eta|)^{1/p} \le c <\infty.
\end{equation}

Let $E:=\cup_{\xi\in\cY_n} A_\xi$ and set
$\omega(x):=\min\{|A_\xi|: \xi\in\cY_n, x\in A_\xi\}$ for $x\in E$.
By~(\ref{est-AA}) it follows that
$$
\sum_{\xi\in\cY_n} |A_\xi|^{-1/p} \ONE_{A_\xi}(x) \le c\omega(x)^{-1/p},
\quad x\in E.
$$
We use this and (\ref{est-norm-F}) to obtain
\begin{align*}
\|g\|_p
&\le cK\|\omega^{-1/p}\|_p
=cK\Big(\int_E\omega^{-1}(x)d\mu(x)\Big)^{1/p}\\
&\le cK \Big(\sum_{\xi\in\cY_n} |A_\xi|^{-1}\int_{M}\ONE_{A_\xi}(x)d\mu(x)\Big)^{1/p}
= cK(\#\cY_n)^{1/p}\le cK n^{1/p},
\end{align*}
which completes the proof.
$\qed$

\smallskip


\noindent
{\bf Proof of Proposition~\ref{prop:embed-Lp} \& Theorem~\ref{thm:jackson}.}
The argument is quite standard, but we shall give it for the sake of self-containment. 
Denote briefly
$a_{\xi}:=\langle f, \tilde\psi_\xi\rangle$
and let
$\{a_{\xi_m}\psi_{\xi_m}\}_{m\ge 1}$ be a rearrangement of the sequence $\{a_\xi\psi_\xi\}_{\xi\in \cX}$
such that
$\|a_{\xi_1}\psi_{\xi_1}\|_p \ge \|a_{\xi_2}\psi_{\xi_2}\|_p \ge \cdots$.
Denote
$G_n:= \sum_{m=1}^n a_{\xi_m}\psi_{\xi_m}$.
It suffices to show that
\begin{equation}\label{jackson-1}
\|f-G_n\|\le cn^{-(1/\tau-1/p)}\cN(f) \quad \hbox{for}\quad n\ge 1.
\end{equation}
Assume $\cN(f)>0$ and let
$
\cM_\nu:=\{m: 2^{-\nu}\cN(f)\le\|a_{\xi_m}\psi_{\xi_m}\|_p < 2^{-\nu+1}\cN(f)\}.
$
Denote $K_\ell:=\#\big(\cup_{\nu\le \ell}\cM_\nu\big) $.
Then (\ref{Btau-norm}) yields
$K_\ell \le 2^{\ell\tau}$, $\ell\ge 0$, and hence
$\# \cM_\nu \le 2^{\nu\tau}$, $\nu\ge 0$.
Let $g_\nu:=\sum_{m\in \cM_\nu}a_{\xi_m}\psi_{\xi_m}$.
Now using (\ref{repres-f}) and Lemma~\ref{lem:nonlinear} we infer
\begin{align*}
\|f-G_{K_\ell}\|_p
&\le \big\|\sum_{\nu>\ell} g_\nu\big\|_p
\le \sum_{\nu>\ell} \|g_\nu\|_p
\le c\sum_{\nu>\ell} 2^{-\nu}\cN(f)(\#\cM_\nu)^{1/p}\\
&\le c\cN(f) \sum_{\nu>\ell}2^{-\nu(1-\tau/p)}
\le c\cN(f)2^{-\ell(1-\tau/p)}
\le c\cN(f)2^{-\ell\tau(1/\tau-1/p)}.
\end{align*}
Therefore,
$
\|f-G_{\lfloor 2^{\ell\tau} \rfloor }\|_p \le c\cN(f)2^{-\ell\tau(1/\tau-1/p)},
$
$\forall \ell \ge 0$,

which implies (\ref{jackson-1}).

The proof of Proposition~\ref{prop:embed-Lp} is contained in the above by simply taking
$G_n$ with no terms, i.e. $G_n=0$.
$\qed$

\smallskip

A major {\em open problem} here is to prove the companion to (\ref{jackson}) Bernstein estimate:
\begin{equation}\label{bernstein1}
\|g\|_{\tB_\tau^{s}} \le c n^{s/d} \|g\|_p
\quad \hbox{for}\quad g \in \Omega_n,
\quad 1<p<\infty.
\end{equation}
This estimate would allow to characterize the rates of nonlinear
$n$-term approximation from $\{\psi_\xi\}_{\xi\in \cX}$ in $\LL^p$ ($1<p<\infty$).

\section{Triebel-Lizorkin spaces}\label{F-spaces}
\setcounter{equation}{0}

To introduce Triebel-Lizorkin spaces we shall use the cutoff functions
$\varphi_0, \varphi\in C^\infty_0(\bR_+)$ 
from the definition of Besov spaces (Definition~\ref{def-B-spaces}).
As there we set $\varphi_j(\lambda):= \varphi(2^{-j}\lambda)$ for $j\ge 1$.

The possibly anisotropic nature of the geometry of $M$ is again the reason for introducing
two types of F-spaces.


\begin{definition}\label{def-F-spaces}
Let $s\in \R$, $0<p<\infty$, and $0<q \le \infty$.

$(a)$
The Triebel-Lizorkin space  $F_{pq}^{s}= F_{pq}^{s}(L)$
is defined as the set of all $f \in \cD'$
such that
\begin{equation}\label{def-F-space1}
\|f\|_{F_{pq}^s} :=
\Big\|\Big(\sum_{j\ge 0} \Big(
2^{js}
|\varphi_j(\sqrt L) f(\cdot)|
\Big)^q\Big)^{1/q}\Big\|_{\LL^p} <\infty.
\end{equation}

$(b)$
The Triebel-Lizorkin space  $\tF_{pq}^{s}= \tF_{pq}^{s}(L)$
is defined as the set of all $f \in \cD'$
such that
\begin{equation}\label{def-F-space2}
\|f\|_{\tF_{pq}^s} :=
\Big\|\Big(\sum_{j\ge 0} \Big(
|B(\cdot, 2^{-j})|^{-s/d}
|\varphi_j(\sqrt L) f(\cdot)|
\Big)^q\Big)^{1/q}\Big\|_{\LL^p} <\infty.
\end{equation}
Above the $\ell^q$-norm is replaced by the $\sup$-norm if $q=\infty$.
\end{definition}

As in the case of  Besov spaces it will be convenient for us to use equivalent definitions
of the $F$-spaces which are based on spectral decompositions that utilize $b^j$ rather than $2^j$,
where $b>1$ is the constant from the definition of the frames in \S\ref{sec:frames}.
Let the functions $\Phi_0, \Phi \in C^\infty$ obey (\ref{con_Phi-0})-(\ref{con_Phi})
and as before set $\Phi_j(\lambda):=\Phi(b^{-j}\lambda)$ for $j\ge 1$.
We define new norms on the $F$-spaces by
\begin{align}
\|f\|_{F_{pq}^s(\Phi)}
&:= \Big\|\Big(\sum_{j\ge 0} \Big(b^{js} |\Phi_j(\sqrt L) f(\cdot)|
\Big)^q\Big)^{1/q}\Big\|_{p}
\quad\hbox{and}\label{def-F-space-Phi}\\
\|f\|_{\tF_{pq}^s(\Phi)}
&:= \Big\|\Big(\sum_{j\ge 0} \Big(|B(\cdot, b^{-j})|^{-s/d}
|\Phi_j(\sqrt{L}) f(\cdot)|\Big)^q\Big)^{1/q}\Big\|_{p}. \label{def-tF-space-Phi}
\end{align}


\begin{proposition}\label{prop:Fspq-independ}
For all admissible indices $\|\cdot\|_{F_{pq}^s}$ and $\|\cdot\|_{F_{pq}^s(\Phi)}$
are equivalent quasi-norms on $\tF_{pq}^s$,
and $\|\cdot\|_{\tF_{pq}^s}$ and $\|\cdot\|_{\tF_{pq}^s(\Phi)}$
are equivalent quasi-norms on $\tF_{pq}^s$.
Therefore, the~definitions of $F_{pq}^s$ and $\tF_{pq}^s$ are independent of
the particular selection of $\varphi_0, \varphi$.
\end{proposition}

\noindent
{\bf Proof.}
We shall only establish the equivalence of
$\|\cdot\|_{\tF_{pq}^s(\Phi)}$ and $\|\cdot\|_{\tF_{pq}^s}$.
As in the proof of Proposition~\ref{prop:Bspq-independ} there exist functions
$\tilde\Phi_0$ and $\tilde\Phi$
with the properties of $\Phi_0$ and $\Phi$ from (\ref{con_Phi-0})-(\ref{con_Phi})
such that
$$
\tilde\Phi_0(\lambda)\Phi_0(\lambda) + \sum_{j\ge 0}\tilde\Phi(b^{-j}\lambda)\Phi(b^{-j}\lambda)=1,
\quad \lambda\in \bR_+.
$$
Setting $\tilde\Phi_j(\lambda):=\tilde\Phi(b^{-j}\lambda)$ for $j\ge 1$
we have
$
\sum_{j\ge 0}\tilde\Phi_j(\lambda)\Phi_j(\lambda)=1,
$
which implies
$f=\sum_{j\ge 0}\tilde\Phi_j(\sqrt{L})\Phi_j(\sqrt{L})f$
in $\cD'$.

Assume $1<b<2$ (the proof in the case $b\ge 2$ is similar) and let $j\ge 1$.
Clearly, there exist $\ell >1$ 
and $m\ge 1$ such that
$[2^{j-1}, 2^{j+1}] \subset [b^{m-1}, b^{m+\ell+1}]$.
Now, precisely as in the proof of Proposition~\ref{prop:Bspq-independ} we have
$$
\varphi_j(\sqrt L)f(x)
= \sum_{\nu=m}^{m+\ell} \varphi_j(\sqrt L)\tilde\Phi_\nu(\sqrt{L})\Phi_\nu(\sqrt{L})f(x)
$$
and for $m\le \nu\le m+\ell$
$$
|\varphi_j(\sqrt L)\tilde\Phi_\nu(\sqrt{L})\Phi_\nu(\sqrt{L})f(x)|
\le \frac{c}{|B(x, b^{-\nu})|}
\int_M\frac{|\Phi_\nu(\sqrt{L})f(y)|}{(1+b^\nu\rho(x, y))^{\sigma-d/2}} d\mu(y).
$$
Let $0<r<\min\{p, q\}$ and choose $\sigma\ge |s|+ d/r+3d/2+1$.
Then just as in the proof of Proposition~\ref{prop:Bspq-independ} we obtain
$$
|B(x, 2^{-j})|^{-s/d}|\varphi_j(\sqrt{L})f(x)|
\le c\sum_{\nu=m}^{m+\ell}
\cM_r\Big(|B(\cdot, b^{-\nu})|^{-s/d}\Phi_\nu(\sqrt{L})f\Big)(x).
$$

A~similar estimate holds for $j=0$.
We use the above in the definition of $\|f\|_{\tF_{pq}^s}$
and the maximal inequality (\ref{max-ineq}) to obtain
\begin{align*}
\|f\|_{\tF_{pq}^s}
&\le  c\Big\|\Big(\sum_{\nu\ge 0}
\Big[\cM_r\big(|B(\cdot, b^{-\nu})|^{-s/d}\Phi_\nu(\sqrt{L})f\big)(\cdot)\Big]^q\Big)^{1/q}\Big\|_p\\
&\le  c\Big\|\Big(\sum_{\nu\ge 0}
\Big[|B(\cdot, b^{-\nu})|^{-s/d}\Phi_\nu(\sqrt{L})f\Big]^q\Big)^{1/q}\Big\|_p
= c\|f\|_{\tF_{pq}^s(\Phi)}.
\end{align*}
In the same way one proves the estimate
$\|f\|_{\tF_{pq}^s(\Phi)} \le c\|f\|_{\tF_{pq}^s}$.
$\qed$


\begin{proposition}\label{prop:Fspq-embedding}
The F-spaces $F_{pq}^s$ and $\tF_{pq}^s$ are quasi-Banach spaces which
are continuously embedded in $\cD'$.
More precisely,
for all admissible indices $s, p, q$,
we have:

$(a)$
If $\mu(M)<\infty$, then
\begin{equation}\label{embed-Fspq-1}
|\langle f, \phi \rangle|
\le c \|f\|_{F_{pq}^s}\cP_{m}^*(\phi),
\quad
f\in F_{pq}^s, \;\; \phi\in\cD,
\end{equation}
when
$2m > d\big(\frac{1}{\min\{p, 1\}}-1\big)-s$,
and
\begin{equation}\label{embed-tFspq-1}
|\langle f, \phi \rangle|
\le c \|f\|_{\tF_{pq}^s}\cP_{m}^*(\phi),
\quad
f\in \tF_{pq}^s, \;\; \phi\in\cD,
\end{equation}
when
$2m > \max\big\{0, d\big(\frac{1}{\min\{p, 1\}}-1\big)-s\big\}$.

$(b)$
If $\mu(M)=\infty$, then
\begin{equation}\label{embed-Fspq-2}
|\langle f, \phi \rangle|
\le c \|f\|_{F_{pq}^s}\cP_{m,\ell}^*(\phi),
\quad
f\in F_{pq}^s, \;\; \phi\in\cD,
\end{equation}
when
$2m > d\big(\frac{1}{\min\{p, 1\}}-1\big)-s$
and $\ell>2d$, and
\begin{equation}\label{embed-tFspq-2}
|\langle f, \phi \rangle|
\le c \|f\|_{\tF_{pq}^s}\cP_{m,\ell}^*(\phi),
\quad
f\in \tF_{pq}^s, \;\; \phi\in\cD,
\end{equation}
when
$2m > \max\big\{0, d\big(\frac{1}{\min\{p, 1\}}-1\big)-s\big\}$
and $\ell>\max\big\{2d, \big|d\big(\frac{1}{p}-1\big)\big|+|s|\big\}$.
\end{proposition}

\noindent
{\bf Proof.} The proof of this proposition is essentially the same as the proof
Proposition~\ref{prop:Bspq-embedding}.
One only has to observe that
$
 \||B(x, 2^{-j})|^{-s/d} \varphi_j(\sqrt L)f\|_p \le \|f\|_{\tF_{pq}^s}
$
and replace $\|f\|_{\tB_{pq}^s}$ by $\|f\|_{\tF_{pq}^s}$ everywhere in the proof of
Proposition~\ref{prop:Bspq-embedding}.
$\qed$

\subsection{Heat kernel characterization of Triebel-Lizorkin spaces}

Our aim is to show that the spaces $F_{pq}^{s}$ and $\tF_{pq}^{s}$ can be equivalently
defined using directly the~heat kernel when $p, q$ are restricted to
$1 < p < \infty$ and $1 < q \le \infty$.


\begin{definition}\label{def:FH-norms}
Given $s \in \R$, let $m$ be the smallest $m\in \bZ_+$ such that $m>s$.
We define
\begin{align*}
&\|f\|_{F_{pq}^{s}(H)}
:= \|e^{-L}f\|_p
+ \Big\|\Big(\int_0^1 \big[t^{-s/2}  |(tL)^{m/2}
e^{-tL}f(\cdot) | \big]^q \frac{dt}{t}\Big)^{1/q}\Big\|_p \quad\hbox{and}\\
&\|f\|_{\tF_{pq}^{s}(H)}
:= \||B(\cdot, 1)|^{-\frac{s}{d}}e^{-L}f\|_p
+ \Big\|\Big(\int_0^1 \Big[|B(\cdot, t^{1/2})|^{-\frac{s}{d}}
|(tL)^{\frac{m}{2}}e^{-tL}f (\cdot)|\Big]^q \frac{dt}{t}\Big)^{\frac{1}{q}}
\Big\|_p
\end{align*}
with the usual modification when $q=\infty$.
\end{definition}


\begin{theorem}\label{thm:Heat-F-charac}
Suppose $s \in \R$, $1 < p <\infty$, $1< q \le \infty$,
and $m>s$, $m\in \bZ_+$ as in the above definition.

$(a)$
If $f\in \cD'$, then
$f\in F_{pq}^{s}$ if and only if $e^{-L}f\in \LL^p$ and
$\|f\|_{F_{pq}^{s}(H)} < \infty$.
Moreover, if $f\in F_{pq}^{s}$, then
$\|f\|_{F_{pq}^{s}}\sim \|f\|_{F_{pq}^{s}(H)}$.

$(b)$
If $f\in \cD'$, then
$f\in \tF_{pq}^{s} \Longleftrightarrow |B(\cdot, 1)|^{s/d}e^{-L}f\in \LL^p$ and
$\|f\|_{\tF_{pq}^{s}(H)} < \infty$.
Moreover, if $f\in \tF_{pq}^{s}$, then
$\|f\|_{\tF_{pq}^{s}}\sim \|f\|_{\tF_{pq}^{s}(H)}$.
\end{theorem}

\noindent
{\bf Proof.}
We shall only prove Part (b). The proof of Part (a) is similar and will be omitted.
The proof bears a lot of similarities with the proof of Theorem~\ref{thm:Heat-B-charac}
and we shall utilize some parts from the latter.

Let $\varphi_0$, $\varphi$, and $\varphi_j$, $j\ge 1$, be Littlewood-Paley functions,
just as in the proof of Theorem~\ref{thm:Heat-B-charac}.
Then  $f=\sum_{j\ge 0} \varphi_j^2(\sqrt L)f$ for $f\in \cD'$ and hence
$$
|B(\cdot, t^{1/2})|^{-s/d}(tL)^{m/2} e^{-tL}f
= \sum_{j\ge 0} |B(\cdot, t^{1/2})|^{-s/d}(tL)^{m/2} e^{-tL}\varphi_j^2(\sqrt L)f
=: \sum_{j\ge 0}F_j.
$$
Now, precisely as in the proof of Theorem~\ref{thm:Heat-B-charac} (see (\ref{est-Fj-B})) we obtain
$$
|F_{j}(x)| \le c\big(1+ (t4^j)^{-\frac{s}{2}}\big)(t4^{j})^{\frac{m}{2}}e^{-t4^{j}}
\int_M \frac{|B(y, 2^{-j})|^{-\frac{s}{d}}|\varphi(2^{-j}\sqrt L)f(y)|}
{|B(x, 2^{-j})|(1+2^j\rho(x, y))^{\sigma-|s|}} d\mu(y).
$$
Choose $r$ and $\sigma$ so that $0< r < \min\{p, q\}$ and $\sigma \ge |s|+\frac{d}{r}+d+1$,
and denote briefly
$h_{j}(t):= \big[(t4^{j})^{\frac{m}{2}}+ (t4^j)^{(m-s)/2}\big]e^{-t4^{j}}$.
Evidently, $\varphi(2^{-j}\sqrt L)f \in \Sigma_{2^{j+1}}$ and applying Lemma~\ref{lem:Peetre}
we get for $j\ge 1$
\begin{align*}
|F_{j}(x)|
&\le \frac{ch_{j}(t)}{|B(x, 2^{-j})|}\int_M \frac{1}{(1+2^j\rho(x, y))^{d+1}} d\mu(y)\\
&  \qquad\qquad\qquad\qquad\qquad\qquad
\times
\sup_{y\in M} \frac{|B(y, 2^{-j})|^{-\frac{s}{d}}|\varphi(2^{-j}\sqrt L)f(y)|}
{(1+2^j\rho(x, y))^{d/r}}\\
& \le ch_{j}(t)\cM_r\Big(|B(\cdot, 2^{-j})|^{-\frac{s}{d}}\varphi(2^{-j}\sqrt L)f\Big)(x).
\end{align*}
Here we used (\ref{tech1}) in estimating the integral,
and $\cM_r$ is the maximal operator, defined in (\ref{def:max-op}).
Hence
$$
|F_{j}(x)| \le ch_{j}(t)\cM_r\Big(|B(\cdot, 2^{-j})|^{-\frac{s}{d}}\varphi_j(\sqrt L)f\Big)(x),
\quad j\ge 1.
$$
Similarly as above we obtain
$$
|F_{0}(x)| \le ch_1(t)\cM_r\Big(|B(\cdot,1)|^{-\frac{s}{d}}
\varphi_0(\sqrt L)f\Big)(x).
$$
Set
$b_j(x):= \cM_r\big(|B(\cdot, 2^{-j})|^{-\frac{s}{d}}\varphi_j(\sqrt L)f\big)(x)$.
Let $1<q<\infty$.
From the above estimates we get
\begin{align*}
&\Big\|\Big(\int_0^1 \big[|B(\cdot, t^{\frac 12})|^{-\frac{s}{d}}(tL)^{\frac{m}{2}} e^{-tL}f\big]^q
\frac{dt}{t}\Big)^{1/q}\Big\|_p
\le c\Big\|\Big(\int_0^1 \Big[\sum_{j\ge 0}h_j(t)b_j(\cdot)\Big]^q
\frac{dt}{t}\Big)^{1/q}\Big\|_p \\
& \le c\Big\|\Big(\sum_{\nu\ge 0}\int_{4^{-\nu-1}}^{4^{-\nu}} \Big[\sum_{j\ge 0}h_j(t)b_j(\cdot)\Big]^q
\frac{dt}{t}\Big)^{1/q}\Big\|_p
\le c\Big\|\Big(\sum_{\nu\ge 0}\Big[\sum_{j\ge 0}a_{j-\nu}b_j(\cdot)\Big]^q
\Big)^{1/q}\Big\|_p.
\end{align*}
Here
$$
a_{j-\nu}:= \max \{h_j(t): t\in[4^{-\nu-1}, 4^{-\nu}]\}
\le c\big(4^{(j-\nu)m/2} + 4^{(j-\nu)(m-s)/2}\big)e^{-4^{j-\nu-1}}
$$
and we set
$a_\nu := (4^{\nu m/2}+ 4^{\nu(m-s)/2})e^{-4^{\nu-1}}$, $\nu\in \bZ$.
We apply Young's inequality to the convolution of the above sequences to obtain
$$
\Big(\sum_{\nu\ge 0}\Big(\sum_{j\ge 0}a_{j-\nu} b_j(x)\Big)^q\Big)^{1/q}
\le \sum_{\nu\in\bZ} a_\nu \Big(\sum_{j\ge 0} b_j(x)^q\Big)^{1/q}
\le c\Big(\sum_{j\ge 0} b_j(x)^q\Big)^{1/q},
$$
where we used that
$\sum_{\nu\in\bZ} a_\nu \le c$ due to $m>s$.
Therefore,
\begin{align*}
&\Big\|\Big(\int_0^1 \big[|B(\cdot, t^{\frac 12})|^{-\frac{s}{d}}(tL)^{\frac{m}{2}} e^{-tL}f\big]^q
\frac{dt}{t}\Big)^{1/q}\Big\|_p
\le c\Big\|
\Big(\sum_{j\ge 0} b_j(\cdot)^q\Big)^{1/q}
\Big\|_p \\
& \qquad\qquad\qquad\qquad
= c\Big\|
\Big(\sum_{j\ge 0}\Big[
\cM_r\Big(|B(\cdot, 2^{-j})|^{-\frac{s}{d}}\varphi_j(\sqrt L)f\Big)
\Big]^q\Big)^{1/q}
\Big\|_p\\
& \qquad\qquad\qquad\qquad
\le c\Big\|
\Big(\sum_{j\ge 0}
\Big(|B(\cdot, 2^{-j})|^{-\frac{s}{d}}\varphi_j(\sqrt L)f(\cdot)\Big)^q\Big)^{1/q}
\Big\|_p
\le c\|f\|_{\tF_{pq}^{s}}.
\end{align*}
Here in the former inequality we used the maximal inequality (\ref{max-ineq}).

It is easier to show that
$
\||B(\cdot, 1)|^{-s/d}e^{-L}f\|_p \le c\|f\|_{\tF_{pq}^{s}}.
$
The proof follows in the footsteps of the above proof and will be omitted.
Combining the above two estimates we get
$\|f\|_{\tF_{pq}^{s}(H)} \le c\|f\|_{\tF_{pq}^{s}}$.
The derivation of this estimate in the case $q=\infty$ is easier and will be omitted.


We next prove an estimate in the opposite direction.
We only consider the case when $1<q<\infty$; the case $q=\infty$ is easier.
Assume that $\varphi_0$, $\varphi$, and $\varphi_j$, $j\ge 1$, are as in the definition
of $F^s_{pq}$ (Definition~\ref{def-F-spaces}).
For $j\ge 1$, we obtain exactly as in the proof of Theorem~\ref{thm:Heat-B-charac} (see (\ref{est-Bphij}))
$$
|B(x, 2^{-j})|^{-s/d}|\varphi_{j}(\sqrt L)f(x)|
\le c\int_M\frac{|B(y, t^{1/2})|^{-s/d}|(tL)^{m/2} e^{-tL}f(y)|}
{|B(x, 2^{-j})|\big(1+2^j\rho(x, y)\big)^{\sigma-|s|}} d\mu(y).
$$
Choose $\sigma \ge |s|+d+1$ and denote briefly
$
F(x, t) := |B(x, t^{1/2})|^{-s/d}|(tL)^{m/2} e^{-tL}f(x)|.
$
Set $S_m:=\{y\in M:2^{m-1}\le 2^j\rho(x, y)<2^m\}$, $S_m\subset B(x, 2^{m-j})$.
Then
\begin{align*}
&|B(x, 2^{-j})|^{-s/d}|\varphi_{j}(\sqrt L)f(x)|
\le c\int_{B(x, 2^{-j})} \cdots + c\sum_{m\ge 1} \int_{S_m} \cdots\\
&\qquad\qquad
\le c\sum_{m\ge 0}\frac{|B(x, 2^{m-j})|}{|B(x, 2^{-j})|2^{m(d+1)}}
\frac{1}{|B(x, 2^{m-j})|}\int_{B(x, 2^{m-j})}F(y, t) d\mu(y)\\
&\qquad\qquad
\le c(\cM_1F(\cdot, t))(x)\sum_{m\ge 0}2^{-m} \le c(\cM_1F(\cdot, t))(x),
\end{align*}
where we used (\ref{doubling}).
Therefore, for any $t\in [4^{-j}, 4^{-j+1}]$ and $x\in M$
$$
|B(x, 2^{-j})|^{-s/d}|\varphi_{j}(\sqrt L)f(x)|
\le c\cM_1(F(\cdot, t))(x),
$$
which yields
$$
|B(x, 2^{-j})|^{-s/d}|\varphi_{j}(\sqrt L)f(x)|^q
\le c\int_{4^{-j}}^{4^{-j+1}} \Big[\cM_1(F(\cdot, t))(x)\Big]^q\frac{dt}{t},
\quad x\in M.
$$
These readily imply
\begin{align*}
&\Big\|\Big(\sum_{j\ge 1} \Big(
|B(\cdot, 2^{-j})|^{-s/d}
|\varphi_j(\sqrt L) f(\cdot)|
\Big)^q\Big)^{1/q}\Big\|_{p}\\
& \qquad\qquad
\le c \Big\|\Big(\sum_{j\ge 1}
\int_{4^{-j}}^{4^{-j+1}}\Big[\cM_1(F(\cdot, t))(\cdot)\Big]^q \frac{dt}{t}
\Big)^{1/q}\Big\|_{p}\\
& \qquad\qquad
= c \Big\|\Big(
\int_0^1\Big[\cM_1(F(\cdot, t))(\cdot)\Big]^q \frac{dt}{t}
\Big)^{1/q}\Big\|_{p}
\le c \Big\|\Big(
\int_0^1 |F(\cdot, t)|^q \frac{dt}{t}
\Big)^{1/q}\Big\|_{p}.
\end{align*}
Here for the latter inequality we used the maximal inequality (\ref{max-ineq-integral}).
One easily obtains
$$
\||B(x, 1)|^{-s/d}\varphi_{0}(\sqrt L)f\|_p
\le c\||B(\cdot, 1)|^{-s/d} e^{-L}f\|_p.
$$
The above estimates imply
$\|f\|_{\tF_{pq}^{s}} \le c\|f\|_{\tF_{pq}^{s}(H)}$
and this completes the proof.
$\qed$

\subsection{Frame decomposition of Triebel-Lizorkin spaces}

Here we present the characterization of the F-spaces $F_{pq}^s$ and $\tF_{pq}^s$ via
the frames $\{\psi_\xi\}_{\xi\in\cX}$, $\{\tilde\psi_\xi\}_{\xi\in\cX}$
from \S\S\ref{natural-frame}-\ref{dual-frame}.
We adhere to the notation from \S\ref{sec:frames}, in particular,
$\cX:= \cup_{j\ge 0} \cX_j$ will denote the sets of the centers of the frame elements and
$\{A_\xi\}_{\xi\in\cX_j}$ will be the associated partitions of $M$.

We first introduce the sequence spaces $f_{pq}^s$ and $\tf_{pq}^s$ associated with
$F_{pq}^s$ and $\tF_{pq}^s$, respectively.


\begin{definition}\label{def-f-spaces}
Suppose $s\in \R$, $0<p<\infty$, and $0<q\le\infty$.

$(a)$
$f_{pq}^s$
is defined as the space of all complex-valued sequences
$a:=\{a_{\xi}\}_{\xi\in \cX}$ such that
\begin{equation}\label{def-f-space}
\|a\|_{f_{pq}^s}
:=\Big\|\Big(\sum_{j\ge 0}b^{jsq}\sum_{\xi \in \cX_j}
\big[|a_{\xi}|\tONE_{A_\xi}(\cdot)\big]^q
\Big)^{1/q}
\Big\|_{\LL^p} <\infty.
\end{equation}

$(b)$
$\tf_{pq}^s$
is defined as the space of all complex-valued sequences
$a:=\{a_{\xi}\}_{\xi\in \cX}$ such that
\begin{equation}\label{def-tf-space}
\|a\|_{\tf_{pq}^s}
:=\Big\|\Big(\sum_{\xi \in \cX}
\big[|A_\xi|^{-s/d}|a_{\xi}|\tONE_{A_\xi}(\cdot)\big]^q
\Big)^{1/q}
\Big\|_{\LL^p} <\infty.
\end{equation}
Above the $\ell^q$-norm is replaced by the $\sup$-norm when $q=\infty$.
Recall that $\tONE_{A_\xi}:=|A_\xi|^{-1/2}\ONE_{A_\xi}$
with $\ONE_{A_\xi}$ being the characteristic function of $A_\xi$.
\end{definition}

As in the case of Besov spaces we shall use the ``analysis" and ``synthesis" operators defined by
\begin{equation}\label{anal_synth_oprts}
S_{\tilde\psi}: f \rightarrow \{\langle f, \tilde\psi_\xi\rangle\}_{\xi \in \cX}
\quad\text{and}\quad
T_{\psi}: \{a_\xi\}_{\xi \in \cX} \rightarrow \sum_{\xi\in \cX}a_\xi\psi_\xi.
\end{equation}
Here the roles of $\{\psi_\xi\}$, $\{\tilde\psi_\xi\}$ are interchangeable.


\begin{theorem}\label{thm:F-character}
Let $s\in \R$, $0< p< \infty$ and $0<q\le \infty$.
$(a)$
The operators
$S_{\tilde\psi}:F_{pq}^s \rightarrow f_{pq}^s$ and $T_{\psi}:f_{pq}^s \rightarrow F_{pq}^s$
are bounded and $T_{\tilde\psi}\circ S_\psi=Id$ on $F_{pq}^s$.
Consequently, $f\in F_{pq}^s$ if and only if
$\{\langle f, \tilde\psi_\xi\rangle\}_{\xi \in \cX}\in f_{pq}^s$,
and if $f\in F_{pq}^s$, then
$\|f\|_{F_{pq}^s} \sim  \|\{\langle f, \tilde\psi_\xi\rangle\}\|_{f_{pq}^s}$.
Furthermore,
\begin{equation}\label{F-norm-equivalence-1}
\|f\|_{F_{pq}^s}
\sim \Big\|\Big(
\sum_{j\ge 0}b^{jsq}\sum_{\xi\in \cX_j}
\big[|\langle f, \tilde\psi_\xi \rangle||\psi_\xi(\cdot)|\big]^q
\Big)^{1/q}\Big\|_{\LL^p}.
\end{equation}

\noindent
$(b)$
The operators
$S_{\tilde\psi}:\tF_{pq}^s \rightarrow \tf_{pq}^s$ and $T_{\psi}: \tf_{pq}^s \rightarrow \tF_{pq}^s$
are bounded and $T_{\tilde\psi}\circ S_\psi=Id$ on $\tF_{pq}^s$.
Hence, $f\in \tF_{pq}^s$ if and only if
$\{\langle f, \tilde\psi_\xi\rangle\}_{\xi \in \cX}\in \tf_{pq}^s$,
and if $f\in F_{pq}^s$, then
$\|f\|_{F_{pq}^s} \sim  \|\{\langle f, \tilde\psi_\xi\rangle\}\|_{\tf_{pq}^s}$.
Furthermore, %
\begin{equation}\label{F-norm-equivalence-2}
\|f\|_{\tF_{pq}^s}
\sim \Big\|\Big(
\sum_{\xi\in \cX}
\big[|B(\xi, b^{-j})|^{-s/d}|\langle f, \tilde\psi_\xi \rangle||\psi_\xi(\cdot)|\big]^q
\Big)^{1/q}\Big\|_{\LL^p}.
\end{equation}
Above the roles of $\psi_\xi$ and $\tilde\psi_\xi$ can be interchanged.
\end{theorem}

\noindent
{\bf Proof.}
We shall only prove Part (b).
Also, we shall only consider the case when $q<\infty$. The case $q=\infty$ is similar.

This proof runs parallel to the proof of Theorem~\ref{thm:B-character}
and we shall borrow a lot from that proof.
We begin by proving the boundedness of the synthesis operator
$T_{\psi}: \tf_{pq}^s \rightarrow \tF_{pq}^s$.
To this end we shall first prove it for finitely supported sequences and
then extend it to the general case.
Let
$a=\{a_\xi\}_{\xi\in\cX}$ be a finitely supported sequence
and set
$f:= T_{\psi}a=\sum_{\xi\in\cX} a_\xi\psi_\xi$.
We shall use the norm on $\tF_{pq}^s$ defined in (\ref{def-tF-space-Phi})
(see Proposition~\ref{prop:Fspq-independ}).

Choose $r$ and $\sigma$ so that
$0<r< \min\{p, q\}$ and $\sigma \ge |s|+d/r+3d/2+1$.
Now, precisely as in the proof of Theorem~\ref{thm:B-character} we get
\begin{align*}
|B(x, b^{-j})|^{-s/d}|\Phi_j(\sqrt{L})f(x)|
\le c\sum_{m=j-1}^{j+1} \cM_r\Big(\sum_{\eta\in\cX_m} Q_\eta\ONE_{A_\eta}\Big)(x)
\quad \hbox{with $\cX_{-1}:=\emptyset$,}
\end{align*}
where $Q_\eta:= |a_\eta||B(\eta, b^{-m})|^{-s/d-1/2}$.
Inserting the above in the definition of $\tF_{pq}^s(\Phi)$ from (\ref{def-tF-space-Phi})
we get
\begin{align*}
\|f\|_{\tF_{pq}^s(\Phi)}
&=\Big\|
\Big(\sum_{j\ge 0}\Big[|B(\cdot, b^{-j})|^{-s/d}|\Phi_j(\sqrt{L})f(\cdot)|\Big]^q\Big)^{1/q}
\Big\|_p\\
&\le c \Big\|
\Big(\sum_{m\ge 0}\Big[\cM_r\Big(\sum_{\eta\in\cX_m} Q_\eta\ONE_{A_\eta}\Big)(\cdot)\Big]^q\Big)^{1/q}
\Big\|_p\\
&\le c \Big\|
\Big(\sum_{m\ge 0}\Big[\sum_{\eta\in\cX_m} Q_\eta\ONE_{A_\eta}\Big]^q\Big)^{1/q}
\Big\|_p\\
&=c \Big\|
\Big(\sum_{m\ge 0}\sum_{\eta\in\cX_m} \Big[|a_\eta||B(\eta, b^{-m})|^{-s/d-1/2}
\ONE_{A_\eta}\Big]^q\Big)^{1/q}
\Big\|_p
\le c\|\{a_\eta\}\|_{\tf_{pq}^s}.
\end{align*}
Here for the second inequality we used the maximal inequality (\ref{max-ineq})
and for the last inequality
that $|A_\eta|\sim |B(\eta, b^{-m})|$ for $\eta\in\cX_m$.
Thus
$\|T_{\psi}a\|_{\tF_{pq}^s(\Phi)} \le c\|a\|_{\tf_{pq}^s}$
for any finitely supported sequence $a=\{a_\xi\}$.
Now, just as in the proof of Theorem~\ref{thm:B-character} we conclude that
$T_{\psi}a=\sum_{\xi\in\cX} a_\xi\psi_\xi$
is well defined for $\{a_\xi\}_{\xi\in\cX} \in \tf_{pq}^s$
and the operator
$T_{\psi}: \tf_{pq}^s \rightarrow \tF_{pq}^s$
is bounded.

\smallskip


We now prove the boundedness of the operator
$S_{\tilde\psi}:\tF_{pq}^s\rightarrow \tf_{pq}^s$.
Let $f\in \tF^s_{pq}$ and choose $r$ so that $0< r< \min\{p, q\}$.
By (\ref{rep-tpsi}) it follows that
$$
\langle f, \tilde\psi_\xi \rangle
= c_\eps |A_\xi|^{1/2} \big[\Gamma_{\lambda_j}f(\xi)
+ S_{\lambda_j}\Gamma_{\lambda_j}f(\xi)\big],
$$
which implies
\begin{align*}
\sum_{\xi\in\cX_j}\big[|A_\xi|^{-s/d}|\langle f, \tilde\psi_\xi \rangle|\tONE_{A_\xi}(x)\big]^q
&\le c \sum_{\xi\in\cX_j}|A_\xi|^{-sq/d}|\Gamma_{\lambda_j}f(\xi)|^q|\ONE_{A_\xi}(x)\\
&+ c\sum_{\xi\in\cX_j}|A_\xi|^{-sq/d}|S_{\lambda_j}\Gamma_{\lambda_j}f(\xi)|^q\ONE_{A_\xi}(x).
\end{align*}
Now, we use that $\Gamma_{\lambda_j}f\in \Sigma_{b^{j+2}}$ and Lemma~\ref{lem:Peetre} to obtain for $x\in M$
\begin{align}\label{est-Gamma-F}
\sum_{\xi\in\cX_j}|A_\xi|^{-sq/d}|\Gamma_{\lambda_j}f(\xi)|^q|\ONE_{A_\xi}(x)
&\le c \sum_{\xi\in\cX_j}
\big[\sup_{y\in A_\xi} |B(y, b^{-j})|^{-s/d}|\Gamma_{\lambda_j}f(y)|\big]^q|\ONE_{A_\xi}(x)\notag\\
&\le c\sum_{\xi\in\cX_j}
\Big(\sup_{y\in A_\xi}\frac{|B(y, b^{-j})|^{-s/d}|\Gamma_{\lambda_j}f(y)|}{(1+b^j\rho(x, y))^{d/r}}\Big)^q
\ONE_{A_\xi}(x)\notag\\
&\le c \Big(\sup_{y\in M}\frac{|B(y, b^{-j})|^{-s/d}|\Gamma_{\lambda_j}f(y)|}{(1+b^j\rho(x, y))^{d/r}}\Big)^q\\
&\le c\big[\cM_r\big(|B(\cdot, b^{-j})|^{-s/d}\Gamma_{\lambda_j}f\big)(x)\big]^q. \notag
\end{align}
On the other hand, as in the proof of Theorem~\ref{thm:B-character} (see (\ref{est-SG})) we obtain
\begin{align*}
|A_\xi|^{-s/d}|S_{\lambda_j}\Gamma_{\lambda_j}f(\xi)|
\le c\sup_{y\in M}\frac{|B(y, b^{-j})|^{-s/d}|\Gamma_{\lambda_j}f(y)|}{(1+b^j\rho(x, y))^{d/r}},
\quad x\in A_\xi,\quad \xi\in \cX_j,
\end{align*}
and hence as above using again Lemma~\ref{lem:Peetre}
\begin{align*}
\sum_{\xi\in\cX_j}|A_\xi|^{-sq/d}|S_{\lambda_j}\Gamma_{\lambda_j}f(\xi)|^q\ONE_{A_\xi}(x)
&\le \sum_{\xi\in\cX_j}
\Big(\sup_{y\in M}\frac{|B(y, b^{-j})|^{-s/d}|\Gamma_{\lambda_j}f(y)|}{(1+b^j\rho(x, y))^{d/r}}\Big)^q
\ONE_{A_\xi}(x)\\
&\le c\big[\cM_r\big(|B(\cdot, b^{-j})|^{-s/d}\Gamma_{\lambda_j}f\big)(x)\big]^q.
\end{align*}
This and (\ref{est-Gamma-F}) yield
$$
\sum_{\xi\in\cX_j}\big[|A_\xi|^{-s/d}|\langle f, \tilde\psi_\xi \rangle|\tONE_{A_\xi}(x)\big]^q
\le c\big[\cM_r\big(|B(\cdot, b^{-j})|^{-s/d}\Gamma_{\lambda_j}f\big)(x)\big]^q.
$$
Inserting this in the $\tf^s_{pq}$-norm (Definition~\ref{def-f-spaces})
and using the maximal inequality (\ref{max-ineq}) we get
\begin{align*}
\|\{\langle f, \tilde\psi_\xi \rangle\}\|_{\tf^s_{pq}}
&= \Big\|\Big(\sum_{j\ge 0}
\sum_{\xi\in\cX_j}\big[|A_\xi|^{-s/d}|\langle f, \tilde\psi_\xi \rangle|\tONE_{A_\xi}(\cdot)\big]^q
\Big)^{1/q}\Big\|_{p}\\
&\le c\Big\|\Big(\sum_{j\ge 0}
\big[\cM_r\big(|B(\cdot, b^{-j})|^{-s/d}\Gamma_{\lambda_j}f\big)(\cdot)\big]^q
\Big)^{1/q}\Big\|_{p}\\
&\le c\Big\|\Big(\sum_{j\ge 0}
\big[|B(\cdot, b^{-j})|^{-s/d}|\Gamma_{\lambda_j}f(\cdot)|\big]^q
\Big)^{1/q}\Big\|_{p}
=c\|f\|_{\tF_{pq}^s(\Gamma)}.
\end{align*}
Hence
$
\|\{\langle f, \tilde\psi_\xi \rangle\}\|_{\tf^s_{pq}}
\le c \|f\|_{\tF^s_{pq}(\Gamma)} \le c \|f\|_{\tF^s_{pq}},
$
where we used that the functions $\Gamma_j$, $j\ge 0$, can be used to define
an equivalent norm in $F^s_{pq}$ (see Proposition~\ref{prop:Fspq-independ}).

Therefore, the operator $S_{\tilde\psi}:\tF_{pq}^s\rightarrow \tf_{pq}^s$ is bounded.

The identity $T_{\psi}\circ S_{\tilde\psi}=Id$ on $\tF_{pq}^s$ follows by
Proposition~\ref{prop:decomp-DD2} (c).
This completes the proof of the theorem.
$\qed$

\subsection{Identification of some Triebel-Lizorkin spaces}\label{sec:identify-Fspaces}

We next show that the Triebel-Lizorkin spaces can be viewed as a generalization of
a certain Sobolev type spaces and, in particular, of $\LL^p$, $1<p<\infty$.

In this part, we again make the additional assumption that the {\em reverse doubling condition}
(\ref{reverse-doubling}) is valid, yielding (\ref{D3}).

\smallskip

\noindent
{\bf Generalized Sobolev spaces.} Let $s\in \R$ and $1\le p\le \infty$.
The space $H_s^p$ is defined as the set of all $f\in \cD'$ such that
\begin{equation}\label{def-Hsp}
\|f\|_{H_s^p}:=\|(\Id+L)^{s/2}f\|_p < \infty.
\end{equation}


\begin{theorem}\label{thm:ident-F-spaces}
The following identification is valid:
$$
F^s_{p2}=H_s^p,\quad s\in\R, \quad 1<p<\infty,
$$
with equivalent norms, and in particular,
$$
F^0_{p2}=H_0^p=\LL^p,\quad 1<p<\infty.
$$
\end{theorem}

\noindent
{\bf \boldmath $L^p$-multipliers.}
To establish the above result we next develop $L^p$ multipliers.


\begin{theorem}\label{thm:multipliers}
Suppose $m\in C^k(\R_+)$ for some $k>d$, $m^{(2\nu+1)}(0)=0$ for $\nu\ge 0$ such that $2\nu+1\le k$,
and $\sup_{\lambda\in \R_+}|\lambda^\nu m^{(\nu)}(\lambda)| <\infty$, $0\le \nu\le k$.
Then the operator $m(\sqrt{L})$ is bounded on $\LL^p$ for $1<p<\infty$.
\end{theorem}

\noindent
{\bf Proof.}
As before choose $\ph_0\in C^{\infty}(\bR_+)$ so that $\supp \ph_0\subset [0, 2]$, $0\le \ph_0 \le 1$, and
$\ph_0(\lambda)=1$ for $\lambda \in [0, 1]$.
Let $\ph(\lambda):=\ph_0(\lambda)-\ph_0(2\lambda)$.

Set
$\ph_j(\lambda):=\ph(2^{-j}\lambda)$, $j\ge 1$.
Clearly,
$\sum_{j \ge 0}\ph_j(\lambda)=1$ for $\lambda \in \bR_+$
and hence
$m(\sqrt{L})=\sum_{j\ge 0}m(\sqrt{L})\varphi_j(\sqrt{L})$.
Set $\omega_j(\lambda):=m(2^j\lambda)\varphi(\lambda)$, $j\ge 1$,
and $\omega_0(\lambda):=m(\lambda)\varphi_0(\lambda)$.
Then $\omega_j(2^{-j}\sqrt{L})= m(\sqrt{L})\varphi_j(\sqrt{L})$, $j\ge 0$.
From the hypothesis of the theorem

it readily follows that
$\sup_{\lambda\in \R_+}|\omega_j^{(\nu)}(\lambda)| \le c_k <\infty$ for $0\le \nu\le k$.
Then by Theorem~\ref{thm:main-local-kernels} (see Remark~\ref{rem:kd})
\begin{equation}\label{multipl-1}
\big|\omega_j(2^{-j}\sqrt{L})(x, y)\big| \le c (|B(x, 2^{-j})||B(y, 2^{-j})|)^{-1/2}\big(1+2^j\rho(x, y)\big)^{-k}
\end{equation}
and whenever $\rho (y, y') \le 2^{-j}$
\begin{equation}\label{multipl-2}
\big|\omega_j(2^{-j}\sqrt{L})(x, y)  -  \omega_j(2^{-j}\sqrt{L})(x,y')\big|
\le \frac{c(2^j\rho(y,y'))^\alpha \big(1+2^j\rho(x, y)\big)^{-k}}
{(|B(x, 2^{-j})||B(y, 2^{-j})|)^{1/2}}.
\end{equation}
We choose $0<\eps \le \alpha$ so that $d+2\eps\le k$.

Denote briefly $m_j(x, y):=\omega_j(2^{-j}\sqrt{L})(x, y)$ and set
$K(x, y):=\sum_{j\ge 0}m_j(x, y)$.
We shall show that $K(x, y)$ is well defined for $x\ne y$ and
$|K(x, y)| \le c|B(y, \rho(x, y))|^{-1}$,
and moreover  $K(x, y)$
obeys the following H\"{o}rmander condition
\begin{equation}\label{multipl-3}
\int_{M\setminus B(y, 2\delta)}|K(x, y)-K(x, y')|d\mu(x) \le c,
\quad \hbox{whenever}\quad y'\in B(y, \delta),
\end{equation}
for all $y\in M$ and $\delta>0$.
To this end it suffices to show that for some $\eps>0$ ($\eps$ from above will do)
\begin{equation}\label{multipl-4}
|K(x, y)-K(x, y')|
\le c\Big(\frac{\rho(y, y')}{\rho(x, y)}\Big)^\eps |B(y, \rho(x, y))|^{-1}
\end{equation}
whenever $\rho(y, y') \le \min\{\rho(x, y), \rho(x, y')\}$,
see \cite{CW}.

Given $x, y, y'\in M$ such that $0<\rho(y, y') \le \min\{\rho(x, y), \rho(x, y')\}$
we pick $\ell, n\in \bZ$ ($\ell \ge n$)
so that
$2^{-\ell-1}<\rho(y, y')\le 2^{-\ell}$ and
$2^{-n-1}<\rho(x, y)\le 2^{-n}$.
Without loss of generality we may assume that $n\ge 1$.
Then we can write
\begin{align*}
|K(x, y)-K(x, y')|
&\le \sum_{j=0}^n |m_j(x, y)-m_j(x, y')| + \sum_{j=n+1}^\ell \dots + \sum_{j\ge\ell+1} \dots \\
&= \Omega_1+\Omega_2+\Omega_3.
\end{align*}
To estimate $\Omega_1$ we note that by (\ref{D3})
$|B(y, 2^{-j})| \ge c(2^{j}\rho(x, y))^{-\bet}|B(y, \rho(x, y))|$ and
$$
|B(x, 2^{-j})| \ge c(2^{j+1}\rho(x, y))^{-\bet}|B(x, 2\rho(x, y))|
\ge c'(2^{j}\rho(x, y))^{-\bet}|B(y, \rho(x, y))|,
\; j\le n.
$$

Now, using (\ref{multipl-2}) we obtain
\begin{align}\label{est-Omega1}
\Omega_1
&\le c\sum_{j=0}^n \frac{(2^j\rho(y,y'))^\alpha}{(|B(x, 2^{-j})| |B(y, 2^{-j})|)^{1/2}}
\le \frac{c\rho(y,y')^\alpha}{|B(y, \rho(x, y))|}\sum_{j=0}^n 2^{j\alpha}(2^j\rho(x, y))^{\bet}\\
&\le c\Big(\frac{\rho(y, y')}{\rho(x, y)}\Big)^\alpha |B(y, \rho(x, y))|^{-1},
\quad (\rho(x, y)\sim 2^{-n}).
\notag
\end{align}
%
%
From (\ref{D2}) it follows that
$|B(y, \rho(x, y))| \le c(1+2^j\rho(x, y))^d|B(y, 2^{-j})|$ and
$$
|B(y, \rho(x, y))| \le |B(x, 2\rho(x, y))|
\le c(1+2^j\rho(x, y))^d|B(x, 2^{-j})|,
\; j\ge n+1.
$$
From these and (\ref{multipl-2}) we get
\begin{align}\label{est-Omega2}
\Omega_2
&\le c|B(y, \rho(x, y))|^{-1}
\sum_{j=n+1}^\ell \frac{(2^j\rho(y,y'))^\eps}{(1+2^j\rho(x, y))^{k-d}}\notag\\
&\le c\Big(\frac{\rho(y, y')}{\rho(x, y)}\Big)^\eps |B(y, \rho(x, y))|^{-1}
\sum_{j\ge n+1}\frac{1}{(1+2^j\rho(x, y))^{k-d-\eps}}\\
&\le c\Big(\frac{\rho(y, y')}{\rho(x, y)}\Big)^\eps |B(y, \rho(x, y))|^{-1}, \notag
\end{align}
where we used that $k-d-\eps\ge \eps>0$ and $\rho(x, y)\sim 2^{-n}$.
%
%
To estimate $\Omega_3$ we write
$$
\Omega_3 \le \sum_{j>\ell} |m_j(x, y)| + \sum_{j>\ell} |m_j(x, y')| =: \Omega_3'+\Omega_3''.
$$
By the above estimates for $|B(y, \rho(x, y))|$ and (\ref{multipl-1}) we get
\begin{align}\label{est-Omega31}
\Omega_3'
&\le c|B(y, \rho(x, y))|^{-1}
\sum_{j>\ell} \frac{1}{(1+2^j\rho(x, y))^{k-d}}
\le c\frac{|B(y, \rho(x, y))|^{-1}}{(2^\ell\rho(x, y))^{2\eps}}\\
&\le c\Big(\frac{\rho(y, y')}{\rho(x, y)}\Big)^\eps |B(y, \rho(x, y))|^{-1}, \notag
\end{align}
where we used that $\rho(y, y')\sim 2^{-\ell}$.
One similarly obtains
\begin{align}\label{est-Omega32}
\Omega_3''
\le c\Big(\frac{\rho(y, y')}{\rho(x, y')}\Big)^\eps |B(y', \rho(x, y'))|^{-1}
\le c\Big(\frac{\rho(y, y')}{\rho(x, y)}\Big)^\eps |B(y, \rho(x, y))|^{-1}.
\end{align}
Here the last inequality follows by (\ref{D2}) using that
$\rho(x, y') \sim \rho(x, y)$, which follows from the condition
$\rho(y, y') \le \min\{\rho(x, y), \rho(x, y')\}$.
Putting together estimates (\ref{est-Omega1})-(\ref{est-Omega32}) we obtain (\ref{multipl-4}).
Therefore, the kernel $K(\cdot, \cdot)$ satisfies the H\"{o}rmander condition (\ref{multipl-3}).

The estimate
$|K(x, y)| \le \sum_{j\ge 0}|m_j(x, y)| \le c|B(y, \rho(x, y))|^{-1}$, $x\ne y$,
follows from (\ref{multipl-1}) similarly as above.

We next show that for any compactly supported function $f\in \LL^\infty$
\begin{equation}\label{ker-m}
m(\sqrt{L})f(x)=\int_MK(x, y)f(y)d\mu(y)
\quad\hbox{ for almost all $x \not\in \supp f$.}
\end{equation}
This and the fact that the kernel $K(\cdot, \cdot)$ satisfies the H\"{o}mander condition (\ref{multipl-3})
and $\|m(\sqrt{L})\|_{2\to 2} <\infty$ entails that $m(\sqrt{L})$ is a generalized Calder\'{o}n-Zygmund operator
and therefore $m(\sqrt{L})$ is bounded on $\LL^p$, $1<p<\infty$ (see \cite{CW}).

In turn, identity (\ref{ker-m}) readily follows from this assertion:
If $f_1, f_2\in L^\infty$ are compactly supported and
$\rho(\supp f_1, \supp f_2)\ge c >0$, then
\begin{align}\label{ker-mM}
\langle m(\sqrt{L}) f_1, f_2\rangle
&=\lim_{N\to \infty}
\int_M\int_M \sum_{j=0}^N m_j(x, y)f_1(y)\overline{f_2(x)} d\mu(y) d\mu(x)\\
&=
\int_M\int_M K(x, y)f_1(y)\overline{f_2(x)} d\mu(y) d\mu(x). \notag
\end{align}
The left-hand side identity in (\ref{ker-mM}) is the same as
$$
\langle m(\sqrt{L}) f_1, f_2\rangle
= \lim_{N\to \infty}\sum_{j=0}^N \langle m(\sqrt{L})\varphi_j(\sqrt{L}) f_1, f_2\rangle,
$$
which follows from the fact that
$m(\sqrt{L})f=\sum_{j\ge 0}m(\sqrt{L})\varphi_j(\sqrt{L})f$ in $\LL^2$ for each $f\in \LL^2$
by the spectral theorem.
The right-hand side identity in (\ref{ker-mM}) follows by
$K(x, y) = \sum_{j\ge 0}m_j(x, y)$ and
$\sum_{j\ge 0}|m_j(x, y)| \le c|B(y, \rho(x, y))|^{-1}$ for $x\ne y$,
applying the Lebesgue dominated convergence theorem.

To derive (\ref{ker-m}) from (\ref{ker-mM}) one argues as follows:
Given $f\in L^\infty$ with compact support and $x\notin \supp f$,
one applies (\ref{ker-mM}) with $f_1:=f$ and $f_2:=|B(x, \delta)|^{-1}\ONE_{B(x, \delta)}$,
where $\delta < \rho(x, \supp f)$.
Then passing to the limit as $\delta \to 0$ one arrives at (\ref{ker-m}).
The proof is complete.
$\qed$

\medskip

\noindent
{\bf Proof of Theorem~\ref{thm:ident-F-spaces}.}
Assume first that $f\in H^p_s$, $s\in \R$, $1<p<\infty$.
Let the functions $\varphi_j\in C^\infty_0(\R_+)$, $j=0, 1, \dots$,
be as in the definition of Triebel-Lizorkin and Besov spaces
with this additional property:
$\sum_{j\ge 0}\varphi_j(\lambda)=1$ for $\lambda\in \R_+$.
Assuming that $\eps:=\{\eps_j\}_{j\ge 0}$ is an arbitrary sequence with $\eps_j=\pm 1$,
we write
$$
T_\eps f:=\sum_{j\ge 0}\eps_j2^{js}\varphi_j(\sqrt L)f
= \sum_{j\ge 0} \omega_j(\sqrt L)(\Id+L)^{s/2}f
= m(\sqrt L)(\Id+L)^{s/2}f,
$$
where
$\omega_j(\lambda):= \eps_j 2^{js}(1+\lambda^2)^{-s/2}\varphi_j(\lambda)$
and
$m(\lambda)=\sum_{j\ge 0} \omega_j(\lambda)$.
Using that $\varphi_j(\lambda)=\varphi(2^{-j}\lambda)$, $j\ge 1$,
with $\varphi\in C^\infty$ and $\supp \varphi\subset [1/2, 2]$
it is easy to see that
$\sup_{\lambda>0}|\lambda^\nu\omega_j^{(\nu)}(\lambda)| \le c_\nu$, $\nu\ge 0$,
with $c_\nu$ a constant independent of $j$ and since 
$\supp \varphi_0\subset [0, 2]$ and
$\supp \varphi_j\subset [2^{j-1}, 2^{j+1}]$, $j\ge 1$, then
$\sup_{\lambda>0}|\lambda^\nu m^{(\nu)}(\lambda)| \le 2c_\nu$.
We now appeal to Theorem~\ref{thm:multipliers} to obtain
$\|T_\eps f\|_p \le c \|(\Id+L)^{s/2}f\|_p$, $1<p<\infty$,
for any sequence $\eps:=\{\eps_j\}_{j\ge 0}=\{\pm 1\}$.
Finally, applying Khintchine's inequality (which involve the Rademacher functions) as usual
we arrive at
$$
\|f\|_{F^{s}_{p2}}
\le c\Big\|\Big(\sum_{j\ge 0} \Big(
2^{js}
|\varphi_j(\sqrt L) f(\cdot)|
\Big)^2\Big)^{1/2}\Big\|_{p}
\le c\|(\Id+L)^{s/2}f\|_p
= c \|f\|_{H^p_s}.
$$

To prove an estimate in the opposite direction, let $f\in F^{s}_{p2}$, $s\in \R$, $1<p<\infty$.
We now assume that $\varphi_j\in C^\infty_0(\R_+)$, $j=0, 1, \dots$,
are as in the definition of Tribel-Lizorkin spaces but
with this additional property:
$\sum_{j\ge 0}\varphi_j^2(\lambda)=1$ for $\lambda\in \R_+$.
Using this we can write
$$
(\Id+L)^{s/2}f
= \sum_{j\ge 0}2^{-js}(\Id+L)^{s/2}\varphi_j(\sqrt L)2^{js}\varphi_j(\sqrt L)f
= \sum_{j\ge 0}\theta_j(\sqrt L)2^{js}\varphi_j(\sqrt L)f,
$$
where $\theta_j(\lambda):= 2^{-js}(1+\lambda^2)^{s/2}\varphi_j(\lambda)$.
Denote
$\bZ_r^+:= \{2k+r: k=0, 1, \dots\}$, $r=0, 1$,
and set
$G_r f:= \sum_{j \in \bZ_r^+}\theta_j(\sqrt L)2^{js}\varphi_j(\sqrt L)f$.
Evidently,
$(\Id+L)^{s/2}f = G_0 f+G_1 f$.
Let $\{\eps_{jr}\}_{j\in \bZ_r}$ be an arbitrary sequence with $\eps_{jr}=\pm 1$.
The supports of $\theta_j$ and $\varphi_k$ do not overlap if $j, k\in \bZ_r$, $j\ne k$, and hence
$\theta_j(\sqrt L)\varphi_k(\sqrt L)\equiv 0$ if $j, k\in \bZ_r$, $j\ne k$.
Therefore,
$$
G_r f:=\sum_{j \in \bZ_r^+}\eps_{jr}\theta_j(\sqrt L)
\sum_{j \in \bZ_r^+}\eps_{jr}2^{js}\varphi_j(\sqrt L)f
=m_r(\sqrt L) \sum_{j \in \bZ_r^+}\eps_{jr}2^{js}\varphi_j(\sqrt L)f,
$$
where
$m_r(\lambda):=\sum_{j \in \bZ_r^+}\eps_{jr}\theta_j(\lambda)$.
As above we have
$\sup_{\lambda>0}|\lambda^\nu\theta_j^{(\nu)}(\lambda)| \le c_\nu$, $\nu\ge 0$,
with $c_\nu$ independent of $j$ and hence
$\sup_{\lambda>0}|\lambda^\nu m_r^{(\nu)}(\lambda)| \le c_\nu$, $\nu\ge 0$.
Applying Theorem~\ref{thm:multipliers} we get for any sequence
$\{\eps_{jr}\}_{j\in \bZ_r}=\{\pm 1\}$
$$
\|G_r f\|_p \le c\Big\|\sum_{j \in \bZ_r^+}\eps_{jr}2^{js}\varphi_j(\sqrt L)f\Big\|_p,
\quad 1<p<\infty.
$$
An application of Khintchine's inequality gives
$$
\|G_r f\|_p \le c\Big\|\Big(\sum_{j \in \bZ_r^+} \Big(
2^{js}
|\varphi_j(\sqrt L) f(\cdot)|
\Big)^2\Big)^{1/2}\Big\|_{p}
\le c\|f\|_{F^{s}_{p2}},
\quad r=0, 1,
$$
which implies
$\|(\Id+L)^{s/2}f\|_p \le \|G_0f\|_p + \|G_1 f\|_p \le c\|f\|_{F^{s}_{p2}}$.
$\qed$

\end{document}